\documentclass[leqno,12pt]{article}
\usepackage{amsmath,amssymb,rotating,a4wide,color}
\usepackage{wasysym}
\usepackage{hyperref}
\usepackage[all,cmtip]{xy}

\usepackage{verbatim}

\newdir^{ (}{{}*!/-3pt/\dir^{(}}
\newdir_{ (}{{}*!/-3pt/\dir_{(}}

\entrymodifiers={+!!<0pt,\fontdimen22\textfont2>}

\def\bigtimes{\mathop{\raise-2pt\hbox{\huge$\times$}}}

\newbox\circbulletbox
\setbox\circbulletbox=\hbox{\,{\large$\circledcirc$}\kern-8.7pt\raise .6pt\hbox{$\bullet$}\,}

\let\le\leqslant
\let\ge\geqslant

\def\circVbig{\hbox{\text{\it\r{V}}}}
\def\circVscript{\hbox{\scriptsize\text{\it\r{V}}}}
\def\circVscriptscript{\mbox{\tiny\text{\it\r{V}}}}
\def\circV{{\mathchoice{\circVbig}{\circVbig}{\circVscript}{\circVscriptscript}{}}}
\def\circVlimits_#1^#2{{\mathchoice%
   {\circVbig{}^{\kern2pt #2}_{\kern-2pt #1}}%
   {\circVbig{}^{\kern2pt #2}_{\kern-2pt #1}}%
   {\scriptstyle\circVscript{}^{\kern1.7pt #2}_{\kern-1pt #1}}%
   {\scriptscriptstyle\circVscriptscript{}^{\kern1.5pt #2}_{\kern-1pt #1}}%
   }}
\def\circVr_#1{\circVlimits_#1^r}
\def\circVs_#1{\circVlimits_#1^s}
\def\circVprime{{\circVlimits_{}^{\prime}}}
\def\circVpprime{{\circVlimits_{}^{\prime\prime}}}

\def\circWbig{\hbox{\text{\it\r{W}}}}
\def\circWscript{\hbox{\scriptsize\text{\it\r{W}}}}
\def\circWscriptscript{\mbox{\tiny\text{\it\r{W}}}}
\def\circW{{\mathchoice{\circWbig}{\circWbig}{\circWscript}{\circWscriptscript}}}
\def\circWlimits_#1^#2{{\mathchoice%
   {\circWbig{}^{\kern2pt #2}_{\kern-2pt #1}}%
   {\circWbig{}^{\kern2pt #2}_{\kern-2pt #1}}%
   {\scriptstyle\circWscript{}^{\kern1.7pt #2}_{\kern-1pt #1}}%
   {\scriptscriptstyle\circWscriptscript{}^{\kern1.5pt #2}_{\kern-1pt #1}}%
   }}

\def\OM{\mathchoice
  {\rlap{\kern3.2pt$\overline{\phantom{L}}$}M}
  {\rlap{\kern3.2pt$\overline{\phantom{L}}$}M}
  {\rlap{\kern2.4pt$\scriptstyle\overline{\phantom{L}}$}M}
  {\rlap{\kern1.8pt$\scriptscriptstyle\overline{\phantom{L}}$}M}}

\def\mycirc{{\kern1pt\circ\kern2pt}}

\def\Cinf{{{\BC}_\infty}}

\def\norm{{\rm norm}}
\def\Eis{\mathop{\rm Eis}\nolimits}
\def\Mod{\mathop{\rm Mod}\nolimits}
\def\Rel{\mathop{\rm Rel}\nolimits}
\def\Div{\mathop{\rm Div}\nolimits}

\def\Aut{\mathop{\rm Aut}\nolimits}

\def\End{\mathop{\rm End}\nolimits}

\def\Proj{\mathop{\rm Proj}\nolimits}
\def\Spec{\mathop{\rm Spec}\nolimits}

\def\deg{\mathop{\rm deg}\nolimits}

\def\mymod{\mathop{\rm mod}\nolimits}
\def\Ker{\mathop{\rm Ker}\nolimits}

\def\Quot{\mathop{\rm Quot}\nolimits}

\def\GL{\mathop{\rm GL}\nolimits}
\def\SL{\mathop{\rm SL}\nolimits}

\def\red{{\rm red}}

\def\univ{{\rm univ}}
\def\id{{\rm id}}

\let\phi\varphi
\let\theta\vartheta
\let\epsilon\varepsilon
\let\setminus\smallsetminus

\let\emptyset\varnothing

\newtheorem{Thm}{Theorem}[section]
\newtheorem{Prop}[Thm]{Proposition}
\newtheorem{Lem}[Thm]{Lemma}

\newtheorem{Cor}[Thm]{Corollary}

\newtheorem{Def}[Thm]{Definition}

\newtheorem{PropDef}[Thm]{Proposition-Definition}

\newtheorem{Rem}[Thm]{Remark}

\newtheorem{Ass}[Thm]{Assumption}

\newtheorem{Cons}[Thm]{Construction}

\newtheorem{Exp}[Thm]{Explanation}

\numberwithin{Thm}{subsection}
\def\UseTheoremCounterForNextEquation{\setcounter{equation}{\value{Thm}}\addtocounter{Thm}{1}}

\def\qed{{\hskip0pt\unskip\unskip\nobreak\hfil\penalty50
          \hskip1em\hbox{}\nobreak\hfil
           {$\square$}
          \parfillskip=0pt\finalhyphendemerits=0
          \par}\medskip}

\newenvironment{Proof}
               {\noindent{\bf Proof.}\ }
               {\qed}

               {\noindent{\bf Proof of #1.}\ }
               {\qed}

\newcommand{\BC}{{\mathbb{C}}}

\newcommand{\BF}{{\mathbb{F}}}
\newcommand{\BG}{{\mathbb{G}}}

\newcommand{\BP}{{\mathbb{P}}}

\newcommand{\BZ}{{\mathbb{Z}}}

\newcommand{\Fa}{{\mathfrak{a}}}
\newcommand{\Fb}{{\mathfrak{b}}}

\newcommand{\Fp}{{\mathfrak{p}}}

\newcommand{\CF}{{\cal F}}

\newcommand{\CI}{{\cal I}}

\newcommand{\CL}{{\cal L}}
\newcommand{\CM}{{\cal M}}

\newcommand{\CO}{{\cal O}}

\newcommand{\CS}{{\cal S}}

\newbox\mybox
\def\arrover#1{\mathrel{
       \setbox\mybox=\hbox spread 1.4em
              {\hfil$\scriptstyle#1$\hfil}
       \vbox{\offinterlineskip\copy\mybox
             \hbox to\wd\mybox{\rightarrowfill}}}}

\def\larrover#1{\mathrel{
       \setbox\mybox=\hbox spread 1.4em
              {\hfil$\scriptstyle#1\vphantom{g}$\hfil}
       \vbox{\offinterlineskip\copy\mybox
             \hbox to\wd\mybox{\leftarrowfill}}}}

\def\ontoover#1{\mathrel{
       \setbox\mybox=\hbox spread 1.4em
              {\hfil$\scriptstyle#1\vphantom{g}$\hfil}
       \vbox{\offinterlineskip\copy\mybox
             \hbox to\wd\mybox{\rightarrowfill\hskip-2.8mm
                               $\rightarrow$}}}}
\def\leftontoover#1{\mathrel{
       \setbox\mybox=\hbox spread 1.4em
              {\hfil$\scriptstyle#1\vphantom{g}$\hfil}
       \vbox{\offinterlineskip\copy\mybox
             \hbox to\wd\mybox{$\leftarrow$\hskip-2.8mm
                               \leftarrowfill}}}}
\let\longto\longrightarrow
\let\into\hookrightarrow
\let\onto\twoheadrightarrow

\def\longinto{\lhook\joinrel\longrightarrow}

\def\isoto{\mathrel{
       \setbox\mybox=\hbox spread 0.9em
              {\hfil$\scriptstyle\sim$\hfil}
       \vbox{\offinterlineskip\copy\mybox
             \hbox to\wd\mybox{\rightarrowfill}}}}

\def\Bigskip{\bigskip\bigskip}

\begin{document}

\title{\strut
\vskip-80pt
Compactification of Drinfeld Moduli Spaces\\
as Moduli Spaces of $A$-Reciprocal Maps\\
and Consequences for Drinfeld Modular Forms}
\author{Richard Pink\\[12pt]
\small Department of Mathematics \\[-3pt]
\small ETH Z\"urich\\[-3pt]
\small 8092 Z\"urich\\[-3pt]
\small Switzerland \\[-3pt]
\small pink@math.ethz.ch\\[12pt]}

%\date{March 17, 2016}

\maketitle

\medskip
\begin{center}
\large In memory of David Goss
\end{center}

\Bigskip

\begin{abstract}
We construct a compactification of the moduli space of Drinfeld modules of rank $r$ and level $N$ as a moduli space of $A$-reciprocal maps. This is closely related to the Satake compactification, but not exactly the same. The construction involves some technical assumptions on $N$ that are satisfied for a cofinal set of ideals~$N$. 
In the special case $A=\BF_q[t]$ and $N=(t^n)$ we obtain a presentation for the graded ideal of Drinfeld cusp forms of level $N$ and all weights and can deduce a dimension formula for the space of cusp forms of any weight. We expect the same results in general, but the proof will require more ideas.
\end{abstract}

{\renewcommand{\thefootnote}{}
%11F52 Modular forms associated to Drinfeld modules
%11G09 Drinfeld modules; higher-dimensional motives, etc.
%14D20 Algebraic moduli problems, moduli of vector bundles
%14D22 Fine and coarse moduli spaces
%14M27 Compactifications; symmetric and spherical varieties
}

\newpage
\tableofcontents
\newpage

%%%%%%%%%%%%%%%%%%%%%%%%%%%%%%%%%%%%%%%%%%%%%%%%%%%%%%%%%%%%%%%%%%%%%%%%
%%%%%%%%%%%%%%%%%%%%%%%%%%%%%%%%%%%%%%%%%%%%%%%%%%%%%%%%%%%%%%%%%%%%%%%%

\section*{Introduction}
\addcontentsline{toc}{section}{Introduction}
\label{Intro}

Consider an admissible coefficient ring $A$ over $\BF_q$ with field of quotients~$F$ and a non-zero proper ideal $N$ of~$A$. Let $M^r_{A,N}$ denote the fine moduli space of Drinfeld $A$-modules of rank $r\ge1$ in generic characteristic with a full level $N$-structure. This is an irreducible smooth affine algebraic variety of finite type and dimension $r-1$ over~$F$. Let $\OM^r_{A,N}$ be its Satake compactification according to \cite{PinkSatake}.

\medskip
{\bf $\BF_q$-reciprocal maps:} 
First consider the special case $A=\BF_q[t]$ and $N=(t)$. Here $\smash{M^r_{A,N}}$ is the base change from $\BF_q$ of the open subscheme $\Omega_V$ of $\smash{\BP^{r-1}_{\BF_q}}$ obtained by removing all $\BF_q$-rational hyperplanes. 
Moreover $\smash{\OM^r_{A,N}}$ is the base change from $\BF_q$ of a certain compactification $Q_V$ of~$\Omega_V$ constructed in \cite{PinkSchieder}. (For $r\ge3$ this $Q_V$ is not isomorphic to the tautological compactification $\BP^{r-1}_{\BF_q}$.) The construction comprises an explicit presentation of the projective coordinate ring $R_V$ of~$Q_V$ and hence of the ring of Drinfeld modular forms of level $N$ and all weights.
Here $Q_V$ is obtained by giving a simple construction of~$R_V$, but it is also a fine moduli scheme, as follows.

Let $V$ be an $\BF_q$-vector space of dimension $r$ and set $\circV := V\setminus\{0\}$. For any $\BF_q$-algebra $R$ a \emph{fiberwise invertible $\BF_q$-reciprocal map} $\rho\colon \smash{\circV}\to R^\times$ is any map of the form $v\mapsto \lambda(v)^{-1}$ for an $\BF_q$-linear map $\lambda\colon V\to R$ satisfying $\lambda(\smash{\circV}\;)\subset R^\times$. Any such map is characterized by the equations:
\begin{itemize}
\item[(a)] $\rho(v) \cdot \rho(w) = \rho(v+w) \cdot (\rho(v) + \rho(w))$ for all $v$, $w\in \circV$ with $v+w\in\circV$, and
\item[(b)] $\alpha\rho(\alpha v) = \rho(v)$ for all $\alpha \in \BF_q^{\times}$ and $v \in \circV$.
\end{itemize}
A general \emph{$\BF_q$-reciprocal map} $\rho\colon \circV\to R$ is defined simply as any map with values in $R$ which satisfies the same equations.
%See Definition \ref{FqRecipDef}. 
This notion is globalized to maps $\circV\to\CL(S)$ with values in an invertible sheaf $\CL$ on a scheme $S$ over~$\BF_q$ in Definition \ref{GenFqRecipDef}. Then $Q_V$ becomes a fine moduli scheme for isomorphism classes of pairs $(\CL,\rho)$ consisting of an invertible sheaf $\CL$ and a fiberwise non-zero $\BF_q$-reciprocal map~$\rho$. 

%%%%%%%%%%%%%%%%%%%%%%%%%%%%%%%%%%%%%

\medskip
{\bf $A$-reciprocal maps:} 
The goal of the article at hand is to generalize this theory to produce a compactification of $M^r_{A,N}$ for general $A$ and~$N$. For this we consider the finite $A$-module $V^r_N := (N^{-1}/A)^{\oplus r}$. To any Drinfeld $A$-module $\phi\colon A\to R[\tau]$ of rank $r$ and any level $N$-structure $\lambda\colon {V^r_N\isoto\phi[N]}$ we associate the fiberwise invertible $\BF_q$-reciprocal map $\rho\colon {\circVr_N\to R^\times}$, ${v\mapsto \lambda(v)^{-1}}$. Our first job was to find a useful additional identity satisfied by $\rho$ which reflects the action of $A$ on~$V^r_N$. This problem had a surprisingly simple solution. Namely, consider the set $\Div(N) := {\{a\in A\mid N\subset(a)\}}$ of divisors of~$N$, and for any $a\in\Div(N)$ consider the $A$-submodule $V^r_a := (Aa^{-1}/A)^{\oplus r} \subset V^r_N$. Then by Proposition \ref{DrinModConv} we have the additional identity:
%(see Definition \ref{ARecipDef})
\begin{itemize}
\item[(b$'$)] $a\rho(av) = \sum_{v'\in V^r_a}\rho(v-v')$ for all $a\in\Div(N)$ and  $v \in V^r_N\setminus V^r_a$.
\end{itemize}
We therefore define a general \emph{$A$-reciprocal map} $\rho\colon \circVr_N\to R$ as any map with values in $R$ that satisfies the conditions (a) and (b$'$).

Then we must show that any fiberwise invertible $A$-reciprocal map ${\circVr_N\to R^\times}$ arises from a unique pair $(\phi,\lambda)$ as above. We achieve this in Proposition \ref{DrinMod} under certain technical conditions on the level $N$ that are collated in Assumption \ref{NAss}. The main requirement is that $\Div(N)$ generates $A$ as an $\BF_q$-algebra, while the other assumptions appear for technical reasons and can perhaps be discarded. By Proposition \ref{NExists} the assumptions are satisfied for a cofinal set of non-zero ideals~$N$. For the following we assume that they hold for~$N$.

Next the notion of $A$-reciprocal maps is globalized to maps $\circVr_N\to\CL(S)$ with values in an invertible sheaf $\CL$ on a scheme $S$ over~$\BF_q$ in Definition \ref{GenARecipDef}. By standard arguments there is a fine moduli scheme $Q_{A,V^r_N}$ of isomorphism classes of pairs $(\CL,\rho)$ consisting of an invertible sheaf $\CL$ and a fiberwise non-zero $A$-reciprocal map~$\rho$, and $Q_{A,V^r_N}$ is projective over~$F$. It also contains an open subscheme $\Omega_{A,V^r_N}$ which is a fine moduli scheme of isomorphism classes of fiberwise invertible $A$-reciprocal maps and therefore naturally isomorphic to $M^r_{A,N}$. Thus $Q_{A,V^r_N}$ constitutes a natural compactification of $M^r_{A,N}$ as a fine moduli scheme.

%%%%%%%%%%%%%%%%%%%%%%%%%%%%%%%%%%%%%

\medskip
{\bf Relation with the Satake compactification:} 
We show that $Q_{A,V^r_N}$ shares many properties with the Satake compactification $\OM^r_{A,N}$. For instance, in Theorem \ref{Dense} we prove that the open subscheme $M^r_{A,N} \cong \Omega_{A,V^r_N}$ is dense in $Q_{A,V^r_N}$. In Theorem \ref{Stratification} we show that $Q_{A,V^r_N}$ is stratified by finitely many locally closed subschemes $\Omega_W$ which are indexed by all non-zero free $A/N$-submodules $W\subset V^r_N$ and are isomorphic to $M^s_{A,N}$ for $1\le s\le r$. 
In Theorem \ref{SatakeProjection} we prove that the Satake compactification $\smash{\OM^r_{A,N}}$ is the normalization of $Q_{A,V^r_N}$ in the function field of $M^r_{A,N}$, and in Proposition \ref{RAVNormalMap} we show that the natural morphism $\pi\colon \smash{\OM^r_{A,N}} \to Q_{A,V^r_N}$ is finite and surjective. 

However, a computation of H\"aberli \cite[Prop.\;7.13, Cor.\;7.28]{Haeberli} implies that in general distinct points from the Satake compactification are identified with each other in $Q_{A,V^r_N}$. 
%(See Remark \ref{NormNotIsomRem}.)
Nevertheless, I expect that this is the only difference. More precisely, following H\"aberli it will be possible to say precisely which points are identified in $Q_{A,V^r_N}$. My fond hope is then that $Q_{A,V^r_N}$ is simply the quotient of $\smash{\OM^r_{A,N}}$ by the resulting equivalence relation on the underlying topological space.
%Equivalent: That the structure sheaf $\CO_{Q_{A,V^r_N}}$ is the subsheaf of $\pi_*\CO_{\OM^r_{A,N}}$ consisting of all local functions which are constant on the reduced subscheme of every fiber. 

%%%%%%%%%%%%%%%%%%%%%%%%%%%%%%%%%%%%%

\medskip
{\bf The projective coordinate ring:} 
Most of our constructions are done in the projective coordinate ring $R_{A,V^r_N}$ underlying $\smash{Q_{A,V^r_N}}$. This ring is given by an explicit presentation in Construction \ref{RAVCons}. The open subscheme $M^r_{A,N} \cong \Omega_{A,V^r_N}$ corresponds to a certain localization $RS_{A,V^r_N}$ of $R_{A,V^r_N}$ which is a regular graded integral domain. 
We expect that $R_{A,V^r_N}$ is itself an integral domain and that the natural homomorphism $R_{A,V^r_N} \to RS_{A,V^r_N}$ is injective, but are not yet able to prove this in general.

But let $R_{A,V^r_N}^\norm$ denote the integral closure of $R_{A,V^r_N}$ in $RS_{A,V^r_N}$. By the above-mentioned result on the Satake compactification this is the projective coordinate ring of $\smash{\OM^r_{A,N}}$. Also, let $I_{A,V^r_N}^\norm \subset R_{A,V^r_N}^\norm$ denote the graded ideal of the reduced boundary $(\smash{\OM^r_{A,N}}\setminus M^r_{A,N})^\red$. 
Then in Theorem \ref{ModCuspFormsIsom} we deduce that $R_{A,V^r_N}^\norm$ is the ring of Drinfeld modular forms of level $N$ and all weights and that $I_{A,V^r_N}^\norm$ is the ideal of all cusp forms therein.

Moreover, in (\ref{AIVNDef}) we construct a certain reduced ideal $I_{A,V^r_N} \subset R_{A,V^r_N}$ such that $I_{A,V^r_N}^\norm$ is the radical of the associated ideal $I_{A,V^r_N}\cdot R_{A,V^r_N}^\norm$. The fond hope expressed above corresponds to the expectation that the 
natural map $I_{A,V^r_N}\to I_{A,V^r_N}^\norm$ is an isomorphism.
Given the explicit presentation of the ring $I_{A,V^r_N}$ this would provide an explicit presentation of the ideal of all cusp forms. As a consequence this might lead to a dimension formula for spaces of cusp forms.

%%%%%%%%%%%%%%%%%%%%%%%%%%%%%%%%%%%%%

\medskip
{\bf A special case:} 
I would not bother writing all this up without more positive results in some new cases. Assume that $A=\BF_q[t]$ and $N=(t^n)$ for some $n\ge1$. In Theorem \ref{InjCor} we then prove that $R_{A,V^r_N}$ is an integral domain and injects into $RS_{A,V^r_N}$.
In Theorem \ref{RnCM} we show that $\smash{R_{A,V^r_N}}$ and hence $\smash{Q_{A,V^r_N}}$ is Cohen-Macaulay. It is therefore reasonable to expect that $R_{A,V^r_N}^\norm$ and $\OM^r_{A,N}$ are Cohen-Macaulay as well and that the same holds for general $A$ and~$N$.
In Theorem \ref{InInNorm} we prove that the natural map $I_{A,V^r_N}\to I_{A,V^r_N}^\norm$ is an isomorphism in this special case.
In Subsection \ref{Specialn5} we deduce a simple dimension formula for the space of Drinfeld cusp forms associated to any arithmetic subgroup $\Gamma<\SL_r(\BF_q[t])$ satisfying $\Gamma(t) < \Gamma < \Gamma_1(t)$.

The methods to attain these results are partly a refinement
% Weiterentwicklung
of methods from \cite{PinkSchieder}. There we had already considered a maximal unipotent subgroup $U<\Aut_{\BF_q}(V)$ and shown that the ring of $U$-invariants $R_V^U$ is isomorphic to a polynomial ring in $\dim_{\BF_q}(V)$ variables over~$\BF_q$ and that $R_V$ is a free $R_V^U$-module with an explicit basis.
In our special case we again consider a maximal subgroup $U<\GL_r(\BF_q[t]/(t^n))$ of $q$-power order, prove that the ring of invariants $\smash{R_{A,V^r_N}^U}$ is isomorphic to a polynomial ring in $r$ variables over~$F$, and show that $R_{A,V^r_N}$ is a free module over $R_{A,V^r_N}^U$ with an explicit basis. 
The method works, because $U$ is a group of $q$-power order acting on an $\BF_q$-vector space, because we can compute its invariants in $RS_{A,V^r_N}$, and because the respective module that we wish to describe happens to be a free module over the group ring $F[U]$. 
Unfortunately this method only succeeds in the case $N=(t^n)$, and proving similar results in the general case will require additional ideas.

%%%%%%%%%%%%%%%%%%%%%%%%%%%%%%%%%%%%%

\medskip
{\bf Outlook:} 
In addition to the expectations mentioned above one can ask which form a dimension formula for Drinfeld cusp forms might take in general. Based on the results from Theorem \ref{DimForm5} in our special case, we can surmise that for any fine congruence subgroup $\Gamma<\SL_r(A)$, the space $S_d(\Gamma)$ of cusp forms of weight $d$ and level $\Gamma$ has dimension 
$$\textstyle c(A,r) \cdot [\SL_r(A):\Gamma]\cdot\binom{d-1}{r-1},$$
where the constant $c(A,r)$ depends only on $A$ and~$r$. 
This constant might involve the class number of~$A$ and/or be related to a version of the Tamagawa number of $\SL_{r,A}$ or $\GL_{r,A}$, as Gekeler suggests.
Moreover, for any two fine congruence subgroups $\Gamma\triangleleft\Gamma'<\SL_r(A)$, the space $S_d(\Gamma)$ should be a free module over the group ring $\BF_q[\Gamma'/\Gamma]$.

In another direction one may ask for a dimension formula for modular forms instead of cusp forms. In \cite[Thm.\;4.1]{PinkSchieder} and \cite[Thm.\;8.4]{PinkSatake} we already gave such a formula in the case $A=\BF_q[t]$ for any subgroup $\Gamma$ satisfying $\Gamma(t)<\Gamma<\Gamma_1(t)$. An explicit presentation of the ring $R_{A,V^r_N}^\norm$ would probably yield a dimension formula in general.
%Theorem \ref{DimForm1} already gives a dimension formula for the homogeneous pieces of $R_{A,V^r_N}$ in our special case, but due to the discrepancy to $R_{A,V^r_N}^\norm$ this does not yet give the desired dimension formula.
Note that for rank $r=2$ and sufficiently large weight a dimension formula for modular forms was already given by Gekeler \cite[\S6]{GekelerBook}, and a formula for cusp forms can be obtained in the same way.

%%%%%%%%%%%%%%%%%%%%%%%%%%%%%%%%%%%%%

\medskip
{\bf Reflection on the notion of $A$-reciprocal maps:} 
The ad hoc definition of $A$-reciprocal maps and their study in H\"aberli's thesis \cite[\S8.2]{Haeberli} was an important encouragement for me. But I consider it as provisional and the new definition proposed in this paper as more useful.
One can check that our conditions imply his by Proposition \ref{RhoFormulas} (b) and Proposition \ref{RingHomo1} and that both definitions lead to moduli schemes with the same underlying reduced subscheme. But the conditions from \cite[Def.\;8.14]{Haeberli} are homogeneous of high degree and will therefore introduce an excess of nilpotent elements in the local rings at the boundary.
By contrast, the new relations (b$'$) above are homogeneous of degree $1$ and cannot be outdone in regard to their degree or their elegance.
%The new definition of $A$-reciprocal maps was not possible without the crucial formulas for $\BF_q$-reciprocal maps that I found in Subsection \ref{NiceForm} and which were not available to H\"aberli or Gekeler.

%I feel that my definition of $A$-reciprocal maps is reasonably elegant, if a little strange. 
Nevertheless, I am not yet sure that the new definition is quite final. It is still open which assumptions on the level $N$ are really necessary, and perhaps the definition should be augmented in order to reduce them.

One should also ask whether a variant of the definition of $A$-reciprocal maps might yield a fine moduli scheme that is \emph{isomorphic} to the Satake compactification $\OM^r_{A,N}$. 
With the results on the ideal of the boundary $I_{A,V^r_N}$ in our special case we seem to be almost there, and we would need just a little more data to distinguish different points at the boundary. Can one discover another property of reciprocal maps that helps to achieve this?
%One might want to incorporate this property into the definition.

%%%%%%%%%%%%%%%%%%%%%%%%%%%%%%%%%%%%%

\medskip
{\bf Relation with other work:} 
In a recent manuscript Gekeler \cite{GekelerHigherRankIV} pursues similar goals with a different approach. For simplicity he restricts himself to the case $A=\BF_q[t]$. Let $\Cinf$ denote the completion of the algebraic closure of the field $\BF_q((t^{-1}))$, let $\Omega^r$ be the Drinfeld period domain of rank $r$ over~$\Cinf$, and let $\Gamma(N) < \SL_r(\BF_q[t])$ be the principal congruence subgroup of level~$N$. Then $\Gamma(N)\backslash\Omega^r$ is one of finitely many irreducible components of $M^r_{A,N}(\Cinf)$. 
Let $\Mod(N)$ be the ring of analytic modular forms of level $N$ and all weights, so that $\Proj(\Mod(N))$ is the Satake compactification of $\Gamma(N)\backslash\Omega^r$, that is, the corresponding irreducible component of $\OM^r_{A,N}(\Cinf)$. 

Gekeler proposes to consider the $\Cinf$-subalgebra $\Eis(N)$ of $\Mod(N)$ that is generated by all Eisenstein series of weight $1$ and to view $\Proj(\Eis(N))$ as a natural compactification of $\Gamma(N)\backslash\Omega^r$. 
In \cite[Cor.\;7.6]{GekelerHigherRankIV} he proves that the natural morphism $\pi\colon \Proj(\Mod(N)) \to \Proj(\Eis(N))$ is bijective. He hopes that $\Mod(N)$ is equal to or at least very close to $\Eis(N)$, but unfortunately has no methods to decide that for rank $r>2$. 
%In rank $2$ there is a result of Cornelissen \cite[Thm.\;8.7]{GekelerHigherRankIV} that $\Mod(N)$ is generated by Eisenstein series of weight $1$ and all cusp forms of weight~$2$. 

%\medskip
To see the relation with our approach recall that $R_{A,V^r_N}^\norm$ is the ring of algebraic Drinfeld modular forms of level $N$ and all weights; hence $R_{A,V^r_N}^\norm \otimes_F\Cinf$ is isomorphic to a finite direct sum of copies of $\Mod(N)$ for all irreducible components of $M^r_{A,N}(\Cinf)$.
Also observe that the generators $[\frac{1}{v}\otimes1]$ of our ring $R_{A,V^r_N}$ from Construction \ref{RAVCons} represent the reciprocals of all non-zero $N$-torsion points of the universal Drinfeld module over $M^r_{A,N}$; hence they correspond to all Eisenstein series of weight $1$ for the group $\Gamma(N)$, for instance by \cite[(15.4)]{BBP3}. Thus the image of the homomorphism $R_{A,V^r_N}\otimes_F\Cinf \to R_{A,V^r_N}^\norm \otimes_F\Cinf$ followed by the projection to any one factor $\Mod(N)$ is precisely the subring $\Eis(N)$. 
In the special case $N=(t^n)$ one can hope that by combining our results on the ideal of the boundary $I_{A,V^r_N}$ with the bijectivity of the morphism $\Proj(\Mod(N)) \to \Proj(\Eis(N))$ one can deduce that $\Mod(N)=\Eis(N)$ in this case.
%: see Theorem \ref{GekelerPink}. 
%This gives further credence to the hopes of Gekeler \cite{GekelerHigherRankIV} and those of the present paper.
%Our methods from Section \ref{Specialn} are the only ones at present that are able to say something about the infinitesimal structure near the boundary.

%%%%%%%%%%%%%%%%%%%%%%%%%%%%%%%%%%%%%

\medskip
{\bf Structure of the paper:} 
The article is composed of three major sections. Section \ref{Fq} is devoted to $\BF_q$-reciprocal maps and can be viewed as a continuation of the article \cite{PinkSchieder}. For clarity I now call $\BF_q$-reciprocal maps what we simply called reciprocal maps in \cite{PinkSchieder}. 
We cover some additional topics with applications to $A$-reciprocal maps. In Subsection \ref{Functoriality} we discuss the functoriality of $\BF_q$-reciprocal maps under homomorphisms of finite dimensional $\BF_q$-vector spaces.
In Subsection \ref{NiceForm} we deduce some nice formulas for $\BF_q$-reciprocal maps, which motivated the above condition (b$'$) for $A$-reciprocal maps, and which are crucial for everything that follows.
In Subsections \ref{FqDescRing} and \ref{FqComplem} we collect some technical results for later use, and in Subsection \ref{FqBoundaryIdeal} we give an explicit description of the ideal of the boundary. The arguments in the last two sections are simpler versions of central arguments from Section \ref{Specialn}; it should be helpful for the reader to study them here first.

Section \ref{A} contains all general definitions and results concerning $A$-reciprocal maps and their moduli schemes, and Section \ref{Specialn} contains our results in the special case  $A=\BF_q[t]$ and $N=(t^n)$. The most notable content of these sections was already summarized above. 

%%%%%%%%%%%%%%%%%%%%%%%%%%%%%%%%%%%%%

\medskip
{\bf Acknowledgements:} 
It is my pleasure to acknowledge interesting and helpful conversations with Ernst Gekeler, Simon H\"aberli, and Maxim Mornev, and many valuable comments on earlier versions of this paper by Simon H\"aberli.

\newpage

%%%%%%%%%%%%%%%%%%%%%%%%%%%%%%%%%%%%%%%%%%%%%%%%%%%%%%%%%%%%%%%%%%%%%%
%%%%%%%%%%%%%%%%%%%%%%%%%%%%%%%%%%%%%%%%%%%%%%%%%%%%%%%%%%%%%%%%%%%%%%

\section{$\BF_q$-reciprocal maps}
\label{Fq}

%%%%%%%%%%%%%%%%%%%%%%%%%%%%%%%%%%%%%%%%%%%%%%%%%%%%%%%%%%%%%%%%%%%%%%

\subsection{Basic $\BF_q$-reciprocal maps}
\label{BasicFqRecip}

We begin by reviewing some basic constructions from \cite[\S1 and \S7]{PinkSchieder}.
Throughout this article we fix a finite field $\BF_q$ of order~$q$. 
For any $\BF_q$-vector space $V$ we abbreviate $\circV := V\setminus\{0\}$. 
Consider a finite dimensional $\BF_q$-vector space $V$ and a commutative $\BF_q$-algebra~$R$.

\begin{Def}\label{FqRecipDef}
%(Compare \cite[Def.\;7.2]{PinkSchieder})
A map $\rho\colon\circV\to R$ is called \emph{$\BF_q$-reciprocal} if 
\begin{itemize}
\item[(a)] $\rho(v) \cdot \rho(w) = \rho(v+w) \cdot (\rho(v) + \rho(w))$ for all $v$, $w\in \circV$ with $v+w\in\circV$, and
\item[(b)] $\alpha\rho(\alpha v) = \rho(v)$ for all $\alpha \in \BF_q^{\times}$ and $v \in \circV$.
\end{itemize}
\end{Def}

\begin{Def}\label{FiberConds}
%(See \cite[Def.\;7.3]{PinkSchieder})
An $\BF_q$-reciprocal map $\rho: \circV\to R$ is called 
\begin{itemize}
\item[(a)] \emph{fiberwise non-zero} if for every $\Fp\in\Spec(R)$ there exists $v\in\circV$ with $\rho(v)\not\in\Fp$.
\item[(b)] \emph{fiberwise invertible} if for every $\Fp\in\Spec(R)$ and every $v\in\circV$ we have $\rho(v)\not\in\Fp$.\\
Equivalently: If for every $v\in\circV$ we have $\rho(v)\in R^\times$.
\end{itemize}
\end{Def}

\begin{Prop}\label{FqLinearMap}
%(See \cite[Prop.\;7.4]{PinkSchieder})
For any $\BF_q$-linear map $\lambda\colon V\to R$ satisfying $\lambda(\circV\;)\subset R^\times$, the map $\rho\colon\circV\to R$, $v\mapsto \lambda(v)^{-1}$, is a fiberwise invertible $\BF_q$-reciprocal map, and any fiberwise invertible $\BF_q$-reciprocal map $\circV\to R$ arises in this way.
\end{Prop}

\begin{Cons}\label{RVCons}
\rm We set
\begin{eqnarray}
\nonumber  S_V   & := & \text{the symmetric algebra of $V$ over $\BF_q$,} \\
\nonumber  K_V   & := & \mbox{the field of quotients of $S_V$,}\\
% \Quot(S_V) \mbox{\ the field of quotients of $S_V$,}\\
\nonumber  R_V & := & \mbox{the $\BF_q$-subalgebra of $K_V$ generated by the elements $\frac{1}{v}$ for all $v \in \circV$,} \\
\nonumber  RS_V & := & \mbox{the $\BF_q$-subalgebra of $K_V$ generated by $R_V$ and $S_V$.}
\end{eqnarray}
Note that these are all integral domains, and $RS_V$ is the localization of $R_V$ obtained by inverting the elements $\frac{1}{v}$ for all $v\in\circV$. 
For any basis $X_1,\ldots,X_r$ of~$V$, the field $K_V$ becomes the field of rational functions $\BF_q(X_1,\ldots,X_r)$ and $R_V$ becomes the $\BF_q$-subalgebra generated by the elements $(\alpha_1X_1+\ldots+\alpha_rX_r)^{-1}$ for all $(\alpha_1,\ldots,\alpha_r)\in\BF_q^r\setminus\{(0,\ldots,0)\}$.
\end{Cons}

By \cite[\S1 and Thm.\;7.12]{PinkSchieder} we have:

\begin{Thm}\label{FqModuli}
%(See \cite[Thms.\;7.12, 1.7]{PinkSchieder})
%(See \cite[\S1 and Thm.\;7.12]{PinkSchieder})
\begin{itemize}
\item[(a)] The map 
$$\rho^\univ\colon \circV\longto R_V,\ v\mapsto\tfrac{1}{v}$$ 
is $\BF_q$-reciprocal.
\item[(b)] For any $\BF_q$-algebra $R$ and any $\BF_q$-reciprocal map $\rho\colon\circV\to R$ there exists a unique $\BF_q$-algebra homomorphism $f\colon R_V\to R$ such that $\rho=f\circ\rho^\univ$.
\item[(c)] This $f$ extends to a ring homomorphism $RS_V\to R$ if and only if $\rho$ is fiberwise invertible.
\end{itemize}
\end{Thm}

%%%%%%%%%%%%%%%%%%%%%%%%%%%%%%%%%%%%%%%%%%%%%%%%%%%%%%%%%%%%%%%%%%%%%%

\subsection{Functoriality}
\label{Functoriality}

Consider a short exact sequence of finite dimensional $\BF_q$-vector spaces
\UseTheoremCounterForNextEquation
\begin{equation}\label{ShExSeq}
\xymatrix{0\ar[r]&V'\ar[r]^i&V\ar[r]^p&V''\ar[r]&0\rlap{.}\\}
\end{equation}

\begin{PropDef}\label{Pullback}
For any $\BF_q$-reciprocal map $\rho\colon\circV\to R$ the~map
$$i^*\rho := \rho\circ i\colon \circVprime\to R,\ v'\mapsto \rho(i(v'))$$
is $\BF_q$-reciprocal. We call it the \emph{pullback of $\rho$ under~$i$}.
\end{PropDef}

\begin{Proof}
Clear from Definition \ref{FqRecipDef}, the injectivity ensuring that $i$ induces a map ${\circVprime\to\circV}$.
\end{Proof}

\begin{PropDef}\label{ExtByZero}
For any $\BF_q$-reciprocal map $\rho\colon\circVprime\to R$ the map
$$i_*\rho\!:\ \circV\to R, \ v \mapsto 
\left\{{\huge\mathstrut}\right.\kern-5pt
\begin{array}{cl}
\rho(v) & \hbox{if $v=i(v')$ for $v'\in\circVprime$,} \\[3pt]
0 & \hbox{if $v\not\in i(\circVprime\,)$,}
\end{array}$$
is $\BF_q$-reciprocal. We call it the \emph{extension by zero} or the \emph{pushforward of $\rho$ under~$i$}.
\end{PropDef}

\begin{Proof}
Clearly $i_*\rho$ satisfies the condition \ref{FqRecipDef} (b). 
It also satisfies \ref{FqRecipDef} (a) whenever $v$, $w$, $v+w$ lie in $i(\circVprime)$. In all other cases at least two of $v$, $w$, $v+w$ lie in $V\setminus i(V')$; hence at least two of the values $(i_*\rho)(v)$, $(i_*\rho)(w)$, $(i_*\rho)(v+w)$ are zero, and the equality in \ref{FqRecipDef} (a) for $i_*\rho$ holds trivially. Thus the extension by zero is $\BF_q$-reciprocal.
\end{Proof}

\begin{Prop}\label{i*Homo}
\begin{itemize}
\item[(a)] The functor $i^*$ is represented by an injective $\BF_q$-algebra homomorphism $\epsilon_i\colon R_{V'}\into R_V$ that sends $\frac{1}{v'}$ to $\frac{1}{i(v')}$ for all $v'\in\circVprime$.
\item[(b)] The functor $i_*$ is represented by a surjective $\BF_q$-algebra homomorphism $\pi_i\colon {R_V\onto R_{V'}}$ that sends $\frac{1}{i(v')}$ to $\frac{1}{v'}$ for all $v'\in\circVprime$ and $\frac{1}{v}$ to $0$ for all $v\in V\setminus i(V')$.
\item[(c)] The kernel of $\pi_i$ is generated by the elements $\frac{1}{v}$ for all $v\in V\setminus i(V')$.
\end{itemize}
\end{Prop}

\begin{Proof}
Let $\rho_V^\univ\colon \circV\to R_V$ and $\rho_{V'}^\univ\colon \circVprime\to R_{V'}$ denote the respective universal $\BF_q$-reciprocal maps. Then $\epsilon_i$ and $\pi_i$ are obtained from the universal property of $(R_{V'},\rho_{V'}^\univ)$ and $(R_V,\rho_V^\univ)$ as the unique $\BF_q$-algebra homomorphisms making the following diagrams commute:
$$\xymatrix{& \ \circVprime\ \ar[dl]_-{\rho_{V'}^\univ} \ar[dr]^-{i^*\rho_V^\univ} & \\
R_{V'} \ar[rr]^-{\epsilon_i} && R_V\rlap{,}}
\qquad\qquad
\xymatrix{& \ \circV\ \ar[dl]_-{\rho_V^\univ} \ar[dr]^-{i_*\rho_{V'}^\univ} & \\
R_V \ar[rr]^-{\pi_i} && R_{V'}\rlap{.}}$$
By construction they represent the functors $i^*$ and $i_*$ and are given on the generators as stated. Since $i^*i_*\rho_{V'}^\univ = \rho_{V'}^\univ$, the universal property of $(R_{V'},\rho_{V'}^\univ)$ implies that $\pi_i\circ\epsilon_i=\id_{R_{V'}}$; hence $\epsilon_i$ is injective and $\pi_i$ is surjective. This proves (a) and (b).

For (c) let $J\subset R_V$ denote the ideal generated by the elements $\frac{1}{v}$ for all $v\in V\setminus i(V')$. Then the factor ring $R_V/J$ represents the functor of all $\BF_q$-reciprocal maps on $\smash{\circV}$ which are identically zero on $V\setminus i(V')$. But these are precisely the extensions by zero of $\BF_q$-reciprocal maps on $\smash{\circVprime}$; hence this functor is already represented by~$R_{V'}$. It follows that $\pi_i$ induces an isomorphism $R_V/J\isoto R_{V'}$; proving (c).
\end{Proof}

%%%%%%%%%%%%%%%%%%%%%%%%%%%%%%%%%%%%%

\begin{PropDef}\label{Quotient}
For any $\BF_q$-reciprocal map $\rho\colon\circV\to R$ the map
$$p_*\rho\colon \circVpprime\to R,\
v'' \mapsto \kern-10pt\sum_{v\in p^{-1}(v'')}\kern-10pt\rho(v)$$
is $\BF_q$-reciprocal. We call it the \emph{pushforward of $\rho$ under~$p$}.
\end{PropDef}

\begin{Proof}
Setting $U:=\Ker(p)$, the condition \ref{FqRecipDef} (a) for $p_*\rho$ is equivalent to the formula
$$\sum_{u\in U}\rho(v+u) \cdot \sum_{u'\in U} \rho(w+u')
\ =\ \sum_{u\in U}\rho(v+w+u) \cdot \sum_{u'\in U}(\rho(v+u') + \rho(w+u'))$$
for all $v$, $w\in V\setminus U$ with $v+w\in V\setminus U$. This equation follows from the condition \ref{FqRecipDef} (a) for~$\rho$ by rearranging and reindexing the sums. Also, condition \ref{FqRecipDef} (b) for $p_*\rho$ immediately follows from that for~$\rho$. Thus $p_*\rho$ is $\BF_q$-reciprocal.
\end{Proof}

\begin{Prop}\label{p*Homo}
The functor $p_*$ is represented by an injective $\BF_q$-algebra homomorphism $\epsilon_p\colon R_{V''}\into R_V$ that sends $\frac{1}{v''}$ to $\sum_{v\in p^{-1}(v'')}\frac{1}{v}$ for all $v''\in\circVpprime$.
\end{Prop}

\begin{Proof}
Let $\rho_{V''}^\univ\colon \circVpprime\to R_{V''}$ denote the universal $\BF_q$-reciprocal map. Then $\epsilon_p$ is obtained from the universal property of $(R_{V''},\rho_{V''}^\univ)$ as the unique $\BF_q$-algebra homomorphisms making the following diagram commute:
$$\xymatrix{& \ \circVpprime\ \ar[dl]_-{\rho_{V''}^\univ} \ar[dr]^-{p_*\rho_V^\univ} & \\
R_{V''} \ar[rr]^-{\epsilon_p} && R_V\rlap{.}}$$
By construction this represents the functor $p_*$ and is given on the generators as stated. 
To finish choose a homomorphism $j\colon V''\into V$ such that $p\circ j=\id_{V''}$. Then the defining formulas in Propositions \ref{ExtByZero} and \ref{Quotient} show that $p_*j_*\rho_{V''}^\univ = \rho_{V''}^\univ$. Thus the universal property of $(R_{V''},\rho_{V''}^\univ)$ implies that $\pi_j\circ\epsilon_p=\id_{R_{V''}}$; hence $\epsilon_p$ is injective, and we are done.
\end{Proof}

\begin{Rem}\label{PushforwardGeneral}
\rm Factoring an arbitrary homomorphism of finite dimensional $\BF_q$-vector spaces as $f=i\circ p$ for a surjection $p$ and an injection~$i$, one can define the \emph{pushforward under $f$} by $f_* := i_*\circ p_*$. This is actually given by the same formula as in Proposition \ref{Quotient} for $f$ in place of~$p$.
\end{Rem}

%We will say more about the pushforward $p_*$ in the next subsection.

%%%%%%%%%%%%%%%%%%%%%%%%%%%%%%%%%%%%%%%%%%%%%%%%%%%%%%%%%%%%%%%%%%%%%%

\subsection{Some nice formulas}
\label{NiceForm}

As before let $R$ be a commutative $\BF_q$-algebra. 
Let $R[\tau]$ denote the ring of \emph{$\BF_q$-linear polynomials over~$R$}, that is, of polynomials of the form $f(X)=\sum'_{i\ge0}u_iX^{q^i}$ with all $u_i\in R$. Setting $\tau(X) := X^q$, we write such a polynomial in the shorter form $f=\sum'_{i\ge0}u_i\tau^i$. The multiplication in $R[\tau]$ is defined as composition $f\circ g$, and the identity element $1$ of $R[\tau]$ is the polynomial $\tau^0=X$. For any $u\in R$ we have $\tau\circ u = u^q\circ\tau$; so in general this ring is non-commutative. 
For any $f=\sum'_{i\ge0}u_i\tau^i \in R[\tau]$ we have $df := \frac{d}{dX}f(X) = u_0$, and the map $d\colon R[\tau]\to R$ is an $\BF_q$-algebra homomorphism.

\medskip
For the rest of this subsection we fix an $\BF_q$-reciprocal map $\rho\colon\circV\to R$. To $\rho$ we associate the polynomial
\UseTheoremCounterForNextEquation
\begin{equation}\label{ExpDef}
e_\rho(X)\ :=\ X\cdot\prod_{v\in\circV}(1-\rho(v)X)\ \in\ R[X].\end{equation}

\begin{Prop}\label{ExpFqLinear}
%For any $\BF_q$-reciprocal map $\rho: \circV\to R$ we have
We have $e_\rho\in R[\tau]$. 
\end{Prop}

\begin{Proof}
By Theorem \ref{FqModuli} it suffices to prove this for the universal $\BF_q$-reciprocal map $\rho^\univ\colon{\circV\to R_V}$. Since $R_V$ is an integral domain, it then suffices to prove the statement over the quotient field of~$R_V$, where $\rho^\univ$ becomes fiberwise invertible.
%by Theorem \ref{FqModuli} (b). 
In view of Proposition \ref{FqLinearMap} it thus suffices to prove the statement for the map $\rho=(\lambda|\circV\;)^{-1}$ associated to any injective $\BF_q$-linear map $\lambda\colon V\into k$ for any field $k$ over~$\BF_q$. In that case
%\UseTheoremCounterForNextEquation
%\begin{equation}\label{ExpDefLambda}
$$e_\rho(X)\ =\ X\cdot\prod_{v\in\circV}\left(1-\frac{X}{\lambda(v)}\right),$$
%\end{equation}
is the \emph{exponential function} associated to the subgroup $\lambda(V)\subset k$, and the statement follows from \cite[Cor.\,1.2.2]{GossBS}.
\end{Proof}

%\begin{Prop}\label{LogDeriv}
%We have the following identity in $R[[X]]$:
%$$\frac{1}{e_\rho(X)}\ =\ \frac{1}{X} + \sum_{v\in\circV} \frac{-\rho(v)}{1-\rho(v)X}.$$
%\end{Prop}
%
%\begin{Proof}
%(Compare Goss \cite[\S6]{GossAlg}.) Since $e_\rho(X)$ is an $\BF_q$-linear polynomial with linear term $X$, we have ${\frac{d}{dX}e_\rho(X)=1}$. Thus the desired identity follows by applying the logarithmic derivative to the formula (\ref{ExpDef}).
%\end{Proof}

\begin{Prop}\label{ERhoIdentities}
We have the following identities in the ring $R[X,e_\rho(X)^{-1}]$:
\begin{itemize}
\item[(a)] $\displaystyle\frac{1}{e_\rho(X)}\ =\ \frac{1}{X} + \sum_{v\in\circV} \frac{-\rho(v)}{1-\rho(v)X}\ $.
\item[(b)] $\displaystyle\prod_{v\in\circV}\frac{1}{1-\rho(v)X}\ =\ \frac{X}{e_\rho(X)}\ =\ -\sum_{v\in\circV}\frac{1}{1-\rho(v)X}$ \ if $V\not=0$.
\end{itemize}
\end{Prop}

\begin{Proof}
(Compare Goss \cite[\S6]{GossAlg}.) Since $e_\rho(X)$ is an $\BF_q$-linear polynomial with linear term~$X$, we have ${\frac{d}{dX}e_\rho(X)=1}$. Applying the logarithmic derivative to the formula (\ref{ExpDef}) thus yields the equation (a). Multiplying it by $X$ implies that
$$\frac{X}{e_\rho(X)}\ =\ 
1 + \sum_{v\in\circV} \frac{1-\rho(v)X-1}{1-\rho(v)X}
\ =\ |V| + \sum_{v\in\circV} \frac{-1}{1-\rho(v)X}.$$
Since $q$ divides $|V|$ if $V\not=0$, this implies (b).
\end{Proof}

%%%%%%%%%%%%%%%%%%%%%%%%%%%%%%%%%%%%%

\medskip
Now we return to the short exact sequence (\ref{ShExSeq}). For simplicity we assume that $i$ is the inclusion of an $\BF_q$-subspace ${V'\into V}$.

\begin{Prop}\label{RhoFormulas}
%Suppose that $i$ is the inclusion of a subspace $V'\into V$. Then for
For any $v\in V\setminus V'$ we have 
\begin{itemize}
\item[(a)] $\displaystyle \Bigl(\sum_{v'\in V'}\rho(v-v')\Bigr) \cdot \Bigl(\prod_{v'\in\circVprime}\rho(v')\Bigr) \ =\ \prod_{v'\in V'}\rho(v-v')$.
\item[(b)] $\displaystyle \Bigl(\sum_{v'\in V'}\rho(v-v')\Bigr) \cdot \Bigl(\prod_{v'\in\circVprime}\bigl(\rho(v)-\rho(v')\bigr)\Bigr) \ =\ \rho(v)^{|V'|}.$
\item[(c)] $\displaystyle \Bigl(\sum_{v'\in V'}\rho(v-v')\Bigr) \cdot e_{i^*\rho}\Bigl(\frac{1}{\rho(v)}\Bigr)\ =\ 1$ \ if \ $\rho(v)\in R^\times$.
\end{itemize}
\end{Prop}

\begin{Proof}
As in the proof of Proposition \ref{ExpFqLinear}, showing these equations reduces to the case that $\rho=(\lambda|\circV\;)^{-1}$ for an injective $\BF_q$-linear map $\lambda\colon V\into k$ to a field~$k$. Applying 
%Proposition \ref{LogDeriv}
Proposition \ref{ERhoIdentities} (a) to $i^*\rho$ in place of $\rho$ then shows that
\UseTheoremCounterForNextEquation
\begin{equation}\label{RhoFormulas1}
\frac{1}{e_{i^*\rho}(X)}\ =\ \frac{1}{X} + \sum_{v'\in\circVprime} \frac{-\rho(v')}{1-\rho(v)X}
\ =\ \frac{1}{X} + \sum_{v'\in\circVprime} \frac{-1}{\lambda(v')-X}
\ =\ \sum_{v'\in V'} \frac{1}{X-\lambda(v')}.
\end{equation}
For any $v\in V\setminus V'$ we have 
\UseTheoremCounterForNextEquation
\begin{equation}\label{RhoFormulas2}
e_{i^*\rho}(\lambda(v))\ =\ \lambda(v)\cdot\prod_{v'\in\circVprime}\left(1-\frac{\lambda(v)}{\lambda(v')}\right)
\ =\ \lambda(v)\cdot\prod_{v'\in\circVprime}\frac{\lambda(v')-\lambda(v)}{\lambda(v')},
\end{equation}
where all factors are non-zero by the injectivity of~$\lambda$.
Thus by multiplying the formula (\ref{RhoFormulas1}) by $e_{i^*\rho}(X)$, substituting $X=\lambda(v)$, and using the additivity of $\lambda$ we deduce that
\UseTheoremCounterForNextEquation
\begin{equation}\label{RhoFormulas3}
1\ =\ \left(\sum_{v'\in V'} \frac{1}{\lambda(v)-\lambda(v')}\right) \cdot e_{i^*\rho}(\lambda(v)) 
\ =\ \left(\sum_{v'\in V'}\rho(v-v')\right) \cdot e_{i^*\rho}\left(\frac{1}{\rho(v)}\right),
\end{equation}
proving (c). Also, by (\ref{RhoFormulas2}), the additivity of $\lambda$, and the fact that $(-1)^{q-1}=1$ in $\BF_q$ we have 
$$e_{i^*\rho}\left(\frac{1}{\rho(v)}\right)
\ =\ e_{i^*\rho}(\lambda(v))\ =\ (-1)^{|\circVprime|} \cdot \frac{\displaystyle\prod\nolimits_{v'\in V'}\lambda(v-v')}{\displaystyle\prod\nolimits_{v'\in\circVprime}\lambda(v')}
\ =\ \frac{\displaystyle\prod\nolimits_{v'\in\circVprime}\rho(v')}{\displaystyle\prod\nolimits_{v'\in V'}\rho(v-v')},$$
which together with (\ref{RhoFormulas3}) proves (a). Finally, we can rewrite (\ref{RhoFormulas2}) also in the form 
$$e_{i^*\rho}\left(\frac{1}{\rho(v)}\right)
\ =\ \frac{1}{\rho(v)}\cdot\prod_{v'\in\circVprime}\left(1-\frac{\rho(v')}{\rho(v)}\right)
\ =\ \frac{1}{\rho(v)^{|V'|}} \cdot 
\prod_{v'\in\circVprime}(\rho(v)-\rho(v')).$$
Combined with (\ref{RhoFormulas3}) this proves (b).
\end{Proof}

\begin{Exp}\label{RhoFormulasRem}
\rm It is well-known that an injective $\BF_q$-linear map $\lambda\colon V\into k$ induces an injective $\BF_q$-linear map $p_*\lambda\colon V''\into k$ by the formula $(p_*\lambda)(p(v)) := e_{i^*(\lambda|\circV\;)^{-1}}(v)$. The formula in Proposition \ref{RhoFormulas} (c) translates this equation into the surprisingly simple formula for the reciprocals
\UseTheoremCounterForNextEquation
\begin{equation}\label{RhoFormulas4}
\frac{1}{(p_*\lambda)(p(v))}
\ =\ \sum_{v'\in V'}\frac{1}{\lambda(v-v')}
\end{equation}
for all $v\in V\setminus V'$. The fact that the right hand side is a polynomial in the values of $(\lambda|\circV\,)^{-1}$ is central for everything that follows in this article. It motivated both the quotient construction for reciprocal maps in Definition \ref{Quotient} and the definition of $A$-reciprocal maps in Definition \ref{ARecipDef}.
\end{Exp}

\begin{Prop}\label{PushForwardInjective}
\begin{itemize}
\item[(a)] If $\rho$ is fiberwise invertible, then so is $p_*\rho$.
\item[(b)] If 
%$\rho|(V\setminus V')$ is fiberwise non-zero, i.e., if 
for every $\Fp\in\Spec(R)$ there exists $v\in V\setminus V'$ with $\rho(v)\not\in\Fp$, then $p_*\rho$ is fiberwise non-zero.
\end{itemize}
\end{Prop}

\begin{Proof}
For any $\Fp\in\Spec(R)$ and any $v\in V\setminus V'$ with $\rho(v)\not\in\Fp$, Proposition \ref{RhoFormulas} (b) and the definition of $p_*\rho$ imply that $(p_*\rho)(p(v))\not\in\Fp$. By applying the universal quantifier $\forall$ or the existential quantifier $\exists$ to $v$ we obtain the respective result.
\end{Proof}

\begin{Prop}\label{Composition}
We have the following identities in $R[\tau]$:
\begin{itemize}
\item[(a)] $e_\rho\circ u = u\circ e_{u\rho}$ for any $u\in R$.
\item[(b)] $e_\rho = e_{p_*\rho}\circ e_{i^*\rho}$.
\end{itemize}
\end{Prop}

\begin{Proof}
(a) follows by direct computation from (\ref{ExpDef}).
%{}From (\ref{ExpDef}) we directly obtain
%$$u\cdot e_{u\rho}(X)\ =\ uX\cdot\prod_{v\in\circV}(1-u\rho(v)X)
%\ =\ uX\cdot\prod_{v\in\circV}(1-\rho(v)\cdot uX)
%\ =\ e_\rho(uX),$$
%proving (a). 
To establish (b) we may, as in the proof of Proposition \ref{ExpFqLinear}, reduce ourselves to the case that $\rho=(\lambda|\circV\;)^{-1}$ for an injective $\BF_q$-linear map $\lambda\colon V\into k$ to a field~$k$. Since $e_{i^*\rho}$ is an $\BF_q$-linear polynomial with the precise set of zeros $\lambda(i(V'))$, there is a unique injective $\BF_q$-linear map $\lambda''\colon V''\into k$ satisfying $\lambda''(p(v)) = e_{i^*\rho}(\lambda(v))$ for all $v\in V$. But by Proposition \ref{RhoFormulas} (c) and the definition of $p_*\rho$ we have $(p_*\rho)(p(v)) = e_{i^*\rho}(\lambda(v))^{-1}$ for all $v\in V\setminus V'$. Thus $p_*\rho = (\lambda''|\circVpprime)^{-1}$. Now the formula (b) is precisely that in \cite[(1.12)]{GekelerPowerSums}.
\end{Proof}

%%%%%%%%%%%%%%%%%%%%%%%%%%%%%%%%%%%%%%%%%%%%%%%%%%%%%%%%%%%%%%%%%%%%%%

\subsection{General $\BF_q$-reciprocal maps}
\label{GenFqRecip}

Now we globalize the concept of $\BF_q$-reciprocal maps following \cite[\S7]{PinkSchieder}. For this we assume that $V\not=0$. Let $S$ be a scheme over~$\BF_q$, let $\CL$ be an invertible sheaf on $S$, and let $\CL(S)$ denote the space of global sections of~$\CL$. For any section $\ell\in\CL(S)$ and any point $s\in $ we let $\ell(s)\in\CL\otimes_{\CO_S}k(s)$ denote the value of $\ell$ over the residue field $k(s)$ of~$s$. The (tensor) product of sections $\ell_1,\ldots,\ell_n\in\CL(S)$ is a section $\ell_1\cdots\ell_n\in\CL^{\otimes n}(S)$, and the inverse of a nowhere vanishing section $\ell\in\CL(S)$ is a section $\ell^{-1}\in\CL^\vee(S)$. 

\begin{Def}\label{GenFqRecipDef}
%(Compare \cite[Def.\;7.2]{PinkSchieder})
A map $\rho\colon\circV\to\CL(S)$ is called \emph{$\BF_q$-reciprocal} if 
\begin{itemize}
\item[(a)] $\rho(v) \cdot \rho(w) = \rho(v+w) \cdot (\rho(v) + \rho(w))$ in $\CL^{\otimes2}(S)$ for all $v$, $w\in \circV$ with $v+w\in\circV$, and
\item[(b)] $\alpha\rho(\alpha v) = \rho(v)$ for all $\alpha \in \BF_q^{\times}$ and $v \in \circV$.
\end{itemize}
\end{Def}

\begin{Def}\label{GenFiberConds}
%(See \cite[Def.\;7.3]{PinkSchieder})
An $\BF_q$-reciprocal map $\rho\colon \circV\to\CL(S)$ is called 
\begin{itemize}
\item[(a)] \emph{fiberwise non-zero} if for every point $s\in S$ there exists $v\in\circV$ with $\rho(v)(s)\not=0$.
\item[(b)] \emph{fiberwise invertible} if for every point $s\in S$ and every $v\in\circV$ we have $\rho(v)(s)\not=0$.
%\\
%Equivalently: If for every $v\in\circV$ the section $\rho(v)$ is invertible.
\end{itemize}
\end{Def}

\begin{Rem}\label{Gen=StdFq}
\rm When $S=\Spec(R)$ and $\CL=\CO_X$, Definitions \ref{GenFqRecipDef} and \ref{GenFiberConds} agree precisely with Definitions \ref{FqRecipDef} and \ref{FiberConds}. 
Conversely, for an arbitrary invertible sheaf $\CL$ consider a covering of $S$ by open affines $U_i=\Spec(R_i)$ and an isomorphism $f_i\colon\CO_S|U_i\isoto\CL|U_i$ for each~$i$. Giving an $\BF_q$-reciprocal map $\rho\colon\circV\to\CL(S)$ is then equivalent to giving $\BF_q$-reciprocal maps $\rho_i\colon\circV\to R_i=\CO_{S}(U_i)$ for all $i$ such that $f_i\circ\rho_i = f_j\circ\rho_j$ over $U_i\cap U_j$ for all $i,j$.

Thus all the results from Subsection \ref{NiceForm} have direct analogues in this more general setting. We also obtain a moduli space, as follows.
\end{Rem}

\begin{Cons}\label{GenRVCons}
\rm We endow the rings $R_V$ and $RS_V$ from Construction \ref{RVCons} with the unique $\BZ$-grading for which the elements $\frac{1}{v}$ are homogeneous of degree $1$ for all $v\in\circV$. 
For any integer $d$ let $R_{V,d}$ and $RS_{V,d}$ denote the respective homogenous parts of degree~$d$. By construction $R_V$ is generated over $\BF_q$ by its homogeneous part of degree~$1$. Thus 
$$Q_V \ :=\ \Proj(R_V)$$
is a projective scheme over $\BF_q$ endowed with a natural very ample invertible sheaf $\CO(1)$ and a natural homomorphism $R_{V,d} \to \CO(d)(Q_V)$ for all $d\in\BZ$. In fact this is an isomorphism by \cite[Cor.\,5.4]{PinkSchieder}.

Since $R_V$ is an integral domain and we have now assumed $V\not=0$, we have $R_V\not=\BF_q$ and $Q_V$ is an integral scheme.
Also, since $RS_V$ is the localization of $R_V$ obtained by inverting a non-empty finite set of homogeneous elements of degree~$1$, the scheme
$$\Omega_V\ :=\ \Proj(RS_V)\ \cong\ \Spec(RS_{V,0})$$ 
is an affine open dense subscheme of~$Q_V$. 
\end{Cons}

\begin{Def}\label{FqRecipIsoms}
Consider two pairs $(\CL, \rho)$ and $(\CL', \rho')$ consisting of an invertible sheaf and an $\BF_q$-reciprocal map. An isomorphism of invertible sheaves $f\colon \CL\isoto\CL'$ satisfying $\rho'=f\circ\rho$ is called an \emph{isomorphism} $(\CL, \rho) \isoto(\CL', \rho')$. If there exists such an isomorphism, the pairs $(\CL, \rho)$ and $(\CL', \rho')$
%, or by abuse of terminology, the $\BF_q$-reciprocal maps $\rho$ and $\rho'$,
are called \emph{isomorphic}.
\end{Def}

If $\rho$ or $\rho'$ is fiberwise non-zero, there exists at most one isomorphism $(\CL, \rho) \isoto(\CL', \rho')$. Thus the isomorphism classes of such pairs form a well-posed moduli problem. By \cite[Thm.\;7.10 and Prop.\,7.11]{PinkSchieder} we have:

\begin{Thm}\label{GenFqModuli}
\begin{itemize}
%\item[(a)] The composite map 
%$$\xymatrix@C-10pt{
%\llap{$\rho^\univ\colon$}\ \circV
%\ar[rrr]^-{\textstyle v\mapsto\frac{1}{v}} &&& 
%R_{V,1} \ar[r] & \CO(1)(Q_V)}$$
\item[(a)] The map 
$$\rho^\univ\colon \circV \longto R_{V,1}=\CO(1)(Q_V),
\ v\mapsto\tfrac{1}{v}$$
is $\BF_q$-reciprocal and fiberwise non-zero.
\item[(b)] For any scheme $S$ over $\BF_q$, any invertible sheaf $\CL$ on~$S$, and any fiberwise non-zero $\BF_q$-reciprocal map $\rho\colon\circV\to \CL(S)$ there exists a unique morphism $f\colon S\to Q_V$ over~$\BF_q$ such that $(\CL,\rho)\cong f^*(\CO(1),\rho^\univ)$.
\item[(c)] This $f$ factors through $\Omega_V$ if and only if $\rho$ is fiberwise invertible.
\end{itemize}
\end{Thm}

%%%%%%%%%%%%%%%%%%%%%%%%%%%%%%%%%%%%%
\medskip
We end this subsection by discussing the effect of the functors from Subsection \ref{Functoriality}. Consider the exact sequence (\ref{ShExSeq}).

\begin{Prop}
\label{i*Morph1}
\begin{itemize}
\item[(a)] The homomorphism $\pi_i\colon {R_V\onto R_{V'}}$ from Proposition \ref{i*Homo} (b) induces a closed embedding $\epsilon_i\colon Q_{V'}\into Q_V$, whose image is defined by the equations $\rho^\univ(v)=0$ for all $v\in V\setminus i(V')$. 
\item[(b)] Consider any fiberwise non-zero $\BF_q$-reciprocal map $\rho\colon\circV\to \CL(S)$ over a scheme $S$ over~$\BF_q$. Then the associated morphism $S\to Q_V$ factors through $\epsilon_i$ if and only if $\rho=i_*\rho'$ for an $\BF_q$-reciprocal map $\rho'\colon\circVprime\to \CL(S)$.
\end{itemize}
\end{Prop}

\begin{Proof}
The description in Proposition \ref{i*Homo} shows that $\pi_i$ is a surjective graded $\BF_q$-algebra homomorphism whose kernel is generated by the elements $\frac{1}{v}$ for all $v\in V\setminus i(V')$. This directly implies (a). Part (b) follows as in the proof of Proposition \ref{i*Homo} (c).
\end{Proof}

\begin{Rem}\label{i*Morph2}
\rm Let $X$ be the closed subscheme of $Q_V$ that is defined by the equations ${\rho^\univ(i(v'))=0}$ for all $v'\in\circVprime$. Then the pullback $i^*\rho^\univ$ is fiberwise non-zero over $Q_V\setminus X$; hence by the universal property of $Q_{V'}$ it corresponds to a morphism $Q_V\setminus X\to Q_{V'}$. In fact, this is the morphism induced by the graded $\BF_q$-algebra homomorphism $\epsilon_i\colon R_{V'}\into R_V$ from Proposition \ref{i*Homo} (a). By the same argument as in the proof of Proposition \ref{i*Homo} the morphism $Q_V\setminus X\to Q_{V'}$ is a left inverse of the embedding $\epsilon_i\colon Q_{V'}\into Q_V$ above.
\end{Rem}

\begin{Rem}\label{p*Morph}
\rm Using the description of $\epsilon_i(Q_{V'})$ in Proposition \ref{i*Morph1} (a), Proposition \ref{PushForwardInjective} (b) implies that the pushforward $p_*\rho^\univ$ is fiberwise non-zero over $Q_V\setminus \epsilon_i(Q_{V'})$. By the universal property of $Q_{V''}$ it therefore corresponds to a morphism $Q_V\setminus \epsilon_i(Q_{V'}) \to Q_{V''}$. In fact, this is the morphism induced by the graded $\BF_q$-algebra homomorphism $\epsilon_p\colon R_{V''}\into R_V$ from Proposition \ref{p*Homo}. 
\end{Rem}

% Are these morphisms flat? 

%%%%%%%%%%%%%%%%%%%%%%%%%%%%%%%%%%%%%%%%%%%%%%%%%%%%%%%%%%%%%%%%%%%%%%

\subsection{Description of the ring}
\label{FqDescRing}

For later use we recall the description of $R_V$ from Section 2 of \cite{PinkSchieder}. Choose a basis $X_1,\ldots,X_r$ of~$V$, and for each $0\le k\le r$ consider the subspace $V_k := \BF_q X_1+\ldots+\BF_q X_k$. For every $1\le k\le r$ consider the finite subsets
\UseTheoremCounterForNextEquation
\begin{eqnarray}
\label{DkDef}
\Delta_k &:=& \biggl\{ \frac{1}{X_k + w} \,\biggm|\, w \in \circV_{k-1} \biggr\} \cup\bigl\{1\bigr\} \quad\hbox{and} \\
\UseTheoremCounterForNextEquation
\label{EkDef}
E_k &:=& \biggl\{ \frac{1}{X_k + w} \,\biggm|\, w \in V_{k-1} \biggr\}
\end{eqnarray}
of~$R_V$, each of cardinality $|V_{k-1}| = q^{k-1}$. Note that the $\Delta_k$ differs from $E_k$ only in that the element $\frac{1}{X_k}$ is replaced by~$1$. Observe that by unique factorization for polynomials in the variables $X_1,\ldots,X_r$ we have bijective maps
$$\xymatrix@R-20pt@C-13pt{
\Delta_1\times\ldots\times\Delta_r \ 
\ar[rr]^-\sim && \ \Delta_1\cdots\Delta_r \ar@{}[r]|-{:=} 
& \bigl\{e_1\cdots e_r \bigm| \forall k\colon e_k\in \Delta_k\bigr\} 
\ar@{}[r]|-{\subset} & R_V \rlap{\hbox{\qquad and}}\\
E_1\times\ldots\times E_r \ 
\ar[rr]^-\sim && \ E_1\cdots E_r \ar@{}[r]|-{:=} 
& \bigl\{e_1\cdots e_r \bigm| \forall k\colon e_k\in E_k\bigr\} 
\ar@{}[r]|-{\subset} & R_V \rlap{.} \\}$$
%of subsets of cardinality $1\cdot q\cdots q^{r-1} = q^{\frac{r(r-1)}{2}}$. 
Let $U$ be the subgroup of $\Aut_{\BF_q}(V)$ which sends each $X_k$ to an element of the coset $X_k+V_{k-1}$. In the given basis this corresponds to the subgroup of all upper triangular matrices in $\GL_r(\BF_q)$ with all diagonal entries~$1$. Then $U$ permutes each set $E_k$ transitively (but it does not act on~$\Delta_k$). Moreover, giving an element $g\in U$ is equivalent to giving the images $g(\frac{1}{X_k})\in E_k$ for all $k$, which can be chosen independently; hence $U$ acts freely transitively on $E_1\cdots E_r$.
It also follows that for each $1\le k\le r$ the element
\UseTheoremCounterForNextEquation
\begin{equation}\label{fkDef}
f_k\ := \sum_{e_k\in E_k} e_k 
\ = \sum_{w \in V_{k-1}} \frac{1}{X_k + w}  \ \in \ R_V.
\end{equation}
is fixed by~$U$. By \cite[Thm.\,2.7, Thm.\,2.11]{PinkSchieder} we have:

\begin{Thm}\label{RVBasis}
\begin{enumerate}
\item[(a)] The elements $f_1,\ldots,f_r$ are algebraically independent over $\BF_q$.
\item[(b)] The ring of $U$-invariants is $R_V^U = \BF_q[f_1,\ldots,f_r]$. 
\item[(c)] The ring $R_V$ is a free module over $R_V^U$ with basis $\Delta_1\cdots \Delta_r$. 
\end{enumerate}
\end{Thm}

%%%%%%%%%%%%%%%%%%%%%%%%%%%%%%%%%%%%%

We will also need the following fact:

\begin{Prop}\label{LocalizeByFk}
We have $RS_V=R_V[f_1^{-1},\ldots,f_r^{-1}]$.
\end{Prop}

\begin{Proof}
For any $1\le k\le r$ we apply Proposition \ref{RhoFormulas} (a) to $v:=X_k$ and $V':= V_{k-1}$ and the universal reciprocal map $\rho^\univ\colon \circV\to R_V$, $v\mapsto\frac{1}{v}$. By the definition (\ref{fkDef}) of $f_k$ we obtain that
\UseTheoremCounterForNextEquation
\begin{equation}\label{LocalizeByFk1}
f_k \cdot \kern-6pt\prod_{v'\in\circV_{k-1}}\kern-4pt\frac{1}{v'} \ =\ \Bigl(\sum_{v'\in V_{k-1}}\frac{1}{X_k-v'}\Bigr) \cdot \Bigl(\prod_{v'\in\circV_{k-1}}\frac{1}{v'}\Bigr) \ = 
\prod_{v'\in V_{k-1}}\frac{1}{X_k-v'}
\end{equation}
in $R_V$. Here all the factors except $f_k$ are already invertible in~$RS_V$; hence $f_k$ is also invertible in~$RS_V$. This proves the inclusion ``$\supset$''. 

For the inclusion ``$\subset$'' it suffices to show that for all $1\le k\le r$ and all $v'\in V_{k-1}$ the element $\frac{1}{x_k-v'}$ is invertible in $R_V[f_1^{-1},\ldots,f_r^{-1}]$. We will achieve this by induction on~$k$. Suppose that we already know it for all values smaller than~$k$. Then in particular the element $\frac{1}{v'}$ is invertible in $R_V[f_1^{-1},\ldots,f_r^{-1}]$ for any $v'\in\circV_{k-1}$. Thus the left hand side of (\ref{LocalizeByFk1}) is invertible in $R_V[f_1^{-1},\ldots,f_r^{-1}]$; hence so is the right hand side, and therefore also each factor on the right hand side, as desired.
\end{Proof}

%%%%%%%%%%%%%%%%%%%%%%%%%%%%%%%%%%%%%%%%%%%%%%%%%%%%%%%%%%%%%%%%%%%%%%

\subsection{A partial boundary}
\label{FqComplem}

We will need a variant of Theorem \ref{RVBasis} that concerns a partial boundary of $\Omega_V$ in~$Q_V$. 
% Strictly speaking I only need that the elements in Theorem \ref{RVJBasis} generate $F_V/J_s$ as $R'$-module, and not that they form a free basis. But I worked out the proof as a preparation for the analogous proof that I will have to do in Section \ref{Specialn}. I find it instructive for the reader to see the argument in a simpler case first.
Fix an arbitrary $\BF_q$-subspace $V'\subset V$. We assume that the basis of $V$ in the preceding subsection was chosen such that $V'=V_s$ for some $0\le s\le r$. Consider the following ideal of~$R_V$:
\UseTheoremCounterForNextEquation
\begin{equation}\label{JDef}
J_s\ :=\ \Bigl(\Bigl\{ \frac{1}{v'} \Bigm| v'\in\circV_s\; \Bigr\}\Bigr).
\end{equation} 
We want to give an explicit description of the factor ring $R_V/J_s$ 
(which happens to be the projective coordinate ring of the closed subscheme $X\subset Q_V$ used in Remark \ref{i*Morph2}).
For this we consider the submodule
\UseTheoremCounterForNextEquation
\begin{equation}\label{MsDef}
M_s\ := \bigoplus_{e\in\Delta_{s+1}\cdots\Delta_r} \kern-10pt
\BF_q[f_{s+1},\ldots,f_r]\cdot e\ \subset\ R_V,
\end{equation} 
which by Theorem \ref{RVBasis} is free with the indicated basis over $\BF_q[f_{s+1},\ldots,f_r]$.

\begin{Lem}\label{MsU}
The submodule $M_s$ is $U$-invariant and $M_s^U = \BF_q[f_{s+1},\ldots,f_r]$.
\end{Lem}

\begin{Proof}
First fix any $s+1\le k\le r$ and observe that $\BF_q[f_k] = \BF_q\oplus\BF_q[f_k]\cdot f_k$. Also note that the set $E_k$ from (\ref{EkDef}) is obtained from $\Delta_k$ on replacing the element $1\in\Delta_k$ by the element $\frac{1}{X_k} = f_k-\sum_{1\not=e_k\in\Delta_k}e_k$. Combining these facts we see that 
\begin{eqnarray*}
\bigoplus_{e_k\in\Delta_k}\BF_q[f_k] {\cdot} e_k
&=& \BF_q\oplus\BF_q[f_k] {\cdot} f_k \oplus
\bigoplus_{1\not=e_k\in\Delta_k}\kern-5pt \BF_q[f_k] {\cdot} e_k \\
&=& \BF_q\oplus\bigoplus_{e_k\in E_k}\BF_q[f_k] {\cdot} e_k.
\end{eqnarray*}
Taking the tensor product over $s+1\le k\le r$ and setting $E_I := \prod_{k\in I}E_k$, we deduce that
\begin{eqnarray*}
M_s &=& \bigoplus_{e\in\Delta_{s+1}\cdots\Delta_r} \kern-10pt
\BF_q[f_{s+1},\ldots,f_r]\cdot e \\
&=& \bigoplus_{I\subset\{s+1,\ldots,r\}}\;
\bigoplus_{e\in E_I} \BF_q[f_k|_{k\in I}]\cdot e.
\end{eqnarray*}
As $U$ permutes each $E_k$ and fixes each~$f_k$, it therefore acts on~$M_s$. Since $U$ acts transitively on $E_1\cdots E_r$, it also acts transitively on $E_I$ for each subset $I\subset\{s+1,\ldots,r\}$. The fact that $f_k=\sum_{e_k\in E_k}e_k$ and the above description of $M_s$ thus imply that 
\begin{eqnarray*}
M_s^U \!&=& \kern-5pt\bigoplus_{I\subset\{s+1,\ldots,r\}}\kern-10pt
\BF_q[f_k|_{k\in I}]\cdot 
\sum_{e\in E_I} e \\
&=& \kern-5pt\bigoplus_{I\subset\{s+1,\ldots,r\}}\kern-10pt
\BF_q[f_k|_{k\in I}]\cdot 
\prod_{k\in I} f_k 
\ \ =\ \ \BF_q[f_{s+1},\ldots,f_r],
\end{eqnarray*}
as desired.
\end{Proof}

\begin{Lem}\label{iV'Uinj}
Let $j\colon V''_s\into V$ denote the inclusion of the subspace that is generated by $X_{s+1},\ldots,X_r$. Then the associated homomorphism $\pi_j\colon R_V\onto R_{V''_s}$ from Proposition \ref{i*Homo} (b) restricts to an embedding 
$$\BF_q[f_{s+1},\ldots,f_r] \;\longinto\; R_{V''_s}.$$
\end{Lem}

\begin{Proof}
For each $s+1\le k\le r$ we have $f_k=\sum_{w\in V_{k-1}}\frac{1}{X_k+w}$, where $X_k+w\in V''_s$ if and only if $w\in V''_s$. Thus 
$$\pi_j(f_k)\ = \sum_{w\in V_{k-1}} \!\pi_j(\tfrac{1}{X_k+w})
\ = \kern-5pt \sum_{w\in V''_s\cap V_{k-1}} \kern-10pt \tfrac{1}{X_k+w}
\ \in\ R_{V''_s}.$$
The elements $\pi_j(f_{s+1}),\ldots,\pi_j(f_r)\in R_{V''_s}$ therefore play the same role for the space $V''_s$ with the basis $X_{s+1},\ldots,X_r$ as the elements $f_1,\ldots,f_r$ play for the space $V$ with the basis $X_1,\ldots,X_r$. By Theorem \ref{RVBasis} (a) they are therefore algebraically independent over~$\BF_q$. This means precisely that the restriction of $\pi_j$ to $\BF_q[f_{s+1},\ldots,f_r]$ is injective.
\end{Proof}

\begin{Thm}\label{RVJBasis}
The projection $\pi\colon R_V\onto R_V/J_s$ induces an isomorphism $M_s\isoto R_V/J_s$.
\end{Thm}

\begin{Proof}
%We use the same strategy as in the proof of Theorem \ref{IVBasis}.
Recall from Theorem \ref{RVBasis} (c) that $R_V$ is generated by $\Delta_1\cdots \Delta_r$ as a module over $\BF_q[f_1,\ldots,f_r]$. For any ${1\le k\le s}$, the definition of $J_s$ together with (\ref{DkDef}) shows that any element of $\Delta_k\setminus\{1\}$ lies in~$J_s$, and with (\ref{fkDef}) it shows that $f_k\in J_s$. Thus $R_V/J_s$ is already generated by $\pi(\Delta_{s+1}\cdots \Delta_r)$ as a module over $\BF_q[f_{s+1},\ldots,f_r]$. In other words the induced map $M_s\to R_V/J_s$ is surjective. 

It remains to show that this map is injective, or equivalently that its kernel $M_s\cap J_s$ is zero. For this observe that $U$ stabilizes~$\circVprime$ and therefore acts on~$J_s$. By Lemma \ref{MsU} it therefore also acts on $M_s\cap J_s$. Since $U$ is a finite group of $q$-power order acting on the $\BF_q$-vector space $M_s\cap J_s$, we have $M_s\cap J_s=0$ if and only if $M_s^U\cap J_s = (M_s\cap J_s)^U=0$. We are therefore reduced to showing that $\pi|M_s^U$ is injective.

But by the definition of $J_s$ the map $\pi_j\colon R_V\onto R_{V''_s}$ from Lemma \ref{iV'Uinj} factors through $\pi\colon R_V\onto R_V/J_s$.  Thus it suffices to show that $\pi_j|M_s^U$ is injective. This is now guaranteed by combining Lemmas \ref{MsU} and \ref{iV'Uinj}.
\end{Proof}

%\begin{Cor}\label{RVJCM}
%The factor ring $R_V/J_s$ is Cohen-Macaulay.
%\end{Cor}
%
%\begin{Proof}
%This follows from Theorem \ref{RVJBasis} exactly as in the proof of \cite[Thm.\,1.7]{PinkSchieder}.
%\end{Proof}

%%%%%%%%%%%%%%%%%%%%%%%%%%%%%%%%%%%%%%%%%%%%%%%%%%%%%%%%%%%%%%%%%%%%%%

\subsection{The ideal of the boundary}
\label{FqBoundaryIdeal}

For each $\BF_q$-subspace $V'\subset V$ consider the embedding map $i_{V'}\colon V'\into V$ and the associated ring homomorphism $\pi_{i_{V'}}\colon R_V \onto R_{V'}$ from \ref{i*Homo} (b). In this subsection we are interested in the ideal
\UseTheoremCounterForNextEquation
\begin{equation}\label{IVDef}
I_V\ := \bigcap_{0\not=V'\subsetneqq V} \kern-5pt\Ker(\pi_{i_{V'}}) \ \subset\ R_V.
\end{equation}
Since each $\pi_{i_{V'}}$ is a homomorphism from $R_V$ to an integral domain, this is a reduced ideal, and by construction it is graded. The associated closed subscheme of $Q_V$ is the reduced subscheme at the boundary
\UseTheoremCounterForNextEquation
\begin{equation}\label{PartialOmegaV}
\partial\Omega_V\ :=\ (Q_V\setminus\Omega_V)^\red 
\ =\ \bigcup_{0\not=V'\subsetneqq V} \kern-5pt \epsilon_{i_{V'}}(Q_{V'}).
\end{equation}
The following results give explicit generators for~$I_V$ in analogy to Theorem \ref{RVBasis} (c). We knew them at the time of writing \cite{PinkSchieder} and included them only as an exercise \cite[Ex.\,8.8]{PinkSchieder} in order to shorten 
%preserve the length of 
the paper.
%and because we did not see the full relevance. 
(Caution: The present ideal $I_V$ is not the ideal $\CI_V$ from \cite[\S6]{PinkSchieder}.)

\begin{Thm}\label{IVBasis}
\begin{enumerate}
\item[(a)] The ideal $I_V$ is a free module over $R_V^U$ with basis $E_1\cdots E_r$.
\item[(b)] The ideal $I_V$ is a free module over the group ring $\BF_q[U]$.
\end{enumerate}
\end{Thm}

\begin{Proof}
By \cite[Lemma 2.10]{PinkSchieder} the set $E_1\cdots E_r$ is the basis of a free $R_V^U$-submodule of~$R_V$. Denoting this submodule by~$M$, it follows that $M$ is a free module over $\BF_q[U]$. It remains to show that $M=I_V$.

For this note first that for any element $e_1\cdots e_r\in E_1\cdots E_r$, the reciprocals $e_1^{-1},\ldots,e_r^{-1}$ form a basis of $V$. Thus for any $\BF_q$-subspace $0\not=V'\subsetneqq V$; at least one of them lies in $V\setminus V'$. By the description of $\pi_{i_{V'}}$ in Proposition \ref{i*Homo} (b) it follows that $e_1\cdots e_r \in \Ker(\pi_{i_{V'}})$. Varying $V'$ this shows that $e_1\cdots e_r \in I_V$, and varying $e_1\cdots e_r$ then implies that $M\subset I_V$.

Next observe that $0\to M\to I_V\to I_V/M\to0$ is a short exact sequence of $\BF_q[U]$-modules. Since $M$ is a free $\BF_q[U]$-module, taking $U$-invariants yields a short exact sequence $0\to M^U\to I_V^U\to (I_V/M)^U\to H^1(U,M)=0$. Also, since $U$ is a finite group of $q$-power order acting on the $\BF_q$-vector space $I_V/M$, we have $I_V/M=0$ if and only if $(I_V/M)^U=0$. To prove that $M=I_V$, by the short exact sequence it is therefore enough to prove that $M^U=I_V^U$.

As the given basis $E_1\cdots E_r$ of $M$ over $R_V^U$ is a single free orbit under~$U$, the submodule $M^U$ is the free $R_V^U$-module generated by the element $\sum E_1\cdots E_r = f_1\cdots f_r$. In other words it is the principal ideal of $R_V^U$ generated by $f_1\cdots f_r$. Since $R_V^U = \BF_q[f_1,\ldots,f_r]$ with algebraically independent $f_1,\ldots,f_r$, this ideal is the intersection of the ideals $R_V^U\cdot f_k$ for all $1\le k\le r$. Thus it suffices to prove that $I_V^U\subset R_V^U\cdot f_k$ for every fixed $1\le k\le r$.

To achieve this consider the subspace $V'$ of codimension $1$ that is generated by all $X_j$ except~$X_k$. Then $\pi_{i_{V'}}(f_k)=0$, while for each $1\le j\le r$ with $j\not=k$ we have
$$\pi_{i_{V'}}(f_j)\ = \sum_{w \in V'\cap V_{j-1}} \frac{1}{X_j + w}  \ \in \ R_{V'}.$$
Thus by Theorem \ref{RVBasis} (a) for $V'$ in place of~$V$, the elements $\pi_{i_{V'}}(f_j)$ for $j\not=k$ are algebraically independent over~$\BF_q$, so the kernel of the homomorphism $\pi_{i_{V'}}|R_V^U$ is the ideal $R_V^U\cdot f_k$. By construction $I_V^U = I_V\cap R_V^U$ is contained in this kernel; hence we are done.
\end{Proof}

\begin{Cor}\label{IVGens}
The ideal $I_V$ is generated by the elements $\frac{1}{v_1\cdots v_r}$ for all bases $v_1,\ldots,v_r$ of~$V$. 
\end{Cor}

\begin{Proof}
For any basis $v_1,\ldots,v_r$ of~$V$, setting $X_i:= v_i$ shows that $\frac{1}{v_1\cdots v_r} \in I_V$ by Theorem \ref{IVBasis}. Conversely, Theorem \ref{IVBasis} shows that $I_V$ is generated by elements of this form.
\end{Proof}

\begin{Cor}\label{IVNotPrincipal}
The ideal sheaf of $\partial\Omega_V\subset Q_V$ is in general not locally principal.
\end{Cor}

\begin{Proof}
Consider any subspace $0\not=V'\subsetneqq V$ and choose a complement $V''\subset V$. Then by \cite[Prop.\,8.3]{PinkSchieder}, there exist open subschemes $U$ and $U'$ and an isomorphism between them making the commutative diagram:
$$\xymatrix{
\Omega_{V'} \ar@{^{ (}->}[r]^-{\epsilon_{i_{V'}}} \ar[d]^\wr
& U \rlap{$\ \,\subset\ Q_V$} \ar[d]^\wr \\
\Omega_{V'}\times\{0\} \ar@{^{ (}->}[r]
& U' \rlap{$\ \,\subset\  \Omega_{V'}\!\times\!\Spec(R_{V''}).$} \\
}\qquad\qquad$$
Moreover, the construction in 
%[loc.\ cit.] 
\cite{PinkSchieder} shows that this isomorphism induces an isomorphism
$$\xymatrix{
\Omega_V \cap U \ar@{^{ (}->}[r] \ar[d]^\wr
& U \rlap{$\ \,\subset\ Q_V$} \ar[d]^\wr \\
\ \ \ \ \llap{$(\Omega_{V'}\times\Spec(RS_{V''}))$} \cap U' \ar@{^{ (}->}[r]
& U' \rlap{$\ \,\subset\  \Omega_{V'}\!\times\!\Spec(R_{V''}).$} \\
}\qquad\qquad$$
Thus near the stratum $\epsilon_{i_{V'}}(\Omega_{V'})\subset Q_V$, the open embedding $\Omega_V\into Q_V$ is isomorphic to the open embedding $\Omega_{V'}\times\Spec(RS_{V''})\into\Omega_{V'}\times\Spec(R_{V''})$. This isomorphism identifies the ideal sheaf of the boundary $\partial\Omega_V\subset Q_V$ with the pullback from $\Spec(R_{V''})$ of the ideal sheaf associated to $I_{V''}\subset R_{V''}$. If $\dim_{\BF_q}(V'')\ge2$, Corollary \ref{IVGens} shows that the graded ideal $I_{V''}$ is not principal. Hence the associated ideal sheaf on $\Spec(R_{V''})$ is not locally principal at the apex $\{0\}$ of the cone $\Spec(R_{V''})$.
\end{Proof}

\newpage

%%%%%%%%%%%%%%%%%%%%%%%%%%%%%%%%%%%%%%%%%%%%%%%%%%%%%%%%%%%%%%%%%%%%%%
%%%%%%%%%%%%%%%%%%%%%%%%%%%%%%%%%%%%%%%%%%%%%%%%%%%%%%%%%%%%%%%%%%%%%%

\section{$A$-reciprocal maps}
\label{A}

%%%%%%%%%%%%%%%%%%%%%%%%%%%%%%%%%%%%%%%%%%%%%%%%%%%%%%%%%%%%%%%%%%%%%%

\subsection{Reminders on Drinfeld modules}
\label{Drin}

In this subsection we briefly recall various notions concerning Drinfeld modules (see for instance \cite[\S3]{PinkSatake}). 

Throughout we fix an admissible coefficient ring $A$ containing $\BF_q$ and set $F:=\Quot(A)$. The \emph{degree} of an element $a\in A$ is the number $\deg_A(a) := \dim_{\BF_q}(A/(a))$ if $a\not=0$, respectively $-\infty$ if $a=0$. 
Consider any commutative $F$-algebra~$R$. By a slight abuse of notation we denote the image in $R$ of an element $a\in A$ again by~$a$.
A \emph{standard Drinfeld $A$-module of rank $r\ge1$ over $R$} is an $\BF_q$-algebra homomorphism $\phi\colon A \to R[\tau]$, $a\mapsto\phi_a$ satisfying for every $a\in A\setminus\{0\}$:
\begin{enumerate}
\item[(a)] $\phi_a = \sum_{i=0}^{r\deg_A(a)}\phi_{a,i}\tau^i$ with $\phi_{a,r\deg_A(a)}\in R^\times$ , and
\item[(b)] $d\phi_a = \phi_{a,0} = a$.
\end{enumerate}
As $R$ is an $F$-algebra, condition (b) means that $\phi$ has \emph{generic characteristic}. 

%A \emph{standard Drinfeld $A$-module of rank $r\ge1$ over $R$} is an $\BF_q$-algebra homomorphism $\phi\colon A \to R[\tau]$, $a\mapsto\phi_a$ satisfying
%\begin{enumerate}
%\item[(a)] $\phi_a = \sum_{i=0}^{r\deg_A(a)}\phi_{a,i}\tau^i$ with $\phi_{a,r\deg_A(a)}\in R^\times$ for every $a\in A\setminus\{0\}$, and
%\item[(b)] $d\phi_a = \phi_{a,0} = a$ for every $a\in A$, in other words, that $\phi$ has \emph{generic characteristic}. 
%\end{enumerate}
More generally consider any scheme $S$ over~$F$. For any line bundle $E$ on $S$ let $\End_{\BF_q}(E)$ denote the ring of $\BF_q$-linear endomorphisms of the commutative group scheme underlying~$E$. Any trivialization of line bundles $\BG_{a,U}\isoto E|U$ over an open affine subscheme $U=\Spec(R)\subset S$ induces an isomorphism $R[\tau]\isoto\End_{\BF_q}(E|U)$. 
A general \emph{Drinfeld $A$-module of rank $r\ge1$ over $S$} is an $\BF_q$-algebra homomorphism $\phi\colon A \to \End_{\BF_q}(E)$, $a\mapsto\phi_a$ which for any trivialization of $E$ over an open affine subscheme becomes a standard Drinfeld $A$-module of rank~$r$.

For any such $(E,\phi)$ and any non-zero ideal $N\subset A$, the intersection $\phi[N] := \bigcap_{a\in N}\Ker(\phi_a)$ is an $A$-module subscheme of $E$ that is finite \'etale over $S$ and whose sections over any geo\-met\-ric point of $S$ form a free $A/N$-module of rank~$r$.
Consider the $A$-module $V^r_N := (N^{-1}/A)^{\oplus r}$, where $N^{-1}\subset F$ denotes the inverse fractional ideal of~$N$. 
This is also a free $A/N$-module of rank~$r$. 
A \emph{level $N$-structure} on $\phi$ is an $A$-module homomorphism $\lambda\colon V^r_N\to \phi[N](S)$ which induces an isomorphism in every fiber.
Observe that for any two non-zero ideals $N\subset N'\subset A$, we have a natural inclusion $V^r_{N'}\subset V^r_N$; hence any level $N$-structure restricts to a level $N'$-structure. In particular we can apply this when $N'=(a)$ for $a\in A\setminus\{0\}$, in which case we abbreviate $V^r_a := V^r_{(a)}$.

An \emph{isomorphism} of triples $(E,\phi,\lambda)\isoto(E',\phi',\lambda')$ as above is an isomorphism of line bundles that is compatible with $\phi$ and~$\lambda$. From now on we assume that $N$ is a proper ideal of~$A$. Then there exists at most one isomorphism $E\isoto E'$ that is compatible with the level $N$-structures $\lambda$ and~$\lambda'$, and so the isomorphism classes of such triples form a well-posed moduli problem. By \cite[\S5]{Drinfeld1} there is a \emph{fine moduli scheme} in the following sense:

\begin{Thm}\label{DrinModuli}
There is a scheme $M^r_{A,N}$ over $F$ and a triple $(E^\univ,\phi^\univ,\lambda^\univ)$ as above over $M^r_{A,N}$ such that:
\begin{itemize}
\item[(a)] For any scheme $S$ over $F$ and any triple $(E,\phi,\lambda)$ as above over~$S$, there exists a unique morphism $f\colon S\to M^r_{A,N}$ over~$F$ such that $(E,\phi,\lambda)\cong f^*(E^\univ,\phi^\univ,\lambda^\univ)$.
\item[(b)] This $M^r_{A,N}$ is an irreducible smooth affine algebraic variety of finite type and dimension $r-1$ over~$F$.
\end{itemize}
\end{Thm}

For a less canonical description of this moduli problem fix any $v_0\in\circVr_N$. Then any triple $(E,\phi,\lambda)$ as above determines an isomorphism of line bundles $\BG_{a,S}\isoto E$, $u\mapsto u\cdot\lambda(v_0)$. Giving an isomorphism class of triples $(E,\phi,\lambda)$ is therefore equivalent to giving a standard Drinfeld $A$-module $\phi'\colon A\to\CO_S(S)[\tau]$ of rank $r$ with a level $N$-structure $\lambda'\colon V^r_N\to \phi'[N](S)$ such that $\lambda'(v_0)=1$.

%%%%%%%%%%%%%%%%%%%%%%%%%%%%%%%%%%%%%%%%%%%%%%%%%%%%%%%%%%%%%%%%%%%%%%

\subsection{Conditions on the level}
\label{LevelConds}

Let $A$ be as in Subsection \ref{Drin}.
%Now we fix an admissible coefficient ring $A$ containing~$\BF_q$, and set $F:=\Quot(A)$. 
For any non-zero ideal $N\subset A$ consider the set of \emph{divisors of~$N$}
\UseTheoremCounterForNextEquation
\begin{equation}\label{DivDef}
\Div(N)\ :=\ \{a\in A\mid N\subset(a)\}.
\end{equation}
We are interested in a non-zero proper ideal $N\subset A$ satisfying the following conditions:

\begin{Ass}\label{NAss}
There exists a subset $D\subset\Div(N)$ with the properties:
\begin{enumerate}
%\item[(a)] $D$ generates $A$ as an $\BF_q$-algebra.
\item[(a)] For any element $a\in A\setminus\{0\}$ there exists an element $b\in A\setminus\{0\}$ which is a product of elements of $D$ such that $\deg_A(a-b)<\deg_A(a)$.
\item[(b)] For any distinct $a$, $b\in D$ we have $ab\in\Div(N)$.
\item[(c)] Any element of $\Div(N)$ is a product of elements in~$D$.
\item[(d)] $N$ is principal.
\end{enumerate}
\end{Ass}

\begin{Rem}\label{NAssRem1}
\rm Condition (a) implies that $\Div(N)$ generates $A$ as an $\BF_q$-algebra. One might hope that this consequence be enough to deduce all our results, but at present all conditions (a--d) are used in technical arguments. Condition (a) is used to show Proposition \ref{BarPhiWeakSep}, conditions (b) and (c) to prove Proposition \ref{RingHomo1}, and condition (d) for Proposition \ref{DrinMod}.
\end{Rem}

\begin{Rem}\label{NAssRem2}
\rm In the case $A=\BF_q[t]$ Assumption \ref{NAss} holds for any non-zero ideal $N\subset (t-\alpha)$ for any $\alpha\in\BF_q$, where $D$ consists of $\BF_q^\times$ and all monic irreducible divisors of~$N$.
\end{Rem}

%\begin{Prop}\label{NExists}
%\begin{enumerate}
%\item[(a)] For any non-zero ideal $N$ satisfying Assumption \ref{NAss}, any non-zero ideal $N'\subset N$ also satisfies Assumption \ref{NAss}.
%\item[(b)] For any maximal ideal $\Fp\subset A$ there exists an ideal $N\not\subset\Fp$ satisfying Assumption \ref{NAss}.
%\end{enumerate}
%\end{Prop}
%
%\begin{Proof}
%Assertion (a) follows by taking the same set~$D$. To prove (b) pick any elements $a_1,\ldots,a_n\in A$ which generate $A$ as an $\BF_q$-algebra. After possibly replacing some $a_i$ by $a_i+1$, we may assume that all $a_i\not\in\Fp$. The ideal $N:= (a_1\cdots a_n)^2$ then has the desired properties with $D := \{a_1,\ldots,a_n\}$.
%\end{Proof}

%CAUTION: In general we cannot achieve that $N+\Fa=A$ for an arbitrary non-zero ideal~$\Fa$. For instance, if $\Fa=\Fp_1\Fp_2$ for distinct maximal ideals $\Fp_i$ with residue fields $\BF_2$, all divisors of $N$ must be congruent to $1$ modulo both $\Fp_i$ and hence modulo $\Fa$; so the $\BF_q$-subalgebra generated by them is contained in $\BF_q+\Fa \subsetneqq A$.

\begin{Prop}\label{NExists}
For any maximal ideal $\Fp$ and any ideal $\Fa$ not contained in~$\Fp$, there exists an ideal $N$ contained in~$\Fa$, but not in~$\Fp$, which satisfies Assumption \ref{NAss}.
\end{Prop}

\begin{Proof}
For any integer $d\ge0$ abbreviate $A_{\le d} := \{a\in A\mid \deg_A(a)\le d\}$. 
%For any integer $d\ge0$ abbreviate $\circAbig_d := \{a\in A\setminus\{0\}\mid \deg_A(a)\le d\}$. 
By Riemann-Roch there exists an integer $d>0$ such that for every element $a\in A\setminus A_{\le d}$ there exists an element $b\in A$ which is a product of elements of $A_{\le d}$ such that $\deg_A(a-b)<\deg_A(a)$. Note that any element $c\in A_{\le d}\cap\Fp$ that is a factor of~$b$ is non-constant, so we can replace $c$ by $c+1$ in the product defining $b$ without destroying the property $\deg_A(a-b)<\deg_A(a)$. Thus the set $A_{\le d}\setminus\Fp$ has the property of $D$ in (a), and all its elements are non-zero.
%Let $a_1,\ldots,a_n$ be all elements of $A\setminus\{0\}$ of degree $\le d$. 

Choose an element $a_0\in\Fa\setminus\Fp$, and let $\Fb\subset A$ be the principal ideal generated by $a_0\cdot\prod_{a\in A_{\le d}\setminus\Fp}a$. Then $\Fb$ is contained in~$\Fa$, but not in~$\Fp$. Let $h$ be the class number of~$A$. We claim that the ideal $N:=\Fb^{2h}$ has the desired properties. 

Indeed, by construction $N$ is contained in~$\Fa$, but not in~$\Fp$. 
%Being a power with exponent~$h$, it is a principal ideal. 
Also, since $\Fb$ is principal, so is~$N$; hence condition (d) holds.

Next let $\Fp_1,\ldots,\Fp_m$ be the distinct prime factors of~$\Fb$. Let $D$ be the set of all $a\in A\setminus\{0\}$ for which $(a)=\Fp_1^{e_1}\cdots\Fp_m^{e_m}$ with all $e_i\le h$. Then, in particular, for each $i$ there exists an element $b_i\in D$ with $(b_i)=\Fp_i^h$. It follows that every non-zero element of $A$ whose prime ideal factorization contains no primes other than $\Fp_1,\ldots,\Fp_m$ is a product of elements of~$D$. In particular this therefore holds for all divisors of~$N$, establishing condition (c). 

Moreover, any element of $D$ is a divisor of~$\Fb^h$; hence any product of two elements of $D$ is a divisor of~$N$, proving condition (b). 

Finally, every element of $A_{\le d}\setminus\Fp$ is a divisor of $N$ and therefore a product of elements of~$D$. Thus $D$ satisfies the condition (a), and we are done.
\end{Proof}

%%%%%%%%%%%%%%%%%%%%%%%%%%%%%%%%%%%%%

\medskip
For our later purposes, Proposition \ref{NExists} guarantees that there are sufficiently many non-zero ideals $N$ satisfying Assumption \ref{NAss}.

%%%%%%%%%%%%%%%%%%%%%%%%%%%%%%%%%%%%%%%%%%%%%%%%%%%%%%%%%%%%%%%%%%%%%%

\subsection{Basic $A$-reciprocal maps}
\label{BasicARecip}

For the following we fix an ideal $N\subset A$ satisfying Assumption \ref{NAss}. 
%Let $A_N$ denote the subring of all elements of $F:=\Quot(A)$ that are regular outside infinity and the divisors of~$N$.
We also fix an integer $r>0$ and define $V^r_N := (N^{-1}/A)^{\oplus r}$ and $V^r_a := (a^{-1}A/A)^{\oplus r}$ as in Subsection \ref{Drin}. Then for any $a\in\Div(N)$ we have $V^r_a\subset V^r_N$.
Consider any $F$-algebra~$R$. 

\begin{Def}\label{ARecipDef}
A map $\rho\colon\circVr_N\to R$ is called \emph{$A$-reciprocal} if 
\begin{itemize}
\item[(a)] $\rho(v) \cdot \rho(w) = \rho(v+w) \cdot (\rho(v) + \rho(w))$ for all $v$, $w\in \circVr_N$ with $v+w\in\circVr_N$, and
\item[(b)] $a\rho(av) = \sum_{v'\in V^r_a}\rho(v-v')$ for all $a\in\Div(N)$ and  $v \in V^r_N\setminus V^r_a$.
\end{itemize}
\end{Def}

\begin{Prop}\label{ARecipIsRecip}
Any $A$-reciprocal map is $\BF_q$-reciprocal.
\end{Prop}

\begin{Proof}
Every $\alpha\in\BF_q^\times$ lies in $\Div(N)$ and satisfies $V^r_\alpha=\{0\}$.
\end{Proof}

\begin{Cons}\label{RAVCons}
\rm Let $R_{A,V^r_N}$ denote the factor ring of $R_{V^r_N}\otimes_{\BF_q}F$ modulo the ideal 
$$\left(\left\{ \left. 
\frac{1}{av}\otimes a-\sum_{v'\in V^r_a}\frac{1}{v-v'}\otimes1
\ \right| \begin{array}{l}
\hbox{all $a\in\Div(N)$,}\\[3pt]
\hbox{all $v \in V^r_N\setminus V^r_a$}
\end{array}
\right\}\right).$$
For any element $f\in R_{V^r_N}\otimes_{\BF_q}F$ let $[f]$ denote its image in $R_{A,V^r_N}$. Let $RS_{A,V^r_N}$ be the localization of $R_{A,V^r_N}$ obtained by inverting the elements $[\frac{1}{v}\otimes1]$ for all $v\in\circVr_N$. 
\end{Cons}

\begin{Thm}\label{AModuli}
\begin{itemize}
\item[(a)] The map 
$$\rho^\univ\colon \circVr_N \longto R_{A,V^r_N},
\ v\mapsto [\tfrac{1}{v}\otimes1]$$
is $A$-reciprocal.
\item[(b)] For any $F$-algebra $R$ and any $A$-reciprocal map $\rho\colon\circVr_N\to R$ there exists a unique $F$-algebra homomorphism $f\colon R_{A,V^r_N}\to R$ such that $\rho=f\circ\rho^\univ$.
\item[(c)] This $f$ extends to a ring homomorphism $RS_{A,V^r_N}\to R$ if and only if $\rho$ is fiberwise invertible.
\end{itemize}
\end{Thm}

\begin{Proof}
Theorem \ref{FqModuli} (a) implies that $\rho^\univ$ is $\BF_q$-reciprocal. It is $A$-reciprocal, because in the construction of $R_{A,V^r_N}$ we have divided out precisely the relations corresponding to \ref{ARecipDef} (b) that make an $\BF_q$-reciprocal map an $A$-reciprocal map. This proves (a).

In the situation of (b), the $\BF_q$-reciprocal map underlying $\rho$ already corresponds to a unique $\BF_q$-algebra homomorphism $f'\colon R_{V^r_N}\to R$ such that $\rho=f'\circ\rho^\univ$ by Theorem \ref{FqModuli} (b). As $R$ is an $F$-algebra, this $f'$ corresponds to a unique $F$-algebra homomorphism $f''\colon R_{V^r_N}\otimes_{\BF_q}F\to R$, which sends $\frac{1}{v}\otimes1$ to $\rho(v)$ for all $v\in\circVr_N$. The condition \ref{ARecipDef} (b) for $\rho$ then implies that the ideal in Construction \ref{RAVCons} lies in the kernel of~$f''$. Thus $f''$ factors through the desired~$f$, proving (b).

Assertion (c) then follows directly from the construction of $RS_{A,V^r_N}$.
\end{Proof}

%%%%%%%%%%%%%%%%%%%%%%%%%%%%%%%%%%%%%

\medskip
Again we finish this subsection by addressing functoriality:

\begin{Prop}\label{AExtZero}
Consider any $1\le s\le r$ and any injective $A$-linear map $i\colon V^s_N\into V^r_N$. Then for any $A$-reciprocal map $\rho\colon\circVs_N\to R$ the extension by zero $i_*\rho\colon \circVr_N\to R$ from Proposition \ref{ExtByZero} is $A$-reciprocal.
%following map is $A$-reciprocal:
%$$i_*\rho\!:\ \circVr_N\to R, \ v \mapsto 
%\left\{\begin{array}{cl}
%\rho(w) & \hbox{if $v=i(w)$ for some $w\in\circVs_N$,} \\[3pt]
%0 & \hbox{otherwise.}
%\end{array}\right.$$
\end{Prop}

\begin{Proof}
By Proposition \ref{ExtByZero} the map $i_*\rho$ is already $\BF_q$-reciprocal. 
%Consider any elements $v$, $w\in \circVr_N$ such that $v+w\in\circVr_N$. If all of $v$, $w$, $v+w$ lie in the image of~$i$, the equation $i_*\rho(v) \cdot i_*\rho(w) = i_*\rho(v+w) \cdot (i_*\rho(v) + i_*\rho(w))$ holds 
%because the corresponding equation holds for~$\rho$. Otherwise at least two of $v$, $w$, $v+w$ lie outside the image of~$i$; hence at least two of the values $i_*\rho(v)$, $i_*\rho(w)$, $i_*\rho(v+w)$ are zero, and the same equation follows trivially. Thus $i_*\rho$ satisfies the condition \ref{ARecipDef} (a).
Next take any $a\in\Div(N)$ and any $v \in V^r_N\setminus V^r_a$. 
Suppose first that $v=i(w)$ for some $w\in\circVs_N$. Then $w\not\in V^s_a$. Moreover, for any $v'\in V^r_a$ we have $v-v'\in i(V^s_N)$ if and only if $v'\in i(V^s_a)$. By the condition \ref{ARecipDef} (b) for $\rho$ and the definition of $i_*\rho$ we deduce that 
$$a\cdot i_*\rho(av)\ =\ a\rho(aw)
\ =\ \sum_{w'\in V^s_a}\rho(w-w')
\ =\ \sum_{v'\in V^r_a}i_*\rho(v-v').$$
Next we observe that the far left and right sides of this equation depend only on the coset $v+V^r_a$. Thus the total equation also holds if $v\in i(V^s_N)+V^r_a$. 
Finally suppose that $v\not\in i(V^s_N)+V^r_a$.
Then for all $v'\in V^r_a$ we have $v-v'\not\in i(V^s_N)$ and hence $i_*\rho(v-v')=0$. On the other hand, since $i(V^s_N)$ is a direct summand of $V^r_N$ as an $A$-module, the assumption $v\not\in i(V^s_N)+V^r_a$ implies that $av\not\in i(V^s_N)$. Thus $i_*\rho(av)=0$ as well, and the total equation holds trivially in this case.
Together this proves that $i_*\rho$ satisfies the condition \ref{ARecipDef} (b). Therefore $i_*\rho$ is $A$-reciprocal.
\end{Proof}

\begin{Prop}\label{i*HomoA}
\begin{itemize}
\item[(a)] The functor $i_*$ on $A$-reciprocal maps is represented by a surjective $F$-algebra homomorphism $\pi_i\colon {R_{A,V^r_N}\onto R_{A,V^s_N}}$ that sends $[\frac{1}{i(v')}\otimes1]$ to $[\frac{1}{v'}\otimes1]$ for all $v'\in\circVs_N$ and $[\frac{1}{v}\otimes1]$ to $0$ for all $v\in V^r_N\setminus i(V^s_N)$.
\item[(b)] The kernel of $\pi_i$ is generated by the elements $[\frac{1}{v}\otimes1]$ for all $v\in V^r_N\setminus i(V^s_N)$.
\end{itemize}
\end{Prop}

\begin{Proof}
Precisely analogous to the proof of Proposition \ref{i*Homo} (b) and (c).
\end{Proof}

%%%%%%%%%%%%%%%%%%%%%%%%%%%%%%%%%%%%%%%%%%%%%%%%%%%%%%%%%%%%%%%%%%%%%%

\subsection{The associated ring homomorphism}
\label{AssRingHom}

Keeping the notation of the preceding subsection, we now fix an $A$-reciprocal map $\rho\colon \allowbreak {\circVr_N\to R}$. For any $a\in\Div(N)$ consider the polynomial 
\UseTheoremCounterForNextEquation
\begin{equation}\label{Phia}
\phi_a(X)\ :=\ a\cdot X\cdot\prod_{v\in\circVr_a}\bigl(1-\rho(v)X\bigr)\ \in\ R[X].
\end{equation}
Comparison with (\ref{ExpDef}) shows that $\phi_a = a\circ e_{\rho|\circVr_a}$; hence $\phi_a\in R[\tau]$ by Proposition \ref{ExpFqLinear}.

\begin{Lem}\label{Phiaa'Lem}
For any $a$, $b\in\Div(N)$ with $ab\in\Div(N)$ we have $\phi_a\circ\phi_b=\phi_{ab}$.
\end{Lem}

\begin{Proof}
By assumption we have a short exact sequence 
$$\xymatrix{0\ar[r]&V^r_b\ar[r]^i&V^r_{ab}\ar[r]^p&V^r_a\ar[r]&0\rlap{,}\\}$$
where $i$ denotes the inclusion $V^r_b\into V^r_{ab}$ and $p$ denotes multiplication by~$b$. By the definition \ref{Quotient} of $p_*$ and the condition \ref{ARecipDef} (b), for every $v\in V^r_{ab}\setminus V^r_b$ we have 
$$p_*(\rho|\circVr_{ab})(bv)
\ =\ \sum_{v'\in V^r_b} \rho(v-v')
\ =\ b\rho(bv).$$
In other words we have $p_*(\rho|\circVr_{ab}) = b\cdot\rho|\circVr_a$. Using Proposition \ref{Composition} we find that
\begin{align*}
\phi_a\circ\phi_b
\ &=\ a\circ e_{\rho|\circVr_a} \circ b\circ e_{\rho|\circVr_b} \\
\ &=\ a\circ b\circ e_{b\cdot\rho|\circVr_a} \circ e_{\rho|\circVr_b} \\
\ &=\ ab\circ e_{p_*(\rho|\circVr_{ab})} \circ e_{i^*(\rho|\circVr_{ab})} \\
\ &=\ ab\circ e_{\rho|\circVr_{ab}} \\
\ &=\ \phi_{ab},
\end{align*}
as desired.
\end{Proof}

\begin{Lem}\label{CommLem1}
Any element of $R[\tau]\tau$ which commutes with $\phi_a$ for some non-constant $a\in\Div(N)$ is zero.
\end{Lem}

\begin{Proof}
Suppose that there is a non-zero $\eta\in R[\tau]\tau$ which commutes with~$\phi_a$. Write $\eta=u\tau^i+($higher terms in $\tau)$ with $u\in R\setminus\{0\}$ and $i\ge1$. By the construction (\ref{Phia}) we have $\phi_a=a+($higher terms in $\tau)$. Thus
\begin{align*}
0\ &=\ \eta\circ\phi_a-\phi_a\circ\eta \\
\ &=\ u\tau^i\circ a - a\circ u\tau^i +(\hbox{higher terms in $\tau$}) \\
\ &=\ u(a^{q^i}-a)\tau^i+(\hbox{higher terms in $\tau$}),
\end{align*}
and hence $u(a^{q^i}-a)=0$. But since $a$ is a transcendental element of~$F$ and $i\ge1$, we have $a^{q^i}-a\in F^\times$. Thus we conclude that $u=0$, contrary to the assumption.
%
%CAUTION: This argument seems to fail if $R$ is only an $A_N$-algebra. This is why I consider only the case of generic characteristic.
\end{Proof}

%\begin{Lem}\label{CommLem2}
%For any $a$, $b\in\Div(t^{-1}N)$ we have $\phi_a\circ\phi_t=\phi_t\circ\phi_a$ and $\phi_a\circ\phi_b=\phi_b\circ\phi_a$.
%\end{Lem}
%
%\begin{Proof}
%The first equation follows from Lemma \ref{Phiaa'Lem} by the computation $\phi_a\circ\phi_t=\phi_{at}=\phi_{ta}=\phi_t\circ\phi_a$. Thus $\phi_a$, and by the same argument also $\phi_b$, commutes with~$\phi_t$. Therefore $\eta := \phi_a\circ\phi_b-\phi_b\circ\phi_a$ also commutes with~$\phi_t$. But the commutator of any two elements of $R[\tau]$ lies in $R[\tau]\tau$. Thus by Lemma \ref{CommLem1} we deduce that $\eta=0$, which proves the second equation.
%\end{Proof}

\begin{Prop}\label{RingHomo1}
The map $\Div(N)\to R[\tau]$, $a\mapsto\phi_a$ extends to a unique $\BF_q$-algebra homo\-mor\-phism $\phi\colon A\to R[\tau]$, $a\mapsto\phi_a$, which satisfies $d\phi_a=a$ for all $a\in A$.
\end{Prop}

\begin{Proof}
Let $D=\{a_1,\ldots,a_n\}$ be the subset of $\Div(N)$ from Assumption \ref{NAss}. Since $D$ generates $A$ as an $\BF_q$-algebra, we have a surjective $\BF_q$-algebra homomorphism 
$$\Pi\colon \BF_q[Y_1,\ldots,Y_n]\onto A, \ P\mapsto P(a_1,\ldots,a_n).$$
Next, for any two distinct elements $a_i$, $a_j$ of $D$ we have $a_ia_j\in\Div(N)$ by assumption. Thus by Lemma \ref{Phiaa'Lem} we have $\phi_{a_i}\circ\phi_{a_j}=\phi_{a_ia_j}=\phi_{a_ja_i}=\phi_{a_j}\circ\phi_{a_i}$. In other words, the elements $\phi_{a_i}\in R[\tau]$ all commute with each other; hence we also have an $\BF_q$-algebra homomorphism 
$$\Phi\colon \BF_q[Y_1,\ldots,Y_n]\to R[\tau], \ P\mapsto P(\phi_{a_1},\ldots,\phi_{a_n}).$$
Since each $\phi_{a_i} = a_i+($higher terms in $\tau)$, for any polynomial $P\in\BF_q[Y_1,\ldots,Y_n]$ we have $P(\phi_{a_1},\ldots,\phi_{a_n}) \allowbreak = P(a_1,\ldots,a_n)+($higher terms in $\tau)$. In particular, for any $P\in\Ker(\Pi)$ we have $P(\phi_{a_1},\ldots,\phi_{a_n}) \in R[\tau]\tau$. But at least one $a_i$ is non-constant, and since the associated $\phi_{a_i}$ commutes with all $\phi_{a_j}$, it commutes with $P(\phi_{a_1},\ldots,\phi_{a_n})$. Using Lemma \ref{CommLem1} we therefore deduce that $P(\phi_{a_1},\ldots,\phi_{a_n})=0$. Varying $P$ this shows that $\Ker(\Pi) \subset \Ker(\Phi)$, which imples that $\Phi=\psi\circ\Pi$ for a unique $\BF_q$-algebra homomorphism $\psi\colon A\to R[\tau]$.

By construction this algebra homomorphism satisfies $\psi_a=\phi_a$ for all $a\in D$. By Assumption \ref{NAss} (c) and Lemma \ref{Phiaa'Lem} the same equality then follows for all $a\in\Div(N)$. This proves that the desired extension $\phi$ exists. The uniqueness follows from the fact that $D$ generates~$A$ as an $\BF_q$-algebra.

Finally, by (\ref{Phia}) the formula $d\phi_a=a$ holds for all $a\in\Div(N)$. Since $\phi$ is an $\BF_q$-algebra homomorphism, the same follows for all $a\in A$.
\end{Proof}

%%%%%%%%%%%%%%%%%%%%%%%%%%%%%%%%%%%%%%%%%%%%%%%%%%%%%%%%%%%%%%%%%%%%%%

\subsection{Constant rank}
\label{ConstantRank}

In this subsection we consider an $A$-reciprocal map $\rho\colon \circVr_N\to R$ satisfying the condition
\UseTheoremCounterForNextEquation
\begin{equation}\label{BB}
\forall v\in\circVr_N\colon\ \rho(v)=0\ \hbox{or}\ \rho(v)\in R^\times.
\end{equation}

\begin{Lem}\label{WASubMod}
The subset
$$W\ :=\ \{0\}\cup\{v\in\circVr_N\mid \rho(v)\in R^\times\}$$
is an $A$-submodule of~$V^r_N$.
\end{Lem}

\begin{Proof}
%(Compare \cite[Prop.\,7.7 (b)]{PinkSchieder}.)
Consider any $v$, $w\in W$. If one or more of $v,w,v+w$ is zero, we directly see that $v+w\in \{v,w,0\}\subset W$. Otherwise we have $\rho(v),\rho(w)\in R^\times$ by construction. Thus by Definition \ref{ARecipDef} (a) we have $\rho(v+w) \cdot (\rho(v) + \rho(w)) = \rho(v) \cdot \rho(w) \in R^\times$; hence $\rho(v+w)\in R^\times$ and so $v+w\in W$. Together this shows that $W+W\subset W$.

Next consider any $a\in\Div(N)$ and any $v\in W$. If $v\in V^r_a$, we have $av=0\in W$. Otherwise we have $\rho(v)\in R^\times$ by construction, and so by Definition \ref{ARecipDef} (b) and Proposition \ref{RhoFormulas} (b) we deduce that 
$$a\rho(av) \cdot \prod_{v'\in\circVr_a}\bigl(\rho(v)-\rho(v')\bigr) 
%\ =\ \Bigl(\sum_{v'\in V^r_a}\rho(v-v')\Bigr) \cdot \Bigl(\prod_{v'\in\circVr_a}\bigl(\rho(v)-\rho(v')\bigr)\Bigr) 
\ =\ \rho(v)^{|V^r_a|}
\ \in\ R^\times.$$
Thus $\rho(av)\in R^\times$ and therefore $av\in W$. This shows that $aW\subset W$ for all $a\in\Div(N)$.

Since $\BF_q^\times\subset\Div(N)$, this implies that $W$ is an $\BF_q$-subspace of~$V^r_N$. As $\Div(N)$ generates $A$ as an $\BF_q$-algebra, it is then also an $A$-submodule.
\end{Proof}

\begin{Lem}\label{LambdaLinear}
The map
$$\lambda\colon W\to R,\ v\mapsto
\left\{\begin{array}{cl}
0 & \hbox{if $v=0$,}\\[3pt]
\rho(v)^{-1} & \hbox{if $v\not=0$,}\\
\end{array}\right.$$
is additive and satisfies $\lambda(av)=\phi_a(\lambda(v))$ for all $v\in W$ and $a\in A$.
\end{Lem}

\begin{Proof}
Since $\rho$ is $\BF_q$-reciprocal by Proposition \ref{ARecipIsRecip}, Proposition \ref{FqLinearMap} implies that $\lambda$ is $\BF_q$-linear. In particular it is additive.

Next consider any $a\in\Div(N)$ and any $v\in W$. If $v=0$, we also have $av=0$ and hence $\lambda(av)=0=\phi_a(0)=\phi_a(\lambda(v))$. 
If $v\in V^r_a\setminus\{0\}$, we still have $av=0$ and hence $\lambda(av)=0$. But from (\ref{Phia}) we then obtain that
$$\phi_a(\lambda(v))\ =\ a\cdot \lambda(v)\cdot\prod_{v'\in\circVr_a}(1-\rho(v')\lambda(v)),$$
where the factor $1-\rho(v)\lambda(v)$ associated to $v'=v$ is zero. Thus again we find that $\lambda(av)=0=\phi_a(\lambda(v))$. 
Suppose now that $v\not\in V^r_a$. Then by combining the formulas in (\ref{Phia}) and Definition \ref{ARecipDef} (b) and Proposition \ref{RhoFormulas} (c) we deduce that 
$$\rho(av) \cdot \phi_a(\lambda(v))
\ =\ a\rho(av) \cdot e_{\rho|\circVr_a}(\lambda(v))\ =\ 1.$$
Thus again we find that $\lambda(av)=\phi_a(\lambda(v))$. 

As this formula holds for all $a\in\Div(N)$, which generate $A$ as an $\BF_q$-algebra, and each $\phi_a$ is also $\BF_q$-linear, the formula then follows for all $a\in A$.
\end{Proof}

\begin{Prop}\label{DrinMod}
\begin{itemize}
\item[(a)] If $W$ is zero, we have $\phi_a=a$ in $R[\tau]$ for all $a\in A$.
\item[(b)] If $W$ is non-zero, it is a free $A/N$-module of some rank ${1\le s\le r}$, and $\phi$ is a standard Drinfeld $A$-module of rank~$s$, and for any isomorphism $i\colon V^s_N\isoto W$ the map $\lambda\circ i$ is a level $N$-structure of~$\phi$.
\end{itemize}
\end{Prop}

\begin{Proof}
Assume first that $W=0$. Then for any $a\in\Div(N)$ the definition (\ref{Phia}) of $\phi_a$ shows that $\phi_a(X)=aX$ and hence $\phi_a=a$ in $R[\tau]$. As $\Div(N)$ generates $A$ as an $\BF_q$-algebra, the equality $\phi_a=a$ then holds for all $a\in A$.

Assume now that $W\not=0$. The desired assertions hold trivially if $R=0$, so we may also assume that $R\not=0$. 
%Then by construction we have $\Ker(\lambda)=0$; hence $\lambda$ is injective. Next by 
By Assumption \ref{NAss} (d) we have $N=(a)$ for some non-constant element $a\in\Div(A)$. The definition (\ref{Phia}) of $\phi_a$ and the assumption (\ref{BB}) imply that 
$$\phi_a(X)\ =\ a\cdot X\cdot\prod_{v\in\circW}\left(1-\frac{X}{\lambda(v)}\right)$$
with highest non-zero coefficient in~$R^\times$. 
%On the other hand we have $aW=0$ and hence $\phi_a(\lambda(W))=0$ by Lemma \ref{LambdaLinear}. Since $W\not=0$ and $\lambda$ is injective, we cannot have $\phi_a=a$; hence $\phi_a\not\in R$. 
As $W\not=0$, we have $\phi_a\not\in R$. 
By general theory (for instance \cite[\S2]{Drinfeld1}) it thus follows that $\phi$ is a Drinfeld $A$-module of some constant rank $s\ge1$ and that $W$ is a free $A/N$-module of rank~$s$. Since $W\subset V^r_N\cong(A/N)^{\oplus r}$, this implies that $s\le r$. Finally, since $\lambda\circ i$ is $A$-linear and injective, it is a level $N$-structure of~$\phi$.
\end{Proof}

%%%%%%%%%%%%%%%%%%%%%%%%%%%%%%%%%%%%%

\medskip
We also have the following converse of Proposition \ref{DrinMod}:

\begin{Prop}\label{DrinModConv}
Let $\phi'\colon A\to R[\tau]$ be a standard Drinfeld $A$-module of rank ${1\le s\le r}$ with a level $N$-structure $\lambda'\colon V^s_N\to R$, and let $i\colon V^s_N\into V^r_N$ be any injective $A$-linear map. Then the map
$$\rho\colon \circVr_N\to R,\ v\mapsto
\left\{\begin{array}{cl}
\lambda'(w)^{-1} & \hbox{if $v=i(w)$ for $w\in V^s_N$,}\\[3pt]
0 & \hbox{if $v\not\in i(V^s_N)$,}
\end{array}\right.$$
is $A$-reciprocal and satisfies condition (\ref{BB}), and we have $(W,\lambda\circ i,\phi) = (i(V^s_r),\lambda',\phi')$.
\end{Prop}

\begin{Proof}
Suppose first that $s=r$. Then $\rho$ is $\BF_q$-reciprocal by Proposition \ref{FqLinearMap}.
Next take any $a\in\Div(N)$. The fact that $\lambda'$ is a level $N$-structure of~$\phi'$ and the assumption $d\phi'_a=a$ imply that 
\UseTheoremCounterForNextEquation
\begin{equation}\label{DrinModConv1}
\phi'_a(X)
\ =\ a\cdot X\cdot\kern-3pt\prod_{v'\in \circVr_a}\kern-3pt\left(1-\frac{X}{\lambda'(v')}\right)
\ =\ a\cdot X\cdot\kern-3pt\prod_{v'\in \circVr_a}\kern-3pt\bigl(1-\rho(v')X\bigr).
\end{equation}
For any $v \in V^r_N\setminus V^r_a$, the $A$-linearity of $\lambda'$ combined with (\ref{DrinModConv1}) shows that 
$$\lambda'(av)\ =\ \phi'_a(\lambda'(v))
\ =\ a\cdot\lambda'(v)\cdot\kern-3pt\prod_{v'\in\circVr_a}\kern-3pt\bigl(1-\rho(v')\lambda'(v)\bigr).$$
Applying Proposition \ref{RhoFormulas} (c) to the $\BF_q$-reciprocal map $\rho$ and the subspace $V^r_a$ thus implies that 
$$\Bigl(\, \sum_{v'\in V^r_a} \rho(v-v') \Bigr) \cdot \frac{\lambda'(av)}{a}\ =\ 1$$
and hence 
$$a\cdot\rho(av)\ = \sum_{v'\in V^r_a} \rho(v-v').$$
Thus $\rho$ satisfies the condition \ref{ARecipDef} (b) and is therefore $A$-reciprocal. 
By construction $\rho$ is fiberwise invertible, so it satisfies condition (\ref{BB}), and by Lemmas \ref{WASubMod} and \ref{LambdaLinear} we have $W=V^r_N$ and $\lambda=\lambda'$. Lastly, for all $a\in\Div(N)$ we have $\phi_a=\phi'_a$ by (\ref{Phia}) and (\ref{DrinModConv1}). As $\Div(N)$ generates $A$ as an $\BF_q$-algebra, this implies that $\phi_a=\phi'_a$ for all $a\in A$, and we are done.

In the general case the same arguments with $s$ in place of $r$ yield a fiberwise invertible $A$-reciprocal map $\rho'\colon \circVs_N\to R$ which returns $\lambda'$ and~$\phi'$. The map $\rho$ in question is then simply the extension by zero of $\rho'$ via $i$; hence it is $A$-reciprocal by Proposition \ref{AExtZero}. By construction this $\rho$ satisfies condition (\ref{BB}), and by Lemmas \ref{WASubMod} and \ref{LambdaLinear} we have $W=i(V^s_N)$ and $\lambda\circ i=\lambda'$. Finally, the formula (\ref{Phia}) shows that $\phi_a$ does not change under extension by zero for $a\in\Div(N)$. As $\Div(N)$ generates $A$ as an $\BF_q$-algebra, this implies that $\phi_a=\phi'_a$ for all $a\in A$, and we are done.
\end{Proof}

\subsection{General $A$-reciprocal maps}
\label{GenARecip}

We keep $A$ and $F=\Quot(A)$ and $N$ as in Subsection \ref{BasicARecip}. Let $S$ be a scheme over $F$ and $\CL$ an invertible sheaf on $S$.

\begin{Def}\label{GenARecipDef}
A map $\rho\colon\circVr_N\to\CL(S)$ is called \emph{$A$-reciprocal} if 
\begin{itemize}
\item[(a)] $\rho(v) \cdot \rho(w) = \rho(v+w) \cdot (\rho(v) + \rho(w))$ in $\CL^{\otimes2}(S)$ for all $v$, $w\in \circVr_N$ with $v+w\in\circVr_N$, and
\item[(b)] $a\rho(av) = \sum_{v'\in V^r_a}\rho(v-v')$ for all $a\in\Div(N)$ and  $v \in V^r_N\setminus V^r_a$.
\end{itemize}
\end{Def}

\begin{Rem}\label{Gen=StdA}
\rm When $S=\Spec(R)$ and $\CL=\CO_X$, Definition \ref{GenARecipDef} agrees precisely with Definition \ref{ARecipDef}. As in Remark \ref{Gen=StdFq}, giving a general $A$-reciprocal map $\circVr_N\to\CL(S)$ reduces to giving compatible $A$-reciprocal maps $\circVr_N\to R_i$ for a suitable covering of $S$ by open affines $U_i=\Spec(R_i)$.
%Thus all the results from Subsections \ref{AssRingHom} and \ref{ConstantRank}have direct analogues in this more general setting. 
\end{Rem}

\begin{Cons}\label{GenRAVCons}
\rm The given grading on $R_{V^r_N}$ induces a grading on $R_{V^r_N}\otimes_{\BF_q}F$, and all the generators of the ideal in Construction \ref{RAVCons} are homogeneous of degree~$1$. Thus the factor ring $R_{A,V^r_N}$ inherits a unique grading, and so does $RS_{A,V^r_N}$.
For any integer $d$ let $R_{A,V^r_N,d}$ and $RS_{A,V^r_N,d}$ denote the respective homogenous parts of degree~$d$. 
By construction $R_{A,V^r_N}$ is generated over $F$ by its homogeneous part of degree~$1$. Thus 
$$Q_{A,V^r_N} \ :=\ \Proj(R_{A,V^r_N})$$
is a projective scheme over $F$ endowed with a natural very ample invertible sheaf $\CO(1)$ and a natural homomorphism $R_{A,V^r_N,n} \to \CO(n)(Q_{A,V^r_N})$ for all $n\in\BZ$.
Note that it comes with a closed embedding
$$Q_{A,V^r_N}\ \into\ Q_{V^r_N}\times_{\Spec(\BF_q)} \Spec(F).$$
Also, since $RS_{A,V^r_N}$ is the localization of $R_{A,V^r_N}$ obtained by inverting a non-empty finite set of elements of degree~$1$, the scheme
$$\Omega_{A,V^r_N} \ :=\ \Proj(RS_{A,V^r_N})\ \cong\ \Spec(RS_{A,V^r_N,0})$$ 
is an affine open subscheme of $Q_{A,V^r_N}$. 
\end{Cons}

%\begin{Def}\label{ARecipIsoms}
%Consider two pairs $(\CL, \rho)$ and $(\CL', \rho')$ consisting of an invertible sheaf and an $A$-reciprocal map. An isomorphism of invertible sheaves $f\colon \CL\isoto\CL'$ satisfying $\rho'=f\circ\rho$ is called an \emph{isomorphism} $(\CL, \rho) \isoto(\CL', \rho')$. If there exists such an isomorphism, the pairs $(\CL, \rho)$ and $(\CL', \rho')$ are called \emph{isomorphic}.
%\end{Def}
%
%If $\rho$ or $\rho'$ is fiberwise non-zero, there exists at most one isomorphism $(\CL, \rho) \isoto(\CL', \rho')$. Thus the isomorphism classes of such pairs form a well-posed moduli problem. 

As in Subsection \ref{GenFqRecip}, for any two pairs $(\CL, \rho)$ and $(\CL', \rho')$ consisting of an invertible sheaf and a fiberwise non-zero $A$-reciprocal map, there exists at most one isomorphism $(\CL, \rho) \isoto(\CL', \rho')$. Thus the isomorphism classes of such pairs form a well-posed moduli problem. 

\begin{Thm}\label{GenAModuli}
\begin{itemize}
\item[(a)] The composite map 
$$\xymatrix@C-10pt{
\llap{$\rho^\univ\colon$}\ \circVr_N
\ar[rrrr]^-{\textstyle v\mapsto [\frac{1}{v}\otimes1]} &&&& 
R_{A,V^r_N,1} \ar[r] & \CO(1)(Q_{A,V^r_N})}$$
is $A$-reciprocal and fiberwise non-zero.
\item[(b)] For any scheme $S$ over $F$, any invertible sheaf $\CL$ on~$S$, and any fiberwise non-zero $A$-reciprocal map $\rho\colon\circVr_N\to \CL(S)$ there exists a unique morphism $f\colon S\to Q_{A,V^r_N}$ over~$F$ such that $(\CL,\rho)\cong f^*(\CO(1),\rho^\univ)$.
\item[(c)] This $f$ factors through $\Omega_{A,V^r_N}$ if and only if $\rho$ is fiberwise invertible.
\end{itemize}
\end{Thm}

\begin{Proof}
That this $\rho^\univ$ is $A$-reciprocal is a direct consequence of Theorem \ref{AModuli}. That it is fiberwise non-zero follows from the fact that its images $\frac{1}{v}\otimes1$ generate the augmentation ideal of $R_{A,V^r_N}$. This proves (a).

In the situation of (b), the $\BF_q$-reciprocal map underlying $\rho$ already corresponds to a unique morphism $f'\colon S\to Q_{V^r_N}=\Proj(R_{V^r_N})$ such that $(\CL,\rho)\cong (f')^*(\CO(1),\rho^\univ)$ by Theorem \ref{GenFqModuli} (b). As $S$ is a scheme over~$F$, this $f'$ corresponds to a unique morphism $f''\colon S\to \Proj(R_{V^r_N}\otimes_{\BF_q}F)$ over~$F$. The condition \ref{GenARecipDef} (b) for $\rho$ then implies that this morphism factors through the closed subscheme $Q_{A,V^r_N} \subset \Proj(R_{V^r_N}\otimes_{\BF_q}F)$ defined by the graded ideal in Construction \ref{RAVCons}. This yields the desired morphism~$f$, proving (b).

Assertion (c) then follows directly from the construction of $\Omega_{A,V^r_N}$.
\end{Proof}

%%%%%%%%%%%%%%%%%%%%%%%%%%%%%%%%%%%%%

\begin{Thm}\label{MOmegaIsom}
There is a natural isomorphism
$$M^r_{A,N}\ \cong\ \Omega_{A,V^r_N}.$$
\end{Thm}

\begin{Proof}
First consider an affine scheme $S=\Spec(R)$ over~$F$. By Propositions \ref{DrinMod} and \ref{DrinModConv} in the case $W=V^r_N$, giving a triple of the form $(\BG_{a,S},\phi,\lambda)$ consisting of a Drinfeld $A$-module and a level $N$-structure over $S$ is equivalent to giving a fiberwise invertible $A$-reciprocal map $\rho\colon\circVr_N\to R$. 
By Theorem \ref{AModuli} this is in turn equivalent to giving an $F$-algebra homomorphism $f\colon RS_{A,V^r_N}\to R$.

Fix any $v_0\in\circVr_N$. Then in the above equivalences we have $\lambda(v_0)=1$ if and only if $\rho(v_0)=1$ if and only if $f([\frac{1}{v_0}\otimes1])=1$. Giving a triple $(\BG_{a,S},\phi,\lambda)$ with ${\lambda(v_0)=1}$ is thus equivalent to giving an $F$-algebra homomorphism $f\colon RS_{A,V^r_N}\to R$ satisfying ${f([\frac{1}{v_0}\otimes1])=1}$. Since $RS_{A,V^r_N} = RS_{A,V^r_N,0} \bigl[ [\frac{1}{v_0}\otimes1]^{\pm1} \bigr]$, the latter is equivalent to giving an $F$-algebra homomorphism $RS_{A,V^r_N,0}\to R$. That in turn is equivalent to giving a morphism $S=\Spec(R)\to \Spec(RS_{A,V^r_N,0}) = \Omega_{A,V^r_N}$ over~$S$.

On the other hand, by the remark following Theorem \ref{DrinModuli}, giving a triple $(\BG_{a,S},\phi,\lambda)$ with $\lambda(v_0)=1$ is equivalent to giving an isomorphism class of triples $(E,\phi,\lambda)$. Therefore $M^r_{A,N}$ and $\Omega_{A,V^r_N}$ represent isomorphic functors on the category of affine schemes over~$F$. By gluing affine schemes it follows that $M^r_{A,N}$ and $\Omega_{A,V^r_N}$ represent isomorphic functors on the category of all schemes over~$F$. This isomorphism of functors induces the desired isomorphism $M^r_{A,N} \cong \Omega_{A,V^r_N}$.
\end{Proof}

\begin{Cor}\label{RSVIntDom}
The ring $RS_{A,V^r_N}$ is a regular integral domain.
\end{Cor}

\begin{Proof}
Since $\Spec(RS_{A,V^r_N,0}) 
%\cong \Omega_{A,V^r_N} 
\cong M^r_{A,N}$ is a regular integral scheme, the ring $RS_{A,V^r_N,0}$ is a regular integral domain. Thus so is $RS_{A,V^r_N}\cong RS_{A,V^r_N,0}\bigl[ [\frac{1}{v_0}\otimes1]^{\pm1} \bigr]$ for any $v_0\in\circVr_N$.
\end{Proof}

%%%%%%%%%%%%%%%%%%%%%%%%%%%%%%%%%%%%%%%%%%%%%%%%%%%%%%%%%%%%%%%%%%%%%%

\subsection{Stratification}
\label{Strat}

\begin{Prop}\label{i*Morph1A}
% analogue of Proposition \ref{i*Morph1}
Consider any injective $A$-linear map $i\colon V^s_N\into V^r_N$ for $1\le s\le r$. 
\begin{itemize}
\item[(a)] The homomorphism $\pi_i\colon {R_{A,V^r_N}\onto R_{A,V^s_N}}$ from Proposition \ref{i*HomoA} induces a closed embedding $\epsilon_i\colon Q_{A,V^s_N} \into Q_{A,V^r_N}$ whose image is defined by the equations $\rho^\univ(v)=0$ for all $v\in V^r_N\setminus i(V^s_N)$.
\item[(b)] Consider any fiberwise non-zero $A$-reciprocal map $\rho\colon\circVr_N\to \CL(S)$ over a scheme $S$ over~$F$. Then the associated morphism $S\to Q_{A,V^r_N}$ factors through $\epsilon_i$ if and only if $\rho=i_*\rho'$ for an $A$-reciprocal map $\rho'\colon\circVs_N\to \CL(S)$.
\end{itemize}
\end{Prop}

\begin{Proof}
The description in Proposition \ref{i*HomoA} shows that $\pi_i$ is a surjective graded $F$-algebra homomorphism whose kernel is generated by the elements $[\frac{1}{v}\otimes1]$ for all $v\in V^r_N\setminus i(V^s_N)$. This directly implies (a). Part (b) follows as in the proof of Proposition \ref{i*Homo} (c).
\end{Proof}

\begin{Rem}\label{StratDef}
\rm The image subscheme $\epsilon_i(Q_{A,V^s_N})$ depends only on the submodule $W := i(V^s_N)$, which can be any non-zero free $A/N$-submodule of~$V^r_N$. The same holds for $\epsilon_i(\Omega_{A,V^s_N})$, which is the locally closed affine subscheme of $Q_{A,V^r_N}$ that is defined by the equations $\rho^\univ(v)=0$ for all $v\in V^r_N\setminus i(V^s_N)$ and the inequalities $\rho^\univ(v)\not=0$ for all $v\in i(\circVs_N)$.
%in analogy to (\ref{BB}) the open conditions $\rho(v)\not=0$ for all $v\in\circW$ and the closed conditions $\rho(v)=0$ for all $v\in V^r_N\setminus W$.
In the following we abbreviate
\UseTheoremCounterForNextEquation
\begin{equation}\label{OmegaWDef}
\Omega_W\ :=\ \epsilon_i(\Omega_{A,V^s_N}).
\end{equation}
By Theorem \ref{MOmegaIsom} any isomorphism $i\colon V^s_N\isoto W$ induces a natural isomorphism
\UseTheoremCounterForNextEquation
\begin{equation}\label{OmegaWIsom}
M^s_{A,N}\ \cong\ \Omega_W.
\end{equation}
\end{Rem}

\begin{Thm}\label{Stratification}
The subschemes $\Omega_W$ for all non-zero free $A/N$-submodules $W\subset V^r_N$ are pairwise disjoint and their union is $Q_{A,V^r_N}$.
\end{Thm}

\begin{Proof}
(Compare H\"aberli \cite[Thm.\;8.16]{Haeberli}.)
Consider any point of $Q_{A,V^r_N}$ over a field~$k$. This point can be represented by a non-zero $A$-reciprocal map $\rho\colon \smash{\circVr_N}\to k$ which is unique up to multiplication by~$k^\times$. Since $k$ is a field, this map necessarily satisfies the condition (\ref{BB}). The subset $W$ associated to $\rho$ by Lemma \ref{WASubMod} is then a free $A/N$-submodule of some rank ${1\le s\le r}$ by Proposition \ref{DrinMod}, and $\rho$ is the extension by zero of an invertible $A$-reciprocal map $\rho'\colon \circVs_N\to k$ under an isomorphism $i\colon V^s_N\isoto W\subset V^r_N$, as in Proposition \ref{AExtZero}. This $\rho'$ then corresponds to a point in $\Omega_{A,V^s_N}$; hence the original point lies in the stratum~$\Omega_W$. This shows that $Q_{A,V^r_N}$ is the union of all $\Omega_W$. 

Conversely, the construction of $\Omega_W$ implies that the points of $\Omega_W$ over a field~$k$ are precisely those whose associated subset from Lemma \ref{WASubMod} is~$W$. Since $W$ depends only on the 
equivalence class of $\rho$ under multiplication by~$k^\times$, it is uniquely associated to the point. This shows that the $\Omega_W$ are pairwise disjoint.
\end{Proof}

\begin{Thm}\label{Dense}
The open subscheme $\Omega_{A,V^r_N}$ is dense in $Q_{A,V^r_N}$.
\end{Thm}

\begin{Proof}
Consider any point of $Q_{A,V^r_N}$ over a field~$k$. Suppose that it lies in the stratum $\Omega_W$ and choose an isomorphism $i\colon V^s_N\isoto W$. Then the point corresponds to a Drinfeld $A$-module $\phi$ of rank $s$ over $k$ with a level $N$-structure $\lambda\colon V^s_N\to k$. Set $R:=k[[x]]$ and $K := k((x))$ for a new variable~$x$. By Tate uniformization as in \cite[\S7]{Drinfeld1} we can deform $\phi$ to a Drinfeld $A$-module $\tilde\phi$ of rank $r$ with coefficients in $R$ which is congruent to $\phi$ modulo $(x)$. After replacing $R$ and $K$ by a finite extension we may without loss of generality assume that $\tilde\phi[N](K) \cong V^r_N$. Reduction modulo $(x)$ then induces an isomorphism $\tilde\phi[N](K)\cap R \isoto \phi[N](k)$. Via $i$ we can therefore extend $\lambda$ to a level $N$-structure $\tilde\lambda\colon V^r_N\to K$ of~$\tilde\phi$. 

Consider now the associated invertible $A$-reciprocal map $\tilde\rho := (\tilde\lambda|\circVr_N)^{-1}$. By construction it lands in~$R$, and its reduction modulo $(x)$ is the non-zero $A$-reciprocal map corresponding to our given point of~$\Omega_W$. By the modular interpretation of $Q_{A,V^r_N}$ in Theorem \ref{GenAModuli} we therefore obtain a morphism $\Spec(R)\to Q_{A,V^r_N}$ which maps the closed point of $\Spec(R)$ to the given point of $\Omega_W$ and the generic point to $\Omega_{A,V^r_N}$. Thus the given point lies in the closure of $\Omega_{A,V^r_N}$, as desired.
%The idea is to deform this point towards the interior $\Omega_{A,V^r_N}$ by means of Tate uniformization as in \cite[\S7]{Drinfeld1}.
%
%For this set $R:=k[[x]]$ for a new variable~$x$, and let $A$ act on $R$ and $K := k((x))$ through~$\phi$. It is well-known that $K/R$ is a torsion free $A$-module of infinite rank. Thus we can select elements $u_{s+1},\ldots,u_r\in K\setminus R$ whose images in $K/R$ are $A$-linearly independent. Let $L$ be the $A$-submodule of $K$ that is generated by $u_{s+1},\ldots,u_r$. Let $e_{NL}\in R[[\tau]]$ denote the exponential function associated to the submodule $NL$. Then Tate uniformization yields a Drinfeld $A$-module $\tilde\phi$ of rank $r$ with coefficients in $R$ such that $\tilde\phi_a\circ e_{NL} = e_{NL}\circ\phi_a$ for all $a\in A$. The construction also yields an isomorphism of $A$-modules $\tilde\phi[N](K) \cong \phi[N](k) \oplus L/NL$ where any non-zero residue class $u+NL \in L/NL$ corresponds to the element $e_{NL}(u)\in K\setminus R$. 
\end{Proof}

\begin{Cor}\label{OmegaWClosure}
For any non-zero free $A/N$-submodule $W\subset V^r_N$, the closure of $\Omega_W$ in $Q_{A,V^r_N}$ is the union of $\Omega_{W'}$ for all non-zero free $A/N$-submodules $W'\subset W$. 
\end{Cor}

\begin{Proof}
Choose any isomorphism $i\colon V^s_N\isoto W$. Then the closure of $\Omega_{A,V^s_N}$ in $Q_{A,V^s_N}$ is $Q_{A,V^s_N}$ by Theorem \ref{Dense}.
Since $\epsilon_i$ is a closed embedding, it follows that the closure of $\Omega_W = \epsilon_i(\Omega_{A,V^s_N})$ in $Q_{A,V^r_N}$ is $\epsilon_i(Q_{A,V^s_N})$. But $Q_{A,V^s_N}$ is the union of its strata associated to all non-zero free $A/N$-submodules of $V^s_N$; and $\epsilon_i$ maps these strata to the strata $\Omega_{W'}$ associated to all non-zero free $A/N$-submodules of~$W$.
\end{Proof}

%%%%%%%%%%%%%%%%%%%%%%%%%%%%%%%%%%%%%%%%%%%%%%%%%%%%%%%%%%%%%%%%%%%%%%

\subsection{Changing the level}
\label{ChangeLevel}

In this subsection we consider two non-zero proper ideals $N\subset N'\subset A$ which both satisfy Assumption \ref{NAss}.
%Then $\Div(N')\subset\Div(N)$ and $V^r_{N'}\subset V^r_N$. 

\begin{Prop}\label{NN'Prop}
For any $A$-reciprocal map $\rho\colon \circVr_N\to R$  or $\rho\colon \circVr_N\to \CL(S)$,
\begin{enumerate}
\item[(a)] the restriction $\rho|\circVr_{N'}$ is an $A$-reciprocal map, and
\item[(b)] $\rho$ is fiberwise non-zero, resp.\ invertible, if and only if $\rho|\circVr_{N'}$ has the same property.
\end{enumerate}
\end{Prop}

\begin{Proof}
Part (a) follows directly from Definition \ref{ARecipDef}, because for any $a\in\Div(N')$ we have $a\in\Div(N)$ and $V^r_a\subset V^r_{N'}\subset V^r_N$. 
For (b) it suffices to consider an $A$-reciprocal map $\rho\colon \circVr_N\to k$ to a field~$k$. By Proposition \ref{DrinMod} the submodule $W\subset V^r_N$ from Lemma \ref{WASubMod} is then a free $A/N$-module of some rank $0\le s\le r$. It follows that $W\cap V^r_{N'}$ is a free $A/N'$-module of the same rank~$s$. Thus $\rho$ is non-zero if and only if $s\ge1$ if and only if $\rho|\circVr_{N'}$ is non-zero, and $\rho$ is invertible if and only if $s=r$ if and only if $\rho|\circVr_{N'}$ is invertible.
\end{Proof}

\begin{Cons}\label{NN'Cons}
\rm By the universal property in Theorem \ref{AModuli} (b), the restriction of $A$-reciprocal maps corresponds to a natural $F$-algebra homomorphism $R_{A,V^r_{N'}} \to R_{A,V^r_N}$ which sends $[\tfrac{1}{v}\otimes1]$ to $[\tfrac{1}{v}\otimes1]$ for all $v\in\smash{\circVr_{N'}}$. By Construction \ref{RAVCons} this induces a commutative diagram of graded $F$-algebras
\UseTheoremCounterForNextEquation
\begin{equation}\label{NN'R}
\vcenter{\xymatrix{
RS_{A,V^r_N} & R_{A,V^r_N} \ar[l] \\
RS_{A,V^r_{N'}} \ar[u] & R_{A,V^r_{N'}} \ar[l] \ar[u] \\}}
\end{equation}
Applying $\Proj$ as in Construction \ref{GenRAVCons} this yields a commutative diagram of schemes over~$F$
\UseTheoremCounterForNextEquation
\begin{equation}\label{NN'Q}
\vcenter{\xymatrix@C-10pt{
M^r_{A,N} \ar@{}[r]|-\cong \ar[d] & \Omega_{A,V^r_N} \ar[d] \ar@{^{ (}->}[rr] && Q_{A,V^r_N} \ar[d] \\
M^r_{A,N'} \ar@{}[r]|-\cong & \Omega_{A,V^r_{N'}} \ar@{^{ (}->}[rr] && Q_{A,V^r_{N'}} \\}}
\end{equation}
where the vertical morphism on the left is induced by the restriction of a level $N$-structure on a Drinfeld $A$-module to a level $N'$-structure. The other two vertical morphisms represent the restriction of isomorphism classes of fiberwise invertible, resp. fiberwise non-zero, $A$-reciprocal maps. Moreover, the `if' part of Proposition \ref{NN'Prop} (b) implies that this diagram is cartesian.
\end{Cons}

\begin{Rem}\label{NN'Rem}
\rm It is known that $M^r_{A,N} \to M^r_{A,N'}$ is a finite \'etale Galois covering with Galois group $\Ker(\GL_r(A/N)\onto\GL_r(A/N'))$. By Corollary \ref{RSVIntDom} it follows that $RS_{A,V^r_{N'}} \to RS_{A,V^r_N}$ is a finite \'etale Galois extension of integral domains. In particular it is injective. 

We actually expect that all homomorphisms in the diagram (\ref{NN'R}) are injective and that $R_{A,V^r_{N'}} \to R_{A,V^r_N}$ is a finite ring extension, but cannot conclude that (yet) at this point.
\end{Rem}

%%%%%%%%%%%%%%%%%%%%%%%%%%%%%%%%%%%%%%%%%%%%%%%%%%%%%%%%%%%%%%%%%%%%%%

\subsection{Satake compactification}
\label{Satake}

%First we recall the axiomatic characterization of the Satake compactification of $M^r_{A,N}$ from \cite[\S3, \S4]{PinkSatake}. 
In this subsection we relate $Q_{A,V^r_N}$ with the Satake compactification $\OM^r_{A,N}$ of $M^r_{A,N}$. For this we recall the axiomatic characterization of $\OM^r_{A,N}$ and its properties from \cite[\S\S3-5]{PinkSatake}.
%To explain the latter we recall some notions from~\cite[\S3]{PinkSatake}. 

\medskip
As before we fix an integer $r\ge1$ and consider a commutative $F$-algebra~$R$. Following \cite[Def.\,3.1]{PinkSatake}, a \emph{standard generalized Drinfeld $A$-module of rank $\le r$ over $R$} is an $\BF_q$-algebra homomorphism $\phi\colon A \to R[\tau]$, $a\mapsto\phi_a$ satisfying for every $a\in A\setminus\{0\}$:
\begin{samepage}
\begin{enumerate}
\item[(a)] $\phi_a = \sum_{i=0}^{r\deg_A(a)}\phi_{a,i}\tau^i$ and for every $\Fp\in\Spec(R)$ there exists $i>0$ with $\phi_{a,i}\not\in\Fp$, and
\item[(b)] $d\phi_a = \phi_{a,0} = a$.
\end{enumerate}
\end{samepage}
More generally consider a scheme $S$ over~$F$ and a line bundle $E$ on~$S$. An arbitrary \emph{generalized Drinfeld $A$-module of rank $\le r$ over $R$} is an $\BF_q$-algebra homomorphism $\phi\colon A \to \End_{\BF_q}(E)$, $a\mapsto\phi_a$ which for any trivialization of $E$ over an open affine subscheme becomes a standard generalized Drinfeld $A$-module of rank $\le r$. Then the fiber over any point of $S$ must be a Drinfeld $A$-module of some rank $1\le s\le r$, but this $s$ can vary over~$S$.

An \emph{isomorphism} of generalized Drinfeld $A$-modules is an isomorphism of line bundles that is compatible with~$\phi$. Following 
\cite[Def.\,3.9]{PinkSatake} we call a generalized Drinfeld $A$-module $(E,\phi)$ over $S$ \emph{weakly separating} if, for any Drinfeld $A$-module $(E',\phi')$ over any field $L$ containing~$F$, at most finitely many fibers of $(E,\phi)$ over $L$-valued points of $S$ are isomorphic to $(E',\phi')$.

Recall from Theorem \ref{DrinModuli} that $M^r_{A,N}$ is an integral affine algebraic variety of finite type over~$F$.
By \cite[Def.\,4.1]{PinkSatake}, any open embedding $M^r_{A,N} \into \OM^r_{A,N}$ into a normal integral proper algebraic variety over~$F$,  such that the universal family $(E^\univ,\phi^\univ)$ on $M^r_{A,N}$ extends to a weakly separating generalized Drinfeld $A$-module $(\bar{E}^\univ,\bar{\phi}^\univ)$ over $\OM^r_{A,N}$, is called a \emph{Satake compactification of $M^r_{A,N}$}. By abuse of terminology we call $(\bar{E}^\univ,\bar{\phi}^\univ)$ the \emph{universal family on $\OM^r_{A,K}$}.
By \cite[Thm.\,4.2]{PinkSatake} such a Satake compactification exists, and it together with its universal family is unique up to unique isomorphism. Moreover, let $\bar\CL$ denote the dual of the 
relative Lie algebra of $\bar{E}^\univ$, which is an invertible sheaf on $\OM^r_{A,N}$. Then $\smash{\OM^r_{A,N}}$ is projective over $F$ and $\bar\CL$ is ample by \cite[Thm.\,5.3]{PinkSatake}.

%%%%%%%%%%%%%%%%%%%%%%%%%%%%%%%%%%%%%

\medskip
Now we return to $A$-reciprocal maps.

\begin{Prop}\label{PhiSCons}
Consider any 
%fiberwise non-zero 
$A$-reciprocal map $\rho\colon\circVr_N\to\CL(S)$ over~$S$. Let $E$ be the line bundle on $S$ whose sheaf of sections is the dual~$\CL^\vee$. Then there is a unique $\BF_q$-algebra homomorphism $\phi\colon A\to\End_{\BF_q}(E)$ which for any trivialization $\BG_{a,U}\isoto E|U$ over an open affine subscheme $U=\Spec(R)\subset S$ induces the homomorphism $A\to \End_{\BF_q}(E|U) \cong R[\tau]$ from Proposition \ref{RingHomo1}.
\end{Prop}

\begin{Proof}
First consider any $a\in\Div(N)$. Then for any open subscheme $U\subset S$, any section $e\in E(U)$, and any $v\in\circVr_N$, the expression $1-\rho(v)e$ is a well-defined section of $\BG_a(U)=\CO_S(U)$. Thus the formula (\ref{Phia}) globalizes to a morphism $\phi_a\colon E\to E$ over~$S$. Since locally over $S$ it is $\BF_q$-linear, it defines an element of $\End_{\BF_q}(E)$. As $A$ is generated by $\Div(N)$, we obtain $\phi_a\in \End_{\BF_q}(E)$ for all $a\in A$.
\end{Proof}

\medskip
We apply Proposition \ref{PhiSCons} to the universal $A$-reciprocal map $\rho^\univ\colon A\to \CO(1)(Q_{A,V^r_N})$ from Theorem \ref{GenAModuli} and obtain a line bundle $E$ on $Q_{A,V^r_N}$ with an $\BF_q$-algebra homomorphism
\UseTheoremCounterForNextEquation
\begin{equation}\label{BarPhi}
\phi\colon A\to \End_{\BF_q}(E).
\end{equation}
%So $E$ is the line bundle on $Q_{A,V^r_N}$ whose invertible sheaf of sections is the dual of the invertible sheaf $\CO(1)$. 
% Conversely $\CL$ is the dual of the relative Lie algebra of~$E$. 

\begin{Prop}\label{BarPhiWeakSep}
This $\phi$ is a weakly separating generalized Drinfeld module of rank~$\le r$.
\end{Prop}

\begin{Proof}
(Compare H\"aberli \cite[Cor.\;8.21]{Haeberli}.)
For any $a\in\Div(N)$ we have $\dim_{\BF_q}(V^r_a) = r\dim_{\BF_q}(A/(a)) = r\deg_A(a)$; hence the formula (\ref{Phia}) shows that 
\UseTheoremCounterForNextEquation
\begin{equation}\label{PhiaDeg}
\phi_a\ = \sum_{i=0}^{r\deg_A(a)}\kern-5pt\phi_{a,i}\tau^i.
\end{equation}
In particular this holds for any element $a$ of the set $D$ from Assumption \ref{NAss}. As the degree is additive in products, the expansion (\ref{PhiaDeg}) follows whenever $a$ is a product of elements of~$D$. It also holds for $a=0$, because the empty sum is zero. We claim that it holds for all $a\in A$.

To prove this we use induction on $\deg_A(a)$. Consider any integer $d\ge0$ and suppose that (\ref{PhiaDeg}) holds for all elements $a\in A$ with $\deg_A(a)<d$. Consider an $a\in A$ with $\deg_A(a)=d$ and choose $b$ as in Assumption \ref{NAss} (a). Then (\ref{PhiaDeg}) holds for $b$ in place of~$a$, and we have $\deg_A(a-b)<\deg_A(a)$. Thus (\ref{PhiaDeg}) holds for $a-b$ by the induction hypothesis; and so it holds for $a=(a-b)+b$ as well, finishing the induction proof of the claim.

Next, from (\ref{OmegaWIsom}) we know that $Q_{A,V^r_N}$ is the union of finitely many strata $\Omega_W$ with isomorphisms $\Omega_W \cong \Omega_{A,V^s_N} \cong M^s_{A,N}$ for varying $1\le s\le r$. Moreover, by construction the pullback of the universal $A$-reciprocal map $\rho^\univ$ from $Q_{A,V^r_N}$ to $\Omega_{A,V^s_N}$ is simply the universal $A$-reciprocal map on $\Omega_{A,V^s_N}$. Thus the pullback of the pair $(E,\phi)$ is simply the Drinfeld $A$-module of rank $s$ over $\Omega_{A,V^s_N}$ that corresponds to the universal $A$-reciprocal map over $\Omega_{A,V^s_N}$. Since $s\ge1$, this and the above claim show that $\phi$ satisfies all conditions for a generalized Drinfeld module of rank $\le r$.
Also, transferring $(E,\phi)$ to $M^s_{A,N}$ yields the universal Drinfeld $A$-module over $M^s_{A,N}$ with the level $N$-structure removed. As any Drinfeld $A$-module over a field possesses only finitely many level $N$-structures, the universal Drinfeld $A$-module over $M^s_{A,N}$ is weakly separating. As we have only finitely many strata altogether, it follows that $(E,\phi)$ is weakly separating.
\end{Proof}

%%%%%%%%%%%%%%%%%%%%%%%%%%%%%%%%%%%%%

\begin{Cons}\label{RAVNormalize}
\rm Let $R_{A,V^r_N}^\norm$ denote the integral closure of $R_{A,V^r_N}$ in $RS_{A,V^r_N}$. As $RS_{A,V^r_N}$ is a normal integral domain by Corollary \ref{RSVIntDom}, so is $R_{A,V^r_N}^\norm$ and we have ${R_{A,V^r_N}^\norm \otimes_{R_{A,V^r_N}}RS_{A,V^r_N}} \allowbreak = RS_{A,V^r_N}$. The ring $R_{A,V^r_N}^\norm$ inherits a grading, so that we can define
$$Q_{A,V^r_N}^\norm \ :=\ \Proj(R_{A,V^r_N}^\norm).$$
The natural isomorphism $M^r_{A,N} \cong \Omega_{A,V^r_N} \subset Q_{A,V^r_N}$ from Theorem \ref{MOmegaIsom} yields an open embedding 
$$M^r_{A,N}\ \longinto\ Q_{A,V^r_N}^\norm.$$
\end{Cons}

\begin{Prop}\label{RAVNormalMap}
The natural morphism $\pi\colon Q_{A,V^r_N}^\norm \to Q_{A,V^r_N}$ is finite and surjective.
%and an isomorphism over $\Omega_{A,V^r_N}$.
\end{Prop}

\begin{Proof}
Since $R_{A,V^r_N}$ is an algebra of finite type over a field, the ring $R_{A,V^r_N}^\norm$ is a finite $R_{A,V^r_N}$-algebra by Noether's theorem \cite[Thm.\;4.14]{Eisenbud}.
Thus $\pi$ is a finite morphism. Therefore its image is closed in $Q_{A,V^r_N}$. As this image contains $\Omega_{A,V^r_N}$, which by Theorem \ref{Dense} is dense in $Q_{A,V^r_N}$, it follows that the image of $\pi$ is $Q_{A,V^r_N}$. Thus $\pi$ is finite and surjective.
\end{Proof}

\begin{Thm}\label{SatakeProjection}
The scheme $Q_{A,V^r_N}^\norm$ is the Satake compactification $\OM^r_{A,N}$ of $M^r_{A,N}$.
\end{Thm}

\begin{Proof}
(Compare H\"aberli \cite[Cor.\;8.22]{Haeberli}.)
By construction $Q_{A,V^r_N}^\norm$ is a normal integral proper algebraic variety over $F$ which contains $\Omega_{A,V^r_N} \cong M^r_{A,N}$ as an open subvariety. Moreover, since $(E,\phi)$ is weakly separating and $\pi$ is finite, the pullback $\pi^*(E,\phi)$ is a weakly separating generalized Drinfeld $A$-module which extends  the universal family $(E^\univ,\phi^\univ)$ on $M^r_{A,N}$. By the uniqueness part of \cite[Thm.\,4.2]{PinkSatake} it follows that $Q_{A,V^r_N}^\norm$ \emph{is} the Satake compactification of $M^r_{A,N}$.
\end{Proof}

\begin{Rem}\label{NormNotIsomRem}
\rm The computation in H\"aberli \cite[Prop.\;7.13, Cor.\;7.28]{Haeberli} 
%\cite[Cor.\;9.15]{Haeberli} 
implies that the fiber over every geometric point at the boundary consists of $|\mathop{\rm Pic}(A)|\cdot |(A/N)^\times/A^\times|$ geometric points.
Usually $\pi\colon Q_{A,V^r_N}^\norm \to Q_{A,V^r_N}$ is therefore not an isomorphism. A fortiori $R_{A,V^r_N} \to R_{A,V^r_N}^\norm$ is not an isomorphism in general.
\end{Rem}

%%%%%%%%%%%%%%%%%%%%%%%%%%%%%%%%%%%%%%%%%%%%%%%%%%%%%%%%%%%%%%%%%%%%%%

\subsection{The ideal of the boundary}
\label{ABoundaryIdeal}

For any $1\le s<r$ and any injective $A$-linear map $i\colon V^s_N\into V^r_N$ consider the composite ring homomorphism 
$$\xymatrix@C+10pt{
\tilde\pi_i\colon R_{A,V^r_N} \ar@{->>}[r]^-{\pi_i}
& R_{A,V^s_N} \ar[r] & RS_{A,V^s_N} \\}$$
where $\pi_i$ is the homomorphism from Proposition \ref{i*HomoA}. 
Also, let $R_{A,V^r_N,+} := \bigoplus_{d>0} R_{A,V^r_N,d}$ denote the augmentation ideal of $R_{A,V^r_N}$.
We are interested in the ideal
\UseTheoremCounterForNextEquation
\begin{equation}\label{AIVNDef}
I_{A,V^r_N}\ :=\ R_{A,V^r_N,+}\cap \bigcap_{\text{all}\ s,\;i} \kern-0pt\Ker(\tilde\pi_i) \ \subset\ R_{A,V^r_N}.
\end{equation}
Since $RS_{A,V^s_N}$ is an integral domain, this ideal is reduced. 
By construction it is graded, so it defines a reduced closed subscheme $\partial\Omega_{A,V^r_N}$ of $Q_{A,V^r_N}$. By Theorem \ref{Stratification} its complement is $\Omega_{A,V^r_N}$; hence 
\UseTheoremCounterForNextEquation
\begin{equation}\label{PartialOmegaAVN}
\partial\Omega_{A,V^r_N}\ =\ (Q_{A,V^r_N}\setminus\Omega_{A,V^r_N})^\red.
\end{equation}
%Next consider the integral closure $R_{A,V^r_N}^\norm$ from Construction \ref{RAVNormalize}. 
We are also interested in the ideal
\UseTheoremCounterForNextEquation
\begin{equation}\label{AIVNnormDef}
I_{A,V^r_N}^\norm\ :=\ \sqrt{I_{A,V^r_N}\cdot R_{A,V^r_N}^\norm} 
\ \ \subset\ R_{A,V^r_N}^\norm.
\end{equation}
By construction this is a reduced graded ideal of $R_{A,V^r_N}^\norm$. Recall from Theorem \ref{SatakeProjection} that $Q_{A,V^r_N}^\norm = \OM^r_{A,N}$ is the Satake compactification of $M^r_{A,N}$.
Thus the closed subscheme associated to $I_{A,V^r_N}^\norm$ is the reduced subscheme at the boundary
\UseTheoremCounterForNextEquation
\begin{equation}\label{PartialOmegaAVNnorm}
\partial M^r_{A,N}\ :=\ (\OM^r_{A,N}\setminus M^r_{A,N})^\red.
\end{equation}
We expect that $I_{A,V^r_N} = I_{A,V^r_N}^\norm$. In Subsection \ref{SpecialABoundaryIdeal} we will prove this in a special case.

%%%%%%%%%%%%%%%%%%%%%%%%%%%%%%%%%%%%%%%%%%%%%%%%%%%%%%%%%%%%%%%%%%%%%%

\subsection{Modular forms and cusp forms}
\label{ModCuspForms}

Let $\CO(1)$ denote the pullback to $\OM^r_{A,N}$ of the very ample invertible sheaf $\CO(1)$ under the morphism $\pi$ from Proposition \ref{RAVNormalMap}. As usual, for any quasicoherent sheaf $\CF$ on $\OM^r_{A,N}$ and any integer $d$ we set $\CF(d) := \CF\otimes\CO(1)^{\otimes d}$. 
Let $\CI\subset\CO_{\OM^r_{A,N}}$ denote the ideal sheaf of the reduced boundary $\partial M^r_{A,N}$ from (\ref{PartialOmegaAVNnorm}). In other words it is the ideal sheaf associated to the graded ideal $I_{A,V^r_N}^\norm\subset R_{A,V^r_N}^\norm$.
For any integer $d$ we call
\UseTheoremCounterForNextEquation
\begin{eqnarray}\label{ModFormsDef}
\Gamma(\OM^r_{A,N},\CO(d)) &&\hbox{the space of \emph{modular forms} and} \\[3pt]
\UseTheoremCounterForNextEquation
\label{CuspFormsDef}
\Gamma(\OM^r_{A,N},\CI(d))\, &&\hbox{the space of \emph{cusp forms}}
\end{eqnarray}
\emph{of rank $r$ and level $N$ and weight~$d$}. 

\medskip
For any integer $d$ let $I_{A,V^r_N,d}^\norm \subset R_{A,V^r_N,d}^\norm$ denote the homogenous parts of degree $d$ of $I_{A,V^r_N}^\norm \subset R_{A,V^r_N}^\norm$.

\begin{Thm}\label{ModCuspFormsIsom}
For any $d\ge1$ we have natural isomorphisms
$$\begin{array}{lll}
R_{A,V^r_N,d}^\norm & \!\isoto\ \Gamma(\OM^r_{A,N},\CO(d)) & 
\hbox{and} \\[5pt]
I_{A,V^r_N,d}^\norm & \!\isoto\ \Gamma(\OM^r_{A,N},\CI(d)). & 
\end{array}$$
\end{Thm}

\begin{Proof}
Since $\OM^r_{A,N}$ is normal with graded coordinate ring $R_{A,V^r_N}^\norm$, by \cite[Ch.II, Ex.\;5.14(a)]{Hartshorne} (whose proof does not require that $R_{A,V^r_N}^\norm$ be generated by elements of degree~$1$) the ring $\bigoplus_{d\ge0}\Gamma(\OM^r_{A,N},\CO(d))$ is the integral closure of $R_{A,V^r_N}^\norm$, and hence equal to $R_{A,V^r_N}^\norm$. This yields the first isomorphism, and that in turn directly implies the second.
\end{Proof}

\newpage

%%%%%%%%%%%%%%%%%%%%%%%%%%%%%%%%%%%%%%%%%%%%%%%%%%%%%%%%%%%%%%%%%%%%%%
%%%%%%%%%%%%%%%%%%%%%%%%%%%%%%%%%%%%%%%%%%%%%%%%%%%%%%%%%%%%%%%%%%%%%%

\section{The special case $A=\BF_q[t]$ and $N=(t^n)$}
\label{Specialn}

%%%%%%%%%%%%%%%%%%%%%%%%%%%%%%%%%%%%%%%%%%%%%%%%%%%%%%%%%%%%%%%%%%%%%%

\subsection{Setup}
\label{SpecialnSetup}

Throughout this section we assume that $A=\BF_q[t]$ and $N=(t^n)$ for some $n\ge1$. Then 
$$\Div(N)\ =\ \{\alpha t^\nu\mid \alpha\in\BF_q^\times,\ 0\le \nu\le n\},$$
so it satisfies Assumption \ref{NAss} with the subset $D := \BF_q^\times\cup\{t\}$. To reduce notation we fix $r\ge1$ and abbreviate $V_n := V^r_{t^n} = (t^{-n}A/A)^{\oplus r}$. By induction on $\nu$ and a short computation we have:

\begin{Lem}\label{AtRecipLem}
A map $\rho\colon V_n\to R$ is $A$-reciprocal if and only if it is $\BF_q$-reciprocal and satisfies 
$$t\rho(tv) = \sum_{v'\in V_1}\rho(v-v')$$ 
for all $v \in V_n\setminus V_1$.
\end{Lem}
%the condition \ref{GenARecipDef} (b) for $a=t$ already implies the condition for $a=t^\nu$ for all $1\le\nu\le n$. Thus the condition \ref{GenARecipDef} (b) for all $a\in D$ implies the condition for all $a\in\Div(N)$.

Next we abbreviate $\tilde R_n := R_{V_n}\otimes_{\BF_q}F$ and let $J_n$ be its ideal from Construction \ref{RAVCons}. As a consequence of Lemma \ref{AtRecipLem} the ideal $J_n$ is already generated by the relations
\UseTheoremCounterForNextEquation
\begin{equation}\label{RelDef}
\Rel_v\ :=\ \frac{1}{tv}\otimes t-\sum_{v'\in V_1}\frac{1}{v-v'}\otimes1
\end{equation}
for all $v \in V_n\setminus V_1$. We abbreviate the factor ring as $R_n := \tilde R_n/J_n = \smash{R_{A,V^r_{t^n}}}$ and denote the projection map by $\pi\colon\allowbreak {\tilde R_n\onto R_n}$. 
We also abbreviate $RS_n := RS_{A,V^r_{t^n}}$.
% and $Q_n := Q_{A,V^r_{t^n}}$ and $\Omega_n := \Omega_{A,V^r_{t^n}}$.
The natural action of $\GL_r(\BF_q[t]/(t^n))$ on $V_n$ induces an action on $\smash{\tilde R_n}$ and $J_n$ and hence on $R_n$ and $RS_n$.

\medskip
Note that in the case $n=1$ there are no relations (\ref{RelDef}); hence $R_1=\tilde R_1=R_{V_1}\otimes_{\BF_q}F$ and $RS_1=RS_{V_1}\otimes_{\BF_q}F$, and the schemes $\Omega_1\subset Q_1$ are obtained from $\Omega_{V_1}\subset Q_{V_1}$ by base change from $\Spec(\BF_q)$ to $\Spec(F)$. In fact, in this case any $\BF_q$-reciprocal map is already $A$-reciprocal by Lemma \ref{AtRecipLem}.
%$$\begin{array}{rl}
%Q_{A,V^r_t} &=\ \phantom{\partial} Q_{V^r_t}\times_{\Spec(\BF_q)}\Spec(F), \\[3pt]
%\Omega_{A,V^r_t} &=\ \phantom{\partial} \Omega_{V^r_t}\times_{\Spec(\BF_q)}\Spec(F),\ \hbox{and} \\[3pt]
%\partial\Omega_{A,V^r_t} &=\ \partial\Omega_{V^r_t}\times_{\Spec(\BF_q)}\Spec(F) \\[3pt]
%\end{array}$$
%for the closed subscheme $\partial\Omega_{V^r_t} := (Q_{V^r_t}\setminus\Omega_{V^r_t})^\red$. 

%%%%%%%%%%%%%%%%%%%%%%%%%%%%%%%%%%%%%%%%%%%%%%%%%%%%%%%%%%%%%%%%%%%%%%

\subsection{Description of the ring}
\label{Specialn2}

The goal of this subsection is to give an explicit description of the ring $R_n$ for arbitrary~$n$. These results are collected in Theorem \ref{AAll1}, but in order to get there, we must first introduce some auxiliary notation.

\medskip
Let $b_1,\ldots,b_r$ denote the standard basis of $\BF_q^{\oplus r}$. Then the elements $X_{k,\nu} := [t^{-\nu}b_k]$ for all ${1\le k\le r}$ and ${1\le\nu\le n}$ form a basis of $V_n$ over~$\BF_q$. We bring these elements in the order
\UseTheoremCounterForNextEquation
\begin{equation}\label{VarOrder}
X_{1,1},X_{2,1},\ldots,X_{r,1},X_{1,2},\ldots,X_{r,2},X_{1,3},\ldots,\ldots,X_{r,n-1},X_{1,n},\ldots,X_{r,n}.
\end{equation}
Then multiplying $X_{k,\nu}$ by $t$ yields $0$ if $\nu=1$, or an earlier element of the list if $\nu>1$. 
In particular $X_{1,n},\ldots,X_{r,n}$ is a basis of $V_n$ over $A/(t^n)$. For any $(k,\nu)$ we let $V'_{k,\nu}$ denote the $\BF_q$-subspace that is generated by all elements occurring strictly before $X_{k,\nu}$. In other words $V'_{k,\nu}$ is the direct sum of $V_{\nu-1}\subset V_n$ and the $\BF_q$-subspace generated by $X_{1,\nu},\ldots,X_{k-1,\nu}$.

\medskip
For every $1\le k\le r$ we consider the finite subsets
%in analogy to (\ref{DkDef}) and (\ref{EkDef}), 
\UseTheoremCounterForNextEquation
\begin{eqnarray}
\label{tDkDef}
\tilde\Delta_k &:=& \biggl\{ \frac{1}{X_{k,n} + w}\otimes1 \,\biggm|\, w \in \circVlimits_{k,n}^{\prime} \biggr\} \cup\bigl\{1\bigr\}
\quad\hbox{and} \\
\UseTheoremCounterForNextEquation
\label{tEkDef}
\tilde E_k &:=& \biggl\{ \frac{1}{X_{k,n} + w}\otimes1 \,\biggm|\, w \in V_{k,n}^{\prime} \biggr\}
\end{eqnarray}
of $\tilde R_n$, each of cardinality $|\circVlimits_{k,n}^{\prime}| = q^{r(n-1)+k-1}$. Observe that by unique factorization for polynomials in the variables $X_{k,\nu}$ we have bijective maps
$$\xymatrix@R-20pt@C-13pt{
\tilde\Delta_1\times\ldots\times\tilde\Delta_r \ 
\ar[rr]^-\sim && \ \tilde\Delta_1\cdots\tilde\Delta_r 
\ar@{}[r]|-{\subset} & \tilde R_n \rlap{\hbox{\qquad and}}\\
\tilde E_1\times\ldots\times\tilde E_r \ 
\ar[rr]^-\sim && \ \tilde E_1\cdots\tilde E_r 
\ar@{}[r]|-{\subset} & \tilde R_n \rlap{.} \\}$$
Let $U<\GL_r(\BF_q[t]/(t^n))$ be the subgroup of matrices which are congruent modulo $(t)$ to an upper triangular matrix with diagonal entries~$1$. Viewing $V_n$ as a space of column vectors and letting $\GL_r(\BF_q[t]/(t^n))$ act on $V_n$ by left multiplication, this is the subgroup of all $g\in\GL_r(\BF_q[t]/(t^n))$ such that $g(X_{k,n})\in X_{k,n}+V'_{k,n}$ for all $1\le k\le r$. In fact, any independent choice of an element of $\smash{X_{k,n}+V'_{k,n}}$ for all $k$ corresponds to a unique element of~$U$. 
For the induced action on $\tilde R_n$ it follows that $U$ acts transitively on $\tilde E_k$ and freely transitively on $\tilde E_1\times\ldots\times\tilde E_r$. 
It also follows that for each $1\le k\le r$ the element
\UseTheoremCounterForNextEquation
\begin{equation}\label{AfkDef}
\tilde f_k
\ := \sum_{e_k\in\tilde E_k} e_k 
\ = \sum_{w \in V'_{k,n}} \frac{1}{X_{k,n} + w}\otimes 1
\ \in \ \tilde R_n
\end{equation}
is fixed by~$U$. 
By Theorem \ref{RVBasis}
%\cite[Thm.\,2.7]{PinkSchieder} 
and tensoring with $F$ the elements $\tilde f_1,\ldots,\tilde f_r$ are algebraically independent over $F$ and the elements of $\tilde\Delta_1\cdots\tilde\Delta_r$ are linearly independent over the subring 
%$\tilde R_n' := 
$F[\tilde f_1,\ldots,\tilde f_r]$. We are interested in the submodule
\UseTheoremCounterForNextEquation
\begin{equation}\label{TMnDef}
\tilde M_n\ := 
\bigoplus_{e\in\tilde\Delta_1\cdots\tilde\Delta_r} \kern-5pt
F[\tilde f_1,\ldots,\tilde f_r]\cdot e
\ \subset\ \tilde R_n.
\end{equation}

\begin{Lem}\label{MAsurj}
The projection $\pi$ 
%$\pi\colon \tilde R_n\onto R_n$
induces a surjection $\tilde M_n\to R_n$.
\end{Lem}

\begin{Proof}
The assertion is equivalent to $\tilde R_n=\tilde M_n+J_n$. By the construction of $J_n$ we must therefore show that for every $f\in\tilde R_n$ there exist $m\in\tilde M_n$ and elements $g_v\in\tilde R_n$ such that
\UseTheoremCounterForNextEquation
\begin{equation}\label{Msurj1}
f\ =\ m + \kern-5pt \sum_{v \in V_n\setminus V_1}\kern-8pt g_v\cdot\Rel_v.
\end{equation}
For this recall that $\tilde R_n=R_{V_n}\otimes_{\BF_q}F$. Also observe that $\tilde\Delta_k=\{e\otimes1\mid e\in \Delta'_k\}$ for the subset $\Delta'_k := \{ \frac{1}{X_{k,n} + w}\mid w \in \circVlimits_{k,n}^{\prime}\}\cup\{1\} \subset R_{V_n}$ and that $\tilde f_k =f'_k\otimes1$ with $f'_k := \sum_{w \in V_{k,n}}\frac{1}{X_{k,n} + w}$. Thus $\tilde M_n = M'_n\otimes_{\BF_q}F$ for the submodule $M'_n := \bigoplus_{e\in\Delta'_1\cdots\Delta'_r} \BF_q[f'_1,\ldots,f'_r]\cdot e \subset R_{V_n}$.
Then the variable $t\in F=\BF_q(t)$ enters into the equation (\ref{Msurj1}) only through the relations $\Rel_v$. Rescaling these to 
$$\Rel'_v\ :=\ \Rel_v\cdot\;t^{-1}\ =\ 
\frac{1}{tv}\otimes 1-\sum_{v'\in V_1}\frac{1}{v-v'}\otimes t^{-1},$$ 
it suffices to show that for every $f'\in R_{V_n}$ there exist elements $m\in M'_n\otimes_{\BF_q}F$ and $g_v\in R_{V_n}\otimes_{\BF_q}F$ such that
\UseTheoremCounterForNextEquation
\begin{equation}\label{Msurj2}
f'\otimes1\ =\ m + \kern-5pt \sum_{v \in V_n\setminus V_1}\kern-8pt g_v\cdot\Rel'_v.
\end{equation}
This is a vector valued inhomogeneous linear equation with coefficients in the ring $\BF_q[t^{-1}]$. In terms of any bases of $R_{V_n}$ and $M'_n$ over~$\BF_q$, it is equivalent to a system of inhomogeneous linear equations with coefficients in $\BF_q[t^{-1}]$. 
Our job is to find a solution in $\BF_q(t)$. But a system of inhomogeneous linear equations over a field has a solution over that field if and only if it has a solution over any overfield. It therefore suffices to find a solution over the completion $\BF_q((t^{-1}))$. In fact, we will show that there exists a solution in $\BF_q[[t^{-1}]]$.

We first find a solution modulo $(t^{-1})$. For this observe that $\Rel'_v\equiv \frac{1}{tv}\otimes 1$ modulo $(t^{-1})$, and that as $v$ runs through $V_n\setminus V_1$, the element $v':=tv$ runs through $\smash{\circV_{n-1}}$. To solve the problem modulo $(t^{-1})$ we must therefore show that for every $f'\in R_{V_n}$ there exist $m'\in M'_n$ and elements $g'_v\in R_{V_n}$ such that
\UseTheoremCounterForNextEquation
\begin{equation}\label{Msurj3}
f'\ =\ m' + \kern-5pt \sum_{v'\in\circV_{n-1}}\kern-5pt g'_v\cdot\frac{1}{v'}.
\end{equation}
This is equivalent to saying that $R_{V_n}=M'_n+J'_n$ for the ideal $J'_n\subset R_{V_n}$ that is generated by the elements $\frac{1}{v'}$ for all $v'\in\circV_{n-1}$, or again to saying that the induced homomorphism $M'_n\to R_{V_n}/J'_n$ is surjective. But that is guaranteed by Theorem \ref{RVJBasis} on identifying the data $(V,f_{s+k},\Delta_{s+k})$ from Subsection \ref{FqDescRing} with the present $(V_n,f'_k,\Delta'_k)$ and the data $(V_s,J_s,M_s)$ from Subsection \ref{FqComplem} with the present $(V_{n-1},J'_n,M'_n)$.

Next consider any $i\ge1$ and suppose that for all $0\le j<i$ we already have $m_j\in M'_n$ and $g_{v,j}\in R_{V_n}$ such that $m=\sum_{j=0}^{i-1} m_j\otimes t^{-j}$ and $g_v=\sum_{j=0}^{i-1} g_{v,j}\otimes t^{-j}$ solve the equation (\ref{Msurj2}) modulo $(t^{-i})$. Then the left hand side minus the right hand side is congruent to $f'_i\otimes t^{-i}$ modulo $(t^{-i-1})$ for some $f'_i\in R_{V_n}$. Solve the equation (\ref{Msurj3}) with $(f'_i,m_i,g_{v,i})$ in place of $(f',m',g'_v)$. Then the elements $m=\sum_{j=0}^{i} m_j\otimes t^{-j}$ and $g_v=\sum_{j=0}^{i} g_{v,j}\otimes t^{-j}$ solve the equation (\ref{Msurj2}) modulo $(t^{-i-1})$. 

Before passing to the limit we must take care of one more problem. Namely, observe that $R_{V_n}$ and its submodule $M'_n$ are graded and that the relations $\Rel'_v$ are homogeneous of degree~$1$. We can therefore decompose the equation (\ref{Msurj2}) into its homogeneous parts. In particular we may assume that $f'$ is homogeneous of some degree~$d$. Then for any solution of (\ref{Msurj2}) modulo $(t^{-i})$, replacing $m$ and $g_v$ by their homogeneous parts of degree~$d$, respectively $d-1$, yields another solution modulo $(t^{-i})$. In the inductive process above, we can therefore always arrange that $m_i$ and $g_{v,i}$ are homogeneous of degree~$d$, respectively $d-1$. Then they all lie in fixed finite dimensional $\BF_q$-subspaces of $M'_n$ and~$R_{V_n}$. This then guarantees that $m:=\sum_{j=0}^\infty m_j\otimes t^{-j}$ lies in $M'_n\otimes_{\BF_q}\BF_q[[t^{-1}]]$ and $g_v:=\sum_{j=0}^\infty g_{v,j}\otimes t^{-j}$ lies in $R_{V_n}\otimes_{\BF_q}\BF_q[[t^{-1}]]$. In other words, we have found a solution of (\ref{Msurj2}) with coefficients in $\BF_q[[t^{-1}]]$, as desired.
\end{Proof}

%%%%%%%%%%%%%%%%%%%%%%%%%%%%%%%%%%%%%

\begin{Lem}\label{fkLem}
The elements $\pi(\tilde f_1),\ldots,\pi(\tilde f_r)\in R_n$ are algebraically independent over~$F$.
\end{Lem}

\begin{Proof}
For any $1\le k\le r$ and $1\le\nu\le n$ consider the element
\UseTheoremCounterForNextEquation
\begin{equation}\label{AfknuDef}
f'_{k,\nu}
\ := \sum_{w \in V'_{k,\nu}} \frac{1}{X_{k,\nu} + w}
\ \in \ R_{V_n}.
\end{equation}
If $\nu>1$, we have $V_1\subset V'_{k,\nu}$, and then the relations $\Rel_{X_{k,\nu} + w}$ from (\ref{RelDef}) show that 
\begin{eqnarray*}
\qquad\qquad f'_{k,\nu}\otimes1
&=& \kern-7pt\sum_{w \in V'_{k,\nu} \mymod V_1} \; \sum_{v'\in V_1} \frac{1}{X_{k,\nu}+w+v'}\otimes1 \\
&\equiv& \kern-7pt\sum_{w \in V'_{k,\nu} \mymod V_1} \frac{1}{tX_{k,\nu}+tw}\otimes t \qquad\qquad\qquad \hbox{modulo $J_n$}\\
&=& \sum_{w' \in V'_{k,\nu-1}} \frac{1}{X_{k,\nu-1}+w'}\otimes t \\
&=& f'_{k,\nu-1}\otimes t
\end{eqnarray*}
and hence $\pi(f'_{k,\nu}\otimes1) = \pi(f'_{k,\nu-1}\otimes t)$. By induction on $\nu$ it follows that $\pi(f'_{k,\nu}\otimes1) = \pi(f'_{k,1}\otimes t^{\nu-1})$ for all $1\le\nu\le n$. 
In particular we have $\pi(\tilde f_k) = \pi(f'_{k,n}\otimes1) = \pi(f'_{k,1}\otimes t^{n-1})$. It therefore suffices to show that the elements $\pi(f'_{1,1}\otimes1),\ldots,\pi(f'_{r,1}\otimes1)$ are algebraically independent over~$F$.

For this observe that (\ref{NN'R}) yields a commutative diagram 
$$\xymatrix{
RS_n & R_n \ar[l] \ar@{}[r]|-{\textstyle\ni} & \pi(f'_{k,1}{\otimes}1) \\
\llap{$RS_{V_1}\otimes_{\BF_q}F\, =\ $}
RS_1 \ar[u] & \ R_1\ \ar[l] \ar[u] 
\ar@{}[r]|-{\raisebox{5pt}{$=\tilde R_1\ni$}} & \ f'_{k,1}{\otimes}1\ \ar@{|->}[u] \\}$$
Here the lower horizontal arrow is injective by Construction \ref{RVCons}, and the left vertical arrow is injective by Remark \ref{NN'Rem}. Also, the elements $f'_{1,1}\otimes1, \ldots, f'_{r,1}\otimes1$ of $\tilde R_1 = R_{V_1}\otimes_{\BF_q}F$ are algebraically independent over~$F$ by Theorem \ref{RVBasis} (a). Together this implies that the images of $\pi(f'_{1,1}\otimes1), \ldots, \pi(f'_{r,1}\otimes1)$ in $RS_n$ are algebraically independent over~$F$. They are therefore themselves algebraically independent over~$F$, as desired.
\end{Proof}

%%%%%%%%%%%%%%%%%%%%%%%%%%%%%%%%%%%%%

\begin{Lem}\label{MnUInvariants}
The submodule $\tilde M_n\subset\tilde R_n$ is stable under $U$ and its $U$-invariants are
$$\tilde M_n^U\ =\ F[\tilde f_1,\ldots,\tilde f_r].$$
\end{Lem}

\begin{Proof}
For any $1\le k\le r$ the set $\tilde E_k$ from (\ref{tEkDef}) is obtained from the set $\tilde\Delta_k$ on replacing the element $1\in\tilde\Delta_k$ by the element $\frac{1}{X_{k,n}} = \tilde f_k-\sum_{1\not=e_k\in\tilde\Delta_k}e_k$. By combining this fact with the decomposition $F[\tilde f_k] = F\oplus F[\tilde f_k]\cdot \tilde f_k$ we find that 
\begin{eqnarray*}
\bigoplus_{e_k\in\tilde\Delta_k}F[\tilde f_k] \cdot e_k
&=& F\oplus F[\tilde f_k] \cdot f_k \oplus
\bigoplus_{1\not=e_k\in\tilde\Delta_k}\kern-5pt F[\tilde f_k] \cdot e_k \\
&=& F\oplus\bigoplus_{e_k\in\tilde E_k}F[\tilde f_k] \cdot e_k.
\end{eqnarray*}
Taking the tensor product over $1\le k\le r$ and setting $\tilde E_I := \prod_{k\in I}\tilde E_k$, we deduce that
%\begin{eqnarray*}
%\tilde M_n &=& \bigoplus_{e\in\tilde\Delta_1\cdots\tilde\Delta_r} \kern-10pt
%F[\tilde f_1,\ldots,\tilde f_r]\cdot e \\
%&=& \bigoplus_{I\subset\{1,\ldots,r\}}\;
%\bigoplus_{e\in\tilde E_I} F[\tilde f_k|_{k\in I}]\cdot e.
%\end{eqnarray*}
\UseTheoremCounterForNextEquation
\begin{equation}\label{MeDecomp}
\tilde M_n\ = \bigoplus_{e\in\tilde\Delta_1\cdots\tilde\Delta_r} \kern-10pt
F[\tilde f_1,\ldots,\tilde f_r]\cdot e 
\ = \bigoplus_{I\subset\{1,\ldots,r\}}\;
\bigoplus_{e\in\tilde E_I} F[\tilde f_k|_{k\in I}]\cdot e.
\end{equation}
Since $U$ permutes each $\tilde E_k$ and fixes each~$\tilde f_k$, this shows that $U$ acts on~$\tilde M_n$. Also, since $U$ acts transitively on $\tilde E_1\cdots\tilde E_r$, it also acts transitively on $\tilde E_I$ for each subset $I\subset\{1,\ldots,r\}$. The fact that $\tilde f_k=\sum_{e\in\tilde E_k}e$ and the above description of $\tilde M_n$ therefore implies that 
\begin{eqnarray*}
\tilde M_n^U \!&=& \kern-5pt\bigoplus_{I\subset\{1,\ldots,r\}}\kern-10pt
F[\tilde f_k|_{k\in I}]\cdot \sum_{e\in\tilde E_I} e \\
&=& \kern-5pt\bigoplus_{I\subset\{1,\ldots,r\}}\kern-10pt
F[\tilde f_k|_{k\in I}]\cdot \prod_{k\in I} \tilde f_k
\ \ =\ \ F[\tilde f_1,\ldots,\tilde f_r],
\end{eqnarray*}
as desired.
\end{Proof}

%%%%%%%%%%%%%%%%%%%%%%%%%%%%%%%%%%%%%

\begin{Lem}\label{MAIsom}
The projection $\pi$ induces an isomorphism $\tilde M_n\isoto R_n$.
\end{Lem}

\begin{Proof}
By Lemma \ref{MAsurj} the induced map $\pi\colon \tilde M_n\to R_n$ is already surjective. It remains to show that it is injective, or equivalently that its kernel $\tilde M_n\cap J_n$ is zero.
For this observe that $\tilde M_n$ is stable under $U$ by Lemma \ref{MnUInvariants}, and recall that $J_n$ is stable under $U$ by construction. Thus the intersection $\tilde M_n\cap J_n$ is also stable under~$U$. Since $U$ is a finite group of $q$-power order acting on an $\BF_q$-vector space, we have $\tilde M_n\cap J_n=0$ if and only if $\tilde M_n^U\cap J_n = (\tilde M_n\cap J_n)^U=0$. We are therefore reduced to showing that $\pi|\tilde M_n^U$ is injective.
But by Lemma \ref{MnUInvariants} we have $\tilde M_n^U = F[\tilde f_1,\ldots,\tilde f_r]$, and Lemma \ref{fkLem} implies that $\pi|F[\tilde f_1,\ldots,\tilde f_r]$ is injective.
\end{Proof}

%%%%%%%%%%%%%%%%%%%%%%%%%%%%%%%%%%%%%

\medskip 
To simplify notation, we now set $f_k:=\pi(\tilde f_k)$ and $\Delta_k :=\pi(\tilde\Delta_k)$ and $E_k := \pi(\tilde E_k)$ for all $1\le k\le r$ and abbreviate $E_I := \prod_{k\in I}E_k$ for any subset $I\subset\{1,\ldots,r\}$.
%In view of Lemma \ref{fkLem} we give a name to the free $F$-subalgebra
%\UseTheoremCounterForNextEquation
%\begin{equation}\label{R'nDef}
%R'_n\ :=\ F[f_1,\ldots,f_r]\ \subset\ R_n.
%\end{equation}
Then we can summarize the results of this subsection as follows:

\begin{Thm}\label{AAll1}
\begin{itemize}
\item[(a)] The elements $f_1,\ldots,f_r\in R_n$ are algebraically independent over~$F$.
\item[(b)] Each $\Delta_k$ and $E_k$ is a subset of $R_n$ of cardinality $q^{r(n-1)+k-1}$.
%$|\circVlimits_{k,n}^{\prime}| = q^{r(n-1)+k-1}$.
\item[(c)] The product induces bijective maps
$$\xymatrix@R-20pt@C-13pt{
\Delta_1\times\ldots\times\Delta_r \ 
\ar[rr]^-\sim && \ \Delta_1\cdots\Delta_r 
\ar@{}[r]|-{\subset} & R_n \rlap{\hbox{\qquad and}}\\
E_1\times\ldots\times E_r \ 
\ar[rr]^-\sim && \ E_1\cdots E_r 
\ar@{}[r]|-{\subset} & R_n \rlap{.} \\}$$
\item[(d)] The ring $R_n$ is a free module over $F[f_1,\ldots,f_r]$ with basis $\Delta_1\cdots\Delta_r$.
\item[(e)] We have $R_n^U = F[f_1,\ldots,f_r]$.
\item[(f)] As an $F$-vector space $R_n$ decomposes as
$$R_n\ = \bigoplus_{I\subset\{1,\ldots,r\}}\;
\bigoplus_{e\in E_I} F[f_k|_{k\in I}]\cdot e.$$
\end{itemize}
\end{Thm}

\begin{Proof}
Assertion (a) is the content of Lemma \ref{fkLem}, and (d) and the parts of (b) and (c) concerning the subsets $\Delta_k$ follow from Lemma \ref{MAIsom}. 
Part (e) follows directly from Lemmas \ref{MnUInvariants} and \ref{MAIsom}. 
The remaining assertions follow by combining Lemma \ref{MAIsom} with the decomposition (\ref{MeDecomp}).
\end{Proof}

%%%%%%%%%%%%%%%%%%%%%%%%%%%%%%%%%%%%%%%%%%%%%%%%%%%%%%%%%%%%%%%%%%%%%%

\subsection{Consequences}
\label{Specialn3}

\begin{Prop}\label{RLocFIsRS}
We have $RS_n = R_n[f_1^{-1},\ldots,f_r^{-1}]$. 
\end{Prop}

\begin{Proof}
By the construction (\ref{RAVCons}) of $R_n$ and $RS_n$ we have $RS_n=RS_{V_n}\otimes_{R_{V_n}}R_n$. 
By Proposition \ref{LocalizeByFk} applied to the elements $X_{k,\nu}\in V_n$ in the order (\ref{VarOrder}), the ring $RS_{V_n}$ is the localization of $R_{V_n}$ obtained by inverting the elements $f'_{k,\nu}$ from (\ref{AfknuDef}) for all $k$ and~$\nu$. Thus $RS_n$ is the localization of $R_n$ obtained by inverting the elements $\pi(f'_{k,\nu}\otimes1)$ for all $k$ and~$\nu$.
But by downward induction on~$\nu$, the equation $\pi(f'_{k,\nu}\otimes1) = \pi(f'_{k,\nu-1}\otimes t)$ from the proof of Lemma \ref{fkLem} implies that
$$\pi(f'_{k,\nu}\otimes1)
%\ =\ \pi(f'_{k,1}\otimes t^{\nu-1})
\ =\ \pi(f'_{k,n}\otimes t^{\nu-n})
\ =\ \pi(\tilde f_k)\cdot t^{\nu-n}
\ =\ f_k\cdot t^{\nu-n}.$$
for all $k$ and~$\nu$. Thus $RS_n$ is the localization of $R_n$ obtained by inverting the elements $f_k$ for all~$k$, as desired.
\end{Proof}

\begin{Thm}\label{InjCor}
The ring $R_n$ is an integral domain and injects into~$RS_n$.
\end{Thm}

\begin{Proof}
Theorem \ref{AAll1} (d) implies that $R_n$ injects into $R_n[f_1^{-1},\ldots,f_r^{-1}]$; hence by Proposition \ref{RLocFIsRS} it injects into~$RS_n$. Since $RS_n$ is an integral domain by Corollary \ref{RSVIntDom}, everything follows.
\end{Proof}

%%%%%%%%%%%%%%%%%%%%%%%%%%%%%%%%%%%%%

\begin{Thm}\label{RnCM}
The ring $R_n$ is Cohen-Macaulay.
\end{Thm}

\begin{Proof}
By Theorem \ref{AAll1} the ring $R_n$ is free of finite rank over the polynomial ring $F[f_1, \ldots, f_r]$. Thus $R_n$ has Krull dimension $r$ and the elements $f_1, \ldots, f_r$ form a regular sequence in~$R_n$ of length~$r$. The same then follows for the localization of $R_n$ at the irrelevant maximal ideal $\bigoplus_{d>0} R_{n,d}$; hence this localization is Cohen-Macaulay. Using \cite[Cor.$\;$2.2.15]{Bruns-Herzog}
% \cite[Exc. 2.1.27.c]{Bruns-Herzog}
it follows that the graded ring $R_n$ itself is Cohen-Macaulay.
\end{Proof}

%%%%%%%%%%%%%%%%%%%%%%%%%%%%%%%%%%%%%

\begin{Thm}\label{RnInvts}
For any $1\le n'\le n$ the natural homomorphism $R_{n'}\to R_n$ from Construction \ref{NN'Cons} induces an isomorphism from $R_{n'}$ to the subring of invariants in $R_n$ under the kernel of the natural surjection $\GL_r(\BF_q[t]/(t^n)) \onto \GL_r(\BF_q[t]/(t^{n'}))$.
\end{Thm}

\begin{Proof}
By induction on $n'$ it suffices to consider the case $n'=n-1$ with $n\ge2$. Let $H$ denote the subgroup of $\GL_r(\BF_q[t]/(t^n))$ in question. Then $H$ consists of all elements $h\in U$ such that $h(X_{k,n})\in X_{k,n}+V_1$ for all $1\le k\le r$. In fact, any independent choice of an element of $X_{k,n}+V_1$ for each $k$ corresponds to a unique element of~$H$. Thus we may write $H=H_1\times\ldots\times H_r$, where each $H_k$ permutes the coset $X_{k,n}+V_1$ simply transitively and fixes $X_{k',n}$ for all $k'\not=k$. Each $H_k$ then also acts trivially on $E_{k'}$ for all $k'\not=k$. It follows that for any subset $I\subset\{1,\ldots,r\}$, the orbits of $H$ on $E_I := \prod_{k\in I}E_k$ are simply the products over $k\in I$ of the orbits of $H_k$ on~$E_k$.
By the relation (\ref{RelDef}), the sum over the $H_k$-orbit of an element $[\frac{1}{X_{k,n}+w}\otimes1] \in E_k$ comes out as
$$\sum_{v'\in V_1}\Bigl[\frac{1}{X_{k,n}+w-v'}\otimes 1\Bigr]
\ =\ \Bigl[\frac{1}{tX_{k,n}+tw}\otimes t\Bigr]
\ =\ \Bigl[\frac{1}{X_{k,n-1}+tw}\otimes 1\Bigr]\cdot t.$$
Here $w$ runs through $V'_{k,n}$; hence the element $tw$ runs through $tV'_{k,n}=V'_{k,n-1}$.
Thus with 
$$\tilde E'_k\ :=\ \biggl\{ \frac{1}{X_{k,n-1} + w}\otimes1 \,\biggm|\, w \in V_{k,n-1}^{\prime} \biggr\}$$
and $E'_k := \pi(\tilde E'_k)$ and $E'_I := \prod_{k\in I}E'_k$, the description of $R_n$ in Theorem \ref{AAll1} (f) implies that 
\UseTheoremCounterForNextEquation
\begin{equation}\label{RnInvts1}
R_n^H\ = \bigoplus_{I\subset\{1,\ldots,r\}}\;
\bigoplus_{e'\in E'_I} F[f_i|_{i\in I}]\cdot e'.
\end{equation}
This coincides with the description of $R_{n-1}$ from Theorem \ref{AAll1} (f) for $n-1$ in place of~$n$. Indeed, we already have $\tilde E'_k\subset\tilde R_{n-1}\subset\tilde R_n$, and the subset $E'_k\subset R_n$ is precisely the image of the corresponding subset of $R_{n-1}$ under the natural homomorphism $R_{n-1}\to R_n$. Moreover, the equation $\pi(f'_{k,n}\otimes1) = \pi(f'_{k,n-1}\otimes t)$ from the proof of Lemma \ref{fkLem} implies that the elements $f_k\in R_n$ and $f_k\in R_{n-1}$ differ only by a factor of~$t$. Thus (\ref{RnInvts1}) implies that the natural homomorphism $R_{n-1}\to R_n$ induces an isomorphism $R_{n-1}\isoto R_n^H$, as desired.
\end{Proof}

%%%%%%%%%%%%%%%%%%%%%%%%%%%%%%%%%%%%%%%%%%%%%%%%%%%%%%%%%%%%%%%%%%%%%%

\subsection{The ideal of the boundary}
\label{SpecialABoundaryIdeal}

Consider the ideal $I_n := I_{A,V^r_{t^n}} \subset R_n$ from (\ref{AIVNDef}). Our first result is entirely analogous to Theorem \ref{IVBasis}.

\begin{Thm}\label{AInBasis}
\begin{enumerate}
\item[(a)] The ideal $I_n$ is a free module over $R_n^U$ with basis $E_1\cdots E_r$.
\item[(b)] The ideal $I_n$ is a free module over the group ring $F[U]$.
\end{enumerate}
\end{Thm}

\begin{Proof}
Theorem \ref{AAll1} implies that $U$ acts freely transitively on $E_1\cdots E_r$ and that $R_n^U=F[f_1,\ldots,f_r]$. Theorem \ref{AAll1} (f) thus shows that $E_1\cdots E_r$ is the basis of a free $R_n^U$-submodule of~$R_n$. Denoting this submodule by~$M$, it follows that $M$ is a free module over $F[U]$. It remains to show that $M=I_n$.

For this note first that for any element $e_1\cdots e_r\in E_1\cdots E_r$, the reciprocals $e_1^{-1},\ldots,e_r^{-1}$ form a basis of $V_n$ over $\BF_q[t]/(t^n)$. Thus for any $1\le s<r$ and any $i\colon V^s_{t^n}\into V^r_{t^n}=V_n$ as above, at least one of $e_1^{-1},\ldots,e_r^{-1}$ lies in $V^r_{t^n}\setminus i(V^s_{t^n})$. By the description of $\pi_i$ in Proposition \ref{i*HomoA} it follows that $e_1\cdots e_r \in \Ker(\pi_i)$. Varying $s$ and $i$ this shows that $e_1\cdots e_r \in I_n$, and varying $e_1\cdots e_r$ then implies that $M\subset I_n$.

Next observe that $0\to M\to I_n\to I_n/M\to0$ is a short exact sequence of $F[U]$-modules. Since $M$ is a free $F[U]$-module, taking $U$-invariants yields a short exact sequence $0\to M^U\to I_n^U\to (I_n/M)^U\to H^1(U,M)=0$. Also, since $U$ is a finite group of $q$-power order acting on the $F$-vector space $I_n/M$ with $\BF_q\subset F$, we have $I_n/M=0$ if and only if $(I_n/M)^U=0$. To prove that $M=I_n$, by the short exact sequence it is therefore enough to prove that $M^U=I_n^U$.

As the given basis $E_1\cdots E_r$ of $M$ over $R_n^U$ is a single free orbit under~$U$, the submodule $M^U$ is the free $R_n^U$-module generated by the element $\sum E_1\cdots E_r = f_1\cdots f_r$. In other words it is the principal ideal of $R_n^U$ generated by $f_1\cdots f_r$. Since $R_n^U = F[f_1,\ldots,f_r]$ with algebraically independent $f_1,\ldots,f_r$, this ideal is the intersection of the ideals $R_n^U\cdot f_k$ for all $1\le k\le r$. Thus it suffices to prove that $I_n^U\subset R_n^U\cdot f_k$ for every fixed $1\le k\le r$.

To achieve this suppose first that $r=1$. Then $I_n$ is the augmentation ideal of $R_n$ by construction (\ref{AIVNDef}); hence $I_n^U$ is the augmentation ideal of $R_n^U=F[f_1]$ and thus equal to $R_n^U\cdot f_1$, and we are done.

Otherwise consider the permutation $\sigma := (k,k+1,\ldots,r)$ and let $i\colon V^{r-1}_{t^n}\into V^r_{t^n}$ denote the $\BF_q[t]/(t^n)$-linear embedding defined by $i(X_{j,n})=X_{\sigma j,n}$ for all $1\le j\le r-1$. Then the image of $i$ is the $\BF_q[t]/(t^n)$-submodule generated by $X_{k',n}$ for all $1\le k'\le r$ with $k'\not=k$. For all $1\le k'\le r$, the definition of $f_{k'}$ shows that 
$$\pi_i(f_{k'}) \ = \sum_{w \in V'_{k',n}} \pi_i\bigl(\bigl[\tfrac{1}{X_{k',n} + w}\otimes 1\bigr]\bigr).$$
In the case $k'=k$ we have $X_{k,n} + w \not\in i(V^{r-1}_{t^n})$ for all $w \in V'_{k,n}$ and therefore $\pi_i(f_k)=0$. 
Otherwise we have $k'=\sigma j$ for some $1\le j\le r-1$ and the description of $\pi_i$ in Proposition \ref{i*HomoA} implies that
$$\pi_i(f_{k'}) \ = 
\kern-5pt\sum_{w'\in i^{-1}(V'_{k',n})} 
\kern-10pt \bigl[\tfrac{1}{X_{j,n} + w'}\otimes 1\bigr].$$
But the definition of $V'_{k',n}$ implies that $i^{-1}(V'_{k',n})$ is precisely the subspace $V'_{j,n} \subset V^{r-1}_{t^n}$ for $r-1$ in place of~$r$. Thus the elements $\pi_i(f_{k'})\in R_{A,V^{r-1}_{t^n}}$ for all $k'\not=k$ are precisely the elements $f_1,\ldots,f_{r-1} \in R_{A,V^{r-1}_{t^n}}$. As they are algebraically independent over~$F$ by Theorem \ref{AAll1} (a), it follows that the kernel of $\pi_i|F[f_1,\ldots,f_r]$ is the ideal $F[f_1,\ldots,f_r]\cdot f_k$. 
By construction $I_n^U = I_n\cap R_n^U$ is contained in this kernel; hence we are done.
\end{Proof}

%%%%%%%%%%%%%%%%%%%%%%%%%%%%%%%%%%%%%

\medskip
Next consider the integral closure $R_n^\norm := R_{A,V^r_{t^n}}^\norm$ from Construction \ref{RAVNormalize} and its ideal $I_n^\norm := I_{A,V^r_{t^n}}^\norm$ from (\ref{AIVNnormDef}).

\begin{Thm}\label{InInNorm}
The homomorphism $R_n\to R_n^\norm$ induces an isomorphism $I_n\isoto I_n^\norm$.
\end{Thm}

\begin{Proof}
Theorem \ref{InjCor} implies that $R_n$ injects into $R_n^\norm$; hence $I_n$ injects into $I_n^\norm$. We identify it with its image and must then show that $I_n=I_n^\norm$.

For this observe that $0\to I_n\to I_n^\norm\to I_n^\norm/I_n\to0$ is a short exact sequence of $F[U]$-modules. Since $I_n$ is a free $F[U]$-module by Theorem \ref{AInBasis} (b), taking $U$-invariants yields a short exact sequence $0\to I_n^U\to (I_n^\norm)^U\to (I_n^\norm/I_n)^U\to H^1(U,I_n)=0$. Also, since $U$ is a finite group of $q$-power order acting on the $F$-vector space $I_n^\norm/I_n$ with $\BF_q\subset F$, we have $I_n^\norm/I_n=0$ if and only if $(I_n^\norm/I_n)^U=0$. To prove that $I_n=I_n^\norm$, by the short exact sequence it is therefore enough to prove that $I_n^U=(I_n^\norm)^U$.

We will determine both sides of this equation separately. In the proof of Theorem \ref{AInBasis} we have already seen that $I_n^U = M^U = R_n^U\cdot f_1\cdots f_r$. Also, Proposition \ref{RLocFIsRS} shows that the ideal $R_n\cdot f_1\cdots f_r\subset R_n$ defines a closed subscheme of $Q_{A,V_n}$ whose support is precisely the boundary. This implies that $I_n = \sqrt{ R_n\cdot f_1\cdots f_r}$ within~$R_n$. By (\ref{AIVNnormDef}) we therefore have
$$I_n^\norm\ =\ \sqrt{R_n^\norm\cdot f_1\cdots f_r}$$
within~$R_n^\norm$. This in turn implies that 
$$(I_n^\norm)^U\ =\ \sqrt{(R_n^\norm)^U\cdot f_1\cdots f_r}$$
within~$(R_n^\norm)^U$. But the construction of $R_n^\norm$ implies that $(R_n^\norm)^U$ is the integral closure of $R_n^U$ in~$RS_n^U$. Since $R_n^U$ is already a regular integral domain by Theorem \ref{AAll1}, we deduce that $(R_n^\norm)^U = R_n^U$. As $R_n^U\cdot f_1\cdots f_r$ is already a reduced ideal in $R_n^U = F[f_1,\ldots,f_r]$, it follows that $(I_n^\norm)^U = R_n^U\cdot f_1\cdots f_r = I_n^U$, as desired.
\end{Proof}

\begin{Rem}\label{InInRem}
\rm Theorem \ref{InInNorm} says that the ideals of the boundary in the graded coordinate rings of $Q_{A,V_n}$ and $Q_{A,V_n}^\norm = \OM^r_{A,t^n}$ are the same. Thus, in a sense $Q_{A,V_n}$ and $\OM^r_{A,t^n}$ differ only in the reduced subschemes at the boundary in that some points are identified. More precisely I expect that $Q_{A,V_n}$ is the quotient of $\smash{\OM^r_{A,t^n}}$ by the resulting equivalence relation on the underlying topological space. Compare Remark \ref{NormNotIsomRem}.
\end{Rem}

%%%%%%%%%%%%%%%%%%%%%%%%%%%%%%%%%%%%%%%%%%%%%%%%%%%%%%%%%%%%%%%%%%%%%%

\subsection{Modular forms and cusp forms}
\label{Specialn5}

For any integer $d$ let $I_{n,d}\subset R_{n,d}$ denote the homogeneous parts of degree $d$ of $I_n\subset R_n$.

\begin{Thm}\label{DimForm1}
For any integer $d\ge1$ we have
$$\dim_F(R_{n,d}) \ = 
\sum_{\emptyset\not=I\subset\{1,\ldots,r\}}\!
{\textstyle\binom{d-1}{|I|-1}}\cdot\prod_{k\in I} q^{r(n-1)+k-1}.$$
\end{Thm}

\begin{Proof}
By Theorem \ref{AAll1} (a), for any subset $I\subset\{1,\ldots,r\}$ the ring $F[f_k|_{k\in I}]$ is isomorphic to a polynomial ring in $|I|$ variables over~$F$. Since each $f_k$ is homogeneous of degree~$1$ and $d\ge1$, the homogeneous part of $F[f_k|_{k\in I}]$ of degree $d-|I|$ is zero if $I=\emptyset$ and has dimension $\binom{d-1}{|I|-1}$ otherwise. Also recall that each element of $E_k$ is homogeneous of degree~$1$. The formula follows by combining all this with Theorem \ref{AAll1} (b) and (c) and (f).
\end{Proof}

\begin{Rem}\label{DimForm1Rem}
\rm Since $R_{n,d}\into R_{n,d}^\norm$ is not an isomorphism in general, Theorems \ref{DimForm1} and \ref{SatakeProjection} together do not yet give us a dimension formula for spaces of modular forms. For cusp forms, on the other hand, we succeed using Theorem \ref{InInNorm}:
\end{Rem}

%%%%%%%%%%%%%%%%%%%%%%%%%%%%%%%%%%%%%

\begin{Thm}\label{DimForm2}
For any $d\ge1$ the space $I_{n,d}$ is free module of rank $\binom{d-1}{r-1}$ over the group ring $F[U]$.
%For any integer $d\ge1$, the dimension over $F$ of the homogeneous part of degree $d$ of $I_n$ is
%$$|U|\cdot \binom{d-1}{r-1}.$$
\end{Thm}

\begin{Proof}
By Theorem \ref{AAll1} the ring $R_n^U = F[f_1,\ldots,f_r]$ is isomorphic to a polynomial ring in $r$ variables over~$F$. Since each $f_k$ is homogeneous of degree~$1$ and $d\ge1$, the homogeneous part of $R_n^U$ of degree $d-r$ has dimension $\binom{d-1}{r-1}$. Also recall that each element of $E_k$ is homogeneous of degree~$1$ and that $U$ acts freely transitively on $E_1\cdots E_r$ by Theorem \ref{AAll1} (b) and (c). Thus the formula follows from Theorem \ref{AInBasis}.
\end{Proof}

\begin{Thm}\label{DimForm3}
For any $d\ge1$ the space of cusp forms $\Gamma(\OM^r_{A,t^n},\CI(d))$ is a free module of rank $\binom{d-1}{r-1}$ over the group ring $F[U]$.
\end{Thm}

\begin{Proof}
Combine Theorems \ref{SatakeProjection} and \ref{ModCuspFormsIsom} and \ref{InInNorm} and \ref{DimForm2}.
\end{Proof}

%%%%%%%%%%%%%%%%%%%%%%%%%%%%%%%%%%%%%

\medskip
Finally consider any subgroup $U'<U$. Then $U'\backslash\OM^r_{A,N}$ is the Satake compactification of a Drinfeld moduli space of some intermediate level associated to~$U'$, and $\Gamma(\OM^r_{A,t^n},\CI(d))^{U'}$ is the space of cusp forms of weight $d$ of that level. From Theorem \ref{DimForm3} we directly deduce:

\begin{Thm}\label{DimForm4}
For any $d\ge1$ we have
$$\dim_F\Gamma(\OM^r_{A,t^n},\CI(d))^{U'}
\ =\ \textstyle [U:U']\cdot\binom{d-1}{r-1}.$$
\end{Thm}

%%%%%%%%%%%%%%%%%%%%%%%%%%%%%%%%%%%%%

We can also turn this into a dimension formula for analytic cusp forms, as follows. Take the natural homomorphism $\kappa\colon \SL_r(\BF_q[t]) \to \GL_r(\BF_q[t]/(t^n))$ and consider the arithmetic subgroups $\Gamma(t^n) := \ker(\kappa)$ and $\Gamma_1(t) := \kappa^{-1}(U)$ of $\SL_r(\BF_q[t])$. Consider an arbitrary subgroup $\Gamma(t^n) < \Gamma < \Gamma_1(t)$. Let $\Cinf$ denote the completion of an algebraic closure of the field $\BF_q((t^{-1}))$. Let $\CS_d(\Gamma)$ denote the space of analytic cusp forms of rank $r$ and weight $d$ and level $\Gamma$ according to \cite[Def.\;6.1]{BBP1}.

\begin{Thm}\label{DimForm5}
For any $d\ge1$ we have
$$\dim_\Cinf \CS_d(\Gamma)
\ =\ \textstyle [\Gamma_1(t):\Gamma]\cdot\binom{d-1}{r-1}.$$
\end{Thm}

\begin{Proof}
We must identify the analytic cusp forms with algebraic cusp forms as in \cite{BBP2}. For this let $\hat A\cong\smash{\prod_{\Fp}A_\Fp}$ denote the profinite completion of~$A$. Let $K(t^n)$ be the kernel of the natural homomorphism $\hat\kappa\colon \GL_r(\hat A) \onto \GL_r(\BF_q[t]/(t^n))$ and set $K := \Gamma\cdot K(t^n)$. 
Then $\hat\kappa$ induces isomorphisms $\Gamma/\Gamma(t^n) \isoto K/K(t^n) \isoto U'$ for a certain subgroup $U'$ of the group $U<\GL_r(\BF_q[t]/(t^n))$ from above. In the case $\Gamma=\Gamma_1(t)$ this subgroup is $U\cap\SL_r(\BF_q[t]/(t^n))$, whose index in $U$ is $q^{n-1}$. In the general case we therefore have
$$[U:U']\ =\ q^{n-1} \cdot [\Gamma_1(t):\Gamma].$$
Next, the fact that $\hat\kappa(K(t^n)) < \hat\kappa(K) < U$ shows that $K(t^n)$ and $K$ are fine open compact subgroups in the sense of \cite[Def.\;1.4]{PinkSatake}. The associated Drinfeld moduli spaces from \cite{PinkSatake} are $M^r_{A,K(t^n)} = M^r_{A,t^n}$ and $M^r_{A,K} = U'\backslash M^r_{A,t^n}$. Moreover $\det(K) = \det(K(t^n))$ is the kernel of the natural homomorphism $\hat A^\times \onto (\BF_q[t]/(t^n))^\times$. Using \cite[Prop.\;8.7]{BBP2} this implies that $M^r_{A,K(t^n)}$ and $M^r_{A,K}$ have the same constant field $F_K$ which is a finite Galois extension of $F$ of degree $[F_K/F] = |\BF_q^\times\backslash \hat A^\times/\det(K)| = q^{n-1}$.

The main point is that \cite[Thm.\;10.9]{BBP2} yields an isomorphism between algebraic and analytic modular forms
$$\Gamma(\OM^r_{A,t^n},\CO(d))\otimes_{F_K}\Cinf\ \stackrel{\sim}{\longto}\ \CM_d(\Gamma(t^n)).$$
This induces an isomorphism between algebraic and analytic cusp forms
$$\Gamma(\OM^r_{A,t^n},\CI(d))\otimes_{F_K}\Cinf\ \stackrel{\sim}{\longto}\ \CS_d(\Gamma(t^n)).$$
By taking $U'$-invariants we obtain an isomorphism
$$\Gamma(\OM^r_{A,t^n},\CI(d))^{U'}\otimes_{F_K}\Cinf\ \stackrel{\sim}{\longto}\ \CS_d(\Gamma(t^n))^{U'}
\ =\ \CS_d(\Gamma).$$
Using Theorem \ref{DimForm4} we deduce that 
\begin{eqnarray*}
\dim_\Cinf\CS_d(\Gamma)
&=& \dim_{F_K} \Gamma(\OM^r_{A,t^n},\CI(d))^{U'} \\
&=& [F_K/F]^{-1}\cdot \dim_F \Gamma(\OM^r_{A,t^n},\CI(d))^{U'} \\
&=& q^{1-n}\cdot \textstyle [U:U']\cdot\binom{d-1}{r-1} \\
&=& \textstyle [\Gamma_1(t):\Gamma]\cdot\binom{d-1}{r-1},
\end{eqnarray*}
as desired.
\end{Proof}

\newpage

%%%%%%%%%%%%%%%%%%%%%%%%%%%%%%%%%%%%%%%%%%%%%%%%%%%%%%%%%%%%%%%%%%%%%%
%%%%%%%%%%%%%%%%%%%%%%%%%%%%%%%%%%%%%%%%%%%%%%%%%%%%%%%%%%%%%%%%%%%%%%

%%%%%%%%%%%%%%%%%%%%%%%%%%%%%%%%%%%%%%%%%%%%%%%%%%%%%%%%%%%%%%%%%%%%%%%%

\end{document}